\documentclass[12pt]{amsart}  
 
\usepackage{amssymb} 
\usepackage{enumerate, amsfonts, latexsym,psfig,epsfig, color} 

\usepackage{epstopdf}

\newtheorem{theorem}{Theorem}[section]

\newtheorem{prop}[theorem]{Proposition} 
\newtheorem{proposition}[theorem]{Proposition}

\newtheorem{lemma} [theorem]{Lemma} 
\newtheorem{corollary} [theorem] {Corollary} 

\newtheorem{notation}[theorem]{Notation} 
 
\theoremstyle{definition} 
\newtheorem{definition}[theorem]{Definition} 
\newtheorem{example}[theorem]{Example} 
 
\theoremstyle{remark} 
\newtheorem{remark}[theorem]{Remark} 
\newtheorem{RIproof}[theorem]{}%Proof of Regular Implosions}
 
\newtheorem{remarks}[theorem]{Remarks} 
 
\numberwithin{equation}{section} 
 
\newfont{\msb}{msbm10 scaled 1200} 
\newfont{\euf}{eufm10 scaled 1200} 
\def\Prf:       {\it Proof. \rm} 
 
\def\Gthree{\mathcal G_3}
\def\bt {\begin{theorem}} 
\def\et         {\hfill $\square$\end{theorem}} 
\def\bl {\begin{lemma}} 
\def\el         {\hfill $\square$\end{lemma}} 
\def\bd         {\begin{definition}} 
\def\ed{\end{definition}} 
\def\bc         {\begin{corollary}} 
\def\ec         {\hfill $\square$ \end{corollary}} 
\def\bea {\begin{eqnarray}} 
\def\eea        {\end{eqnarray}} 
\def\be {\begin{eqnarray*}} 
\def\ee {\end{eqnarray*}} 
\def\bp {\begin{proposition}} 
\def\ep {\end{proposition}}  
\def\bex {\begin{example}} 
\def\eex {\end{example}} 
\def\height{\text{\rm{time}}} 
\def\Pin{\Pi} 
  
\def \Cal{\mathcal} 
\def \Bbb{\mathbb} 
\def \S{\Sigma} 
\def\D{\Delta} 
\def\E{\cal E} 
\def\ssm{\smallsetminus} 
 
\def\area{\text{\rm{Area}}}

\def\<{\langle} 
\def\>{\rangle} 
\def\|{{\,|\! |\, }}

\def\P{\Cal P} 
\def\cal{\Cal} 
\def\e{\varepsilon}

\def\-{\underline}

\def\G{\Gamma} 
\def\A{\Cal A} 
 
\def\time{\text{\rm{time}}}

\def\iso{\cong} 
\def\supp{\text{\rm{Supp}}} 
\def\pre{\text{\rm{Pre}}} 
\def\lpl{left para-linear } 
\def\rpl{right para-linear } 
\def\n  {\,|\partial\Delta|} 
\def\vecZ{\mathcal Z} 
\def\tz{{T_0}} % A1A5Lemma 
\def\life{\text{\rm Life}} 
\def\ttt{T_1} 
\def\mess  {\text{\rm Mess}} 
\def\R{\Cal R} 
\def\AFourC     {2C_1 + 6\ll + 2B(5T_0 + 6\ttt + 2) + 2MC_4(6\ttt +
8T_0 + 3) + (B+3)(3\ttt + 2T_0)M +5M+2} %The length of sum |A4|

\def\T{\mathcal T} 
\def\cmm{C_{(\mu,\mu')}(2)} 
\def\vin        {\in_v}  % `is a virtual member of' 
\def\subT{\chi(\Pin_{\T})} 
\def\CT{\chi_c(\T)} %this used to be \chi(\T) 
\def\down{\text{\rm{down}}} 
\def\bonus{\text{\rm{bonus}}} 
\def\dlong{\text{\rm{D}}\Lambda} 
 
\def\plT{p_l({\mathcal T})} 
\def\tplT{\tilde p_l(\T)}
\def\Bb         {(B+3)(3\ttt + 2T_0)M + 6B\ttt + 4BT_0 + 2\ll + 2B +
5M + 1} % The total contribution of all bonus' (without the $n$).
\def\F{\mathcal F} 
\def\ptmm{\hat t_1(\mu,\mu')} 
\def\prmm{\hat\rho(\mu,\mu')} 
\def\pEmm{\hat{\text{\euf T}}  (\mu,\mu')} 

\def\pTmm{\hat\T(\mu,\mu')} 

\def\pT{\hat\T}

\def\xT{x({\mathcal T})} %first non-constant edge 
\def\prT{p_r({\mathcal T})}  
\def\tprT{\tilde p_r(\T)}
\def\ttwo{t_2(\T)} 
\def\ET{{\text{\euf T}}}
\def\Emm{{\text{\euf T}}(\mu,\mu')}
\def\tone{t_1(\T)} 
\def\I{\mathcal I} 
\def\r{\rho} 
\def\rT{\rho_\T} 
 
\def\bonusT{\text{\rm{bonus}}(\T)} 
\def\K  { 2C_0 + 2K_1 + 2B + 1} % |S_0| \leq Kn... 

\def\ll {\lambda_0} 
 
\def\L{\Cal L} 
\def\R{\Cal R}
\def\QT{Q(\T)}

\catcode`\@=11 
\def\serieslogo@{\relax} 
\def\@setcopyright{\relax} 
\catcode`\@=12

\begin{document}

\title[Mapping tori of free group automorphisms]
{The quadratic isoperimetric inequality
for mapping tori of free group automorphisms I:
Positive automorphisms\\
}
 
% author  information 
\author[Martin R. Bridson]{Martin R.~Bridson} 
\address{Martin R.~Bridson\\
Mathematics Department\\
Imperial College of Science, Technology and Medicine\\
180 Queen's Gate \\
London SW7 2BZ\\
UK} 
\email{m.bridson@ic.ac.uk}  
 
% author  information 
\author[Daniel Groves]{Daniel Groves}
\address{Daniel Groves\\
Merton College\\
Oxford OX1 4JD\\
UK} 
\email{grovesd@maths.ox.ac.uk} 
 
\date{29 July, 2003}
 
 
 
\thanks{The first author's work was supported in part
by an EPSRC Advanced Fellowship} 
 
\begin{abstract} {If $F$ is a finitely generated free group
and $\phi$ is a positive automorphism of $F$
then  $F\rtimes_\phi\mathbb Z$
satisfies a quadratic isoperimetric inequality.} 
\end{abstract}

\maketitle

Associated to an   automorphism $\phi$
of any group $G$ one has the algebraic {\em mapping torus}
$G\rtimes_\phi\mathbb Z$. In this paper we shall be concerned
with the case where $G$ is a finitely generated
free group, denoted $F$. We seek to understand the complexity
of  word problems for the
 groups $F\rtimes_\phi\mathbb Z$ as measured
by their Dehn functions. 
 
The class of groups of the form $F\rtimes_\phi\mathbb Z$
has been the subject of intensive investigation  
in recent years and a rich structure has begun to
emerge  in keeping with the subtlety of the
classification of free group automorphisms \cite{BFH}, \cite{BFH2}
\cite{BH2}, \cite{FH}, \cite{Lu}, \cite{Sela}.  (See \cite{BestICM} and
the references therein.) 
Bestvina--Feighn and Brinkmann proved that if $F\rtimes_\phi\mathbb Z$
doesn't contain
a free abelian subgroup of rank two then
it is hyperbolic \cite{BF}, \cite{Brink},
i.e. its Dehn function
is linear. Epstein and Thurston \cite{E+}
proved that if $\phi$ is induced by a surface automorphism
(in the sense discussed below) then $F\rtimes_\phi\mathbb Z$ is automatic
and hence 
has
a quadratic Dehn function. The question
of whether or not all non-hyperbolic groups of the form
$F\rtimes_\phi\mathbb Z$  have quadratic Dehn functions has
attracted a good deal of attention.

Recall that an automorphism $\phi$ of a finitely generated free
group $F$ is called {\em positive} if there is a 
basis $a_1,\dots,a_n$ for $F$ such that the reduced
word representing each $\phi(a_i)\in F$ contains no inverses $a_j^{-1}$.

\renewcommand{\thethmspec}{{\bf{Main Theorem. \kern-.3em}}}

\smallskip

\noindent
\begin{thethmspec} \label{MainThmPos}
{\em Let $F$ be a finitely generated free group.
If $\phi$ is a positive  automorphism of $F$, then
$F\rtimes_\phi\mathbb Z$ satisfies a  quadratic isoperimetric inequality.}
\end{thethmspec} 
 
\smallskip

Modulo a simple change in the interpretation of the symbols used, the
proof of this theorem extends {\em verbatim}  to 
automorphisms $\phi$ that have a  power that admits a train track
representative. Not all automorphisms of free groups admit such
representatives. Nevertheless, in a subsequent article
\cite{BGrovesII} we use the relative train-track technology developed
by Bestvina, Feighn and Handel (\cite{BH2} and \cite{BFH}) and the
architecture of the proof of the Main Theorem to establish the
quadratic isoperimetric inequality for all groups of the form $F
\rtimes \Bbb Z$.

Much of our modern understanding of the automorphisms of
free groups has been guided by the  analogy with
automorphisms of surface groups, i.e. mapping classes
of surfaces of finite type. This analogy  provides a
useful reference point when considering the word problems of
mapping tori. 
 
A self-homeomorphism of a compact surface $S$ defines an outer
automorphism of $\pi_1S$ and hence a semidirect product
$\pi_1S\rtimes_\phi\mathbb Z$. This group  is the fundamental
group of a compact   3-manifold, namely the
mapping torus $M_\phi$ of the homeomorphism. By using
Thurston's Geometrization
Theorem for Haken manifolds, Epstein and Thurston \cite{E+} were
able to prove that
$\pi_1S\rtimes_\phi\mathbb Z$
 is an automatic group and hence its Dehn function is
either linear or quadratic. If $S$ has boundary then
only the quadratic case arises. A more geometric explanation for the existence of 
a quadratic isoperimetric inequality in the bounded case
comes from the fact that
 $M_\phi$ supports a
metric of non-positive curvature, as does any irreducible
3-manifold with non-empty boundary \cite{mb-shs}, \cite{leeb}. 
 
If $S$ has boundary, then $\pi_1S$ is free. Thus the
foregoing considerations give many examples of free-by-cyclic
groups that have quadratic Dehn functions. But
there are many types of  free group automorphisms that
do not arise from surface automorphisms, for example
those $\phi$ that do not have a power leaving any
non-trivial conjugacy class invariant, and those $\phi$
for which there is a word $w\in F$ such that the function $n\mapsto |\phi^n(w)|$
grows like a super-linear polynomial. 
 
The non-automaticity of certain  $F\rtimes_\phi\mathbb Z$ provides
a more subtle obstruction to realising $\phi$ as a surface automorphism:
in contrast to the
Epstein-Thurston Theorem, Brady, Bridson and Reeves  \cite{BB},
\cite{BR} showed that certain mapping tori $F_3\rtimes \mathbb Z$ are
not automatic, for example that associated to the  automorphism
$[a\mapsto a,\, b\mapsto ab,\, c\mapsto a^2c]$.
Such examples show that one cannot proceed  via automaticity in order
to prove the Main Theorem. Nor can one rely on non-positive
curvature,  because Gersten \cite{Ge} showed that the above example
$F_3\rtimes\mathbb Z$ is not the fundamental group of any compact
non-positively curved space. Thus one needs a new approach to the
quadratic isoperimetric inequality.

A technique for dealing with classes of linearly growing automorphisms is described by
Brady and Bridson in \cite{BB}, and Macura  \cite{Mac} developed techniques for
dealing with polynomially growing automorphisms.  But these techniques apply only
to  restricted classes of automorphisms and do not speak to the core 
problem of establishing the quadratic isoperimetric inequality
for mapping tori of general free group automorphisms. 
In the present article and its sequel
we attack  this core problem
directly, undertaking a detailed analysis of the geometry of van Kampen diagrams over the
natural presentations of free-by-cyclic groups.

This paper is organised as follows.
In Section \ref{vanKampSection} we
recall some basic definitions associated to Dehn functions.
In Sections \ref{BCSection} and \ref{time}
 we record some simple but important observations
concerning the large-scale behaviour of  
the van Kampen
diagrams associated to free-by-cyclic groups and in
particular the geometry of {\em corridor} subdiagrams. (The automorphisms
considered up to this point are not assumed to be positive.)
These observations lead us to a strategy for proving the Main Theorem based
on the geometry of the {\em time flow of corridors}. In Section  \ref{StrategySection}
we state a sharper version of the Main Theorem adapted to this strategy
and reduce to the study of automorphisms with  stability
properties  that regulate the evolution of corridors. 
  In Section \ref{PrefFutSec} we develop the notion of
 {\em preferred future} which allows
 us to trace the trajectory of
  $1$-cells in the corridor flow.
  
The estimates that we establish in Sections 5 and 6 reduce us to the
nub of the difficulties that one faces in trying to prove the Main
Theorem, namely
the possible existence of  large blocks of ``constant letters". A
sketch of the strategy that we shall use to overcome this problem is presented
in Section 7. The three main ingredients in this strategy are
the elaborate global cancellation arguments in Section 8, the machinery
of {\em teams} developed in Section 9, and the {\em bonus scheme}
developed in Section 10 to accommodate a final tranche of cancellation
phenomena whose quirkiness eludes the grasp of teams. In a brief final
section we gather our many estimates to establish the bound required for
the Main Theorem. A glossary of constants is included for the reader's convenience. 

\section{Van Kampen Diagrams} \label{vanKampSection}

We recall some basic definitions and facts concerning 
Dehn functions and van Kampen diagrams.

\subsection{Dehn Functions and Isoperimetric Inequalities} 
 
Given a finitely presented group $G=\langle \mathcal A \mid 
\mathcal R \rangle$
and a word $w $ in the generators $\mathcal A^{\pm 1}$ that
represents $1\in G$,
one defines
$$
\area (w)  = \\
 \min\big\{ N \in {\mathbb N}^+ \; | \;
\exists\text{ equality }
w = \prod^N_{j=1}u_j^{-1}r_j u_j \text{ in $F(\mathcal A)$
with } r_j \in \mathcal R^{\pm 1} \big\}\, .
$$

The {\it Dehn function} $\delta(n)$ of the finite 
presentation $\langle \mathcal A \mid \mathcal R\rangle$ is defined by 
$$ 
\delta(n) \; = \; \max\{\text{\rm{Area}}(w) \; 
|\; w \in \text{\rm{ker}}(F(\mathcal A) \twoheadrightarrow G), 
\; |w| \leq n \, \} \, , 
$$ 
where $|w|$ denotes the length of the   word $w$. 
Whenever two presentations 
define isomorphic (or indeed  quasi-isometric) 
groups, the Dehn functions of 
the finite presentations   
are equivalent under the relation 
$\simeq$ that identifies functions  
$[0,\infty)\to [0,\infty)$ that only differ by a quasi-Lipschitz 
distortion of their domain and their range.

For any constants $p,q\ge 1$, one sees that 
$n\mapsto n^p$ is $\simeq$ equivalent to $n\mapsto n^q$ 
only if $p=q$. Thus it makes sense to say that the 
``Dehn function of a group" is $\simeq n^p$. 
 
A group $\G$ is said to {\em satisfy a quadratic isoperimetric 
inequality} if its Dehn function is $\simeq n$ or 
$\simeq n^2$. A result of  Gromov \cite{Gromov}, detailed 
proofs of which were given by several authors, states that if  
a  Dehn function is subquadratic, then it is linear --- 
see \cite[III.H]{BH} for a discussion, proof and references.

See \cite{steer} for a thorough and  elementary account of 
what is known about  Dehn functions and an  
explanation of their connection  with filling 
problems in Riemannian geometry.  
\smallskip  
 
\subsection{Van Kampen diagrams}\label{vkD} 
 
According to van Kampen's lemma (see \cite{vK}, 
\cite{LS} or \cite[I.8A]{BH})   
an equality $w = \prod^N_{j=1}u_jr_ju_j^{-1}$ in the 
free group $\mathcal A$, with $N=\area(w)$, 
can be portrayed by  
a finite, 1-connected, 
combinatorial 2-complex with basepoint, embedded in $\mathbb R^2$. Such a complex is 
called a {\em van Kampen diagram} for $w$; its oriented 1-cells   
are labelled by elements of $\mathcal A^{\pm 1}$; 
the boundary label on each 2-cell (read with clockwise 
orientation from one of its vertices) is an element  
of $\mathcal R^{\pm 1}$;  and the boundary cycle of the 
complex (read with positive orientation from the basepoint) 
is the word $w$; 
the number 
of 2-cells in the  diagram   is $N$.  Conversely, any van Kampen diagram with $M$ 
2-cells gives rise 
to an equality in $F(\mathcal A)$ expressing the word 
labelling the boundary cycle 
of the diagram as a product of $M$  
conjugates of the defining relations.  
Thus 
$\text{\rm{Area}}(w)$ is the minimum number of 2-cells among all 
van Kampen diagrams  
for $w$. If a van Kampen diagram $\Delta$ for $w$ has $\text{\rm{Area}}(w)$ 
2-cells, then $\Delta$ is a called a {\em least-area} diagram. If 
the underlying 2-complex is homeomorphic to a 2-dimensional 
disc, then the van Kampen diagram is called a {\em disc diagram}. 
 
We use the term {\em area} to describe the number of 2-cells in a 
van Kampen diagram, and write $\text{\rm{Area }} \Delta$. We write 
$\partial \Delta$ to denote the boundary cycle of the diagram; we write 
$|\partial\Delta|$ to denote the length of this cycle.

Note that associated to a van Kampen diagram $\Delta$ with basepoint $p$ 
one has a morphism of 
labelled, oriented graphs $h_\Delta: (\Delta^{(1)},p)\to (\mathcal C_\A, 1)$,  where 
$\mathcal C_\A$ is the Cayley graph associated to the choice of 
generators $\A$ for $G$. The map $h_\Delta$ takes $p$ to the 
identity vertex  $1\in \mathcal C_\A$ and preserves the labels on oriented edges. 
 
We shall need the following simple observations. 
 
\begin{lemma} If a van Kampen diagram 
$\Delta$ is least-area, then every simply-connected 
subdiagram of $\Delta$ is also least-area. 
\end{lemma} 
 
Recall that a function $f:\mathbb N\to [0,\infty)$ 
is {\em sub-additive} if $f(n+m)\le f(n) + f(m)$ 
for all $n,m\in\mathbb N$. For example, given $r\ge 1,\, k>0$, 
the function $n\mapsto kn^r$ is sub-additive. 
 
\begin{lemma}  \label{disc} 
Let $f:\mathbb N\to [0,\infty)$ be a sub-additive function and let $\P$ 
be a finite presentation of a group. 
If  $\text{\rm{Area }} \Delta\le f(|\partial\Delta|)$ for every 
least-area disc diagram $\Delta$ over $\P$, then the Dehn 
function of $\P$ is $\le f(n)$. 
\end{lemma}

\subsection{Presenting $F\rtimes\mathbb Z$} 
 
We shall establish the quadratic bound required for 
the Main Theorem by examing the nature of van Kampen 
diagrams over the following natural (aspherical) presentations of
free-by-cyclic groups. 
 
Given a finitely generated free group $F$ and 
an automorphism $\phi$ of $F$, we fix a basis 
$a_1,\dots,a_m$ for $F$, write $u_i$ to denote 
the reduced word equal to $\phi(a_i)$ in $F$, and 
present $ 
F\rtimes_\phi\mathbb Z$ by 
\begin{equation}\label{presentation} 
\P\iso \langle a_1,\dots,a_m,t\mid 
t^{-1}a_1tu_1^{-1},\dots, t^{-1}a_mtu_m^{-1}\rangle. 
\end{equation} 
We shall work exclusively with this presentation. 
\medskip

\begin{figure}[htbp] 
\begin{center} 
  
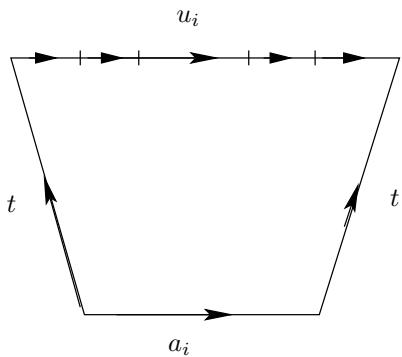 
  
\caption{A $2$-cell in a van Kampen diagram for $F \rtimes_{\phi} \mathbb Z$.} 
\label{figure:2cell} 
\end{center} 
\end{figure}

\subsection{Time and $t$-Corridors with naive tops} 
 
The use of $t$-corridors as a tool for investigating 
van Kampen diagrams has 
become well-established in recent years. In the 
setting of van Kampen diagrams over the above presentation,  
$t$-corridors are easily described. 
 
Consider a van Kampen diagram $\Delta$ over 
the above presentation $\P$ and focus on an edge in the boundary 
$\partial \Delta$ that is labelled 
$t^{\pm 1}$ (read with positive orientation from the basepoint).  
If this edge lies in the boundary 
of a 2-cell, then the boundary  cycle of this 2-cell has 
the form $t^{-1}a_itu_i^{-1}$ (read with suitable orientation from 
a suitable point, see Figure \ref{figure:2cell}). In particular, there 
is  exactly one other edge  
in the boundary of the 2-cell that is labelled $t$; crossing 
this edge we enter another 2-cell with a similar boundary 
label, and iterating the argument we get a chain of 2-cells 
running across the diagram; this chain terminates at an edge of 
$\partial \Delta$ which (following the orientation of $\partial \Delta$ 
in the direction of our original edge labelled $t^{\pm 1}$) is labelled 
$t^{\mp 1}$. This chain of 2-cells is called a {\it{$t$-corridor}}. 
The edges labelled $t$ that we crossed in the above description 
are called the {\em vertical} edges of the corridor.  
The vertical edge on $\partial \Delta$ labelled $t^{-1}$ is 
called the {\em initial} end of the corridor, and at the other end one 
has the {\em terminal} edge. 
 
Formally, one should define a $t$-corridor to be a combinatorial map 
to $\Delta$ from a suitable subdivision of $[0,1]\times [0,1]$: the 
initial edge is the restriction of this map to $\{0\}\times [0,1]$; the 
vertical edges are the images of the 1-cells of the form $\{s\}\times [0,1]$, 
oriented so that the edge joining $(s,0)$ to $(s,1)$ is labelled $t$.
The {\em naive top} of the corridor is the edge-path obtained by restricting 
the above map to $[0,1]\times\{1\}$, and the {\em bottom} is the restriction 
to  $[0,1]\times\{0\}$. 
\smallskip 
 
\noindent{\bf Left/Right Terminology:} The orientation of a disc 
diagram induces an orientation on its corridors. Whenever 
we focus on an individual corridor, we shall regard its 
 initial as being {\em left}most  and its terminal 
edge as being {\em right}most. (This is just a suggestive way of saying 
 that the 
corridor map from $[0,1]\times (0,1)\subset \mathbb R^2$ to $\Delta\subset 
\mathbb R^2$ is 
orientation-preserving.) 
 
\smallskip

\medskip 
 
\begin{figure}[htbp] 
\begin{center} 
  
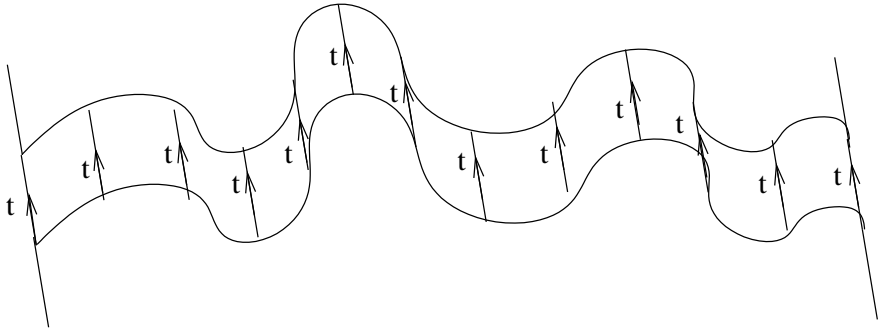 
  
\caption{A $t$-corridor} 
\label{figure:corridor} 
\end{center} 
\end{figure}

See \cite{BG} for a  detailed account of $t$-corridors. 
Here we shall need only the following easy facts: 
\begin{enumerate} 
\item 
distinct $t$-corridors 
have disjoint interiors;  
\item 
if $\sigma$ is the edge-path in $\Delta$ running along 
the (naive) top or bottom  of a $t$-corridor, then $\sigma$ is 
 labelled 
by a word in the letters $\mathcal A^{\pm 1}$ 
that is equal in $F\rtimes\mathbb Z$   to the 
words labelling the subarcs of $\partial \Delta$ 
which share the  endpoints of $\sigma$ (given appropriate 
orientations); 
\item if we are in a least-area diagram  
then the word on the bottom of the corridor is freely reduced; 
\item the number of 2-cells in the 
$t$-corridor is the length of the word labelling the 
bottom side. 
\item In subsection 1.2 we described the map $h_\Delta$ 
associated to a van Kampen diagram. This map 
sends vertices of $\Delta$ to vertices of the Cayley graph $\mathcal C_\A$, 
i.e. elements of $F\rtimes \langle t \rangle$.  If the initial vertex of a  
directed edge in $\D$ is sent to an element of the form $wt^j$, with 
$w\in F$, then the edge is defined to occur at {\bf time} $j$. Note that the 
vertical edges of a fixed corridor all occur at the same time. 
\end{enumerate}     
 
We will
consider the {\em dynamics} 
of the automorphism $\phi$.

\begin{definition} 
[Time and Length]  Item (5) above 
implies that the time of each 
$t$-corridor $S$ is well-defined;
we denote it $\height(S)$.  
 
 We define the {\em length} of a  
corridor $S$ to be the number of 2-cells that it 
contains, which is equal to the number of 1-cells along its bottom. 
We write $|S|$ to denote the length of $S$. 
\end{definition}

\subsection{Conditioning the Diagram} 
 
We are working with the following presentation of $F\rtimes_\phi\mathbb Z$ 
$$ 
\mathcal P = \langle a_1,\dots,a_m,t\mid 
t^{-1}a_1tu_1^{-1},\dots, t^{-1}a_mtu_m^{-1}\rangle. 
$$ 
 
In the light of Lemma \ref{disc}, in order to prove the main 
theorem it suffices to consider only {\em disc diagrams}. Therefore, 
henceforth we shall assume that all diagrams are topological discs. 
We shall also assume that all of the discs considered are 
{\em least-area} diagrams for freely reduced words.  
 
\begin{lemma} 
Every least-area disc diagram  over $\P$ is  
the union of its $t$-corridors. 
\end{lemma} 
 
\begin{proof} Since the diagram is a disc, every 1-cell lies in 
the boundary of some 2-cell. The boundary of each 2-cell 
contains two edges labelled $t$. Consider the equivalence relation 
on 2-cells generated by $e\sim e'$ if the boundaries of $e$ and $e'$ share an 
edge labelled $t$. Each equivalence class forms either a $t$-corridor 
or else a $t$-ring, i.e. the closure of an annular sub-diagram 
whose internal and external cycles are labelled by a word in the 
generators of $F$. If the latter case arose, then since 
$F$ is a free group, the word $u$ on the external 
cycle would be freely equal to the empty word (since it contains no edges 
 labelled $t$). This would contradict the hypothesis that the diagram 
  is least-area, because one could reduce its area 
by excising the simply-connected sub-diagram bounded by this cycle, 
replacing it  with the zero-area diagram for $u$ over the free 
presentation of $F$. 
\end{proof}

\subsection{Folded Corridors}

In the light of the above lemma, we see that the diagrams $\Delta$ that we 
need to consider are essentially determined once one knows which 
pairs of boundary edges are connected by $t$-corridors. However, there 
remains a slight ambiguity arising from the fact that free-reduction in 
the free group is not a canonical process (e.g.  $x = (xx^{-1})x = x (x^{-1}x)$).

\begin{figure}[htbp] 
\begin{center} 
  
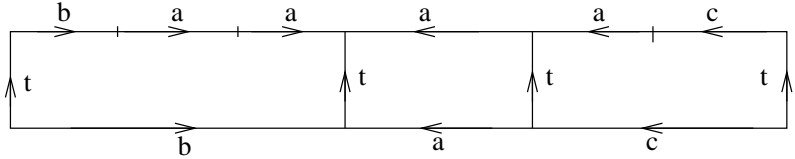 
  
\caption{An unfolded corridor} 
\label{figure:Fold1} 
\end{center} 
\end{figure}

\begin{figure}[htbp] 
\begin{center} 
  
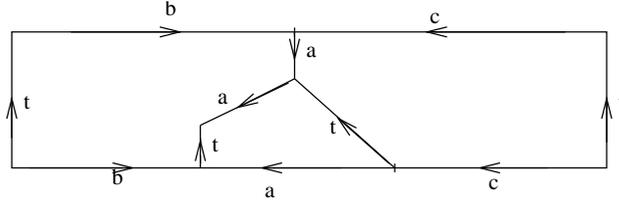 
 
\caption{The  corresponding  unfolded corridor.} 
\label{figure:Fold2} 
\end{center} 
\end{figure}

To avoid this ambiguity, we  fix a least area disc diagram $\Delta$  
and assume that its corridors are {\em folded} in the sense of \cite{B-plms}. 
The topological closure  $T\subset\Delta$ of each corridor is a combinatorial 
disc. 
The hypothesis ``least area" alone 
forces the label on the {\em bottom} of the corridor 
to be a {\em freely reduced} word in the letters $a_i^{\pm 1}$. 
We define the 
 {\em top} of the (folded) corridor to be the  
 injective edge-path that remains when one deletes from the 
frontier of $T$ 
 the bottom and ends of the corridor. The word labelling 
 this path is the freely reduced word in $F$ that equals the 
 label on the naive top of the corridor. Note that, unlike the 
bottom of the corridor, the top may fail to intersect  the closure of some 
2-cells --- see Figures \ref{figure:Fold1} and \ref{figure:Fold2} (where the automorphism is $a \mapsto a, 
b \mapsto ba^2, c \mapsto ca$). 
 
\begin{notation} 
We write $\top (S)$ and $\bot (S)$, respectively, to denote the top and bottom 
of a folded corridor $S$. 
\end{notation} 
 
\smallskip 
 
{\centerline{ 
{\em Henceforth we shall refer to folded $t$-corridors simply as ``corridors".}}}

\bigskip

\medskip 
 
\subsection{Naive Expansion and Death} 
 
For each generator $a_i\in F$ we have the reduced word 
$u_i=\phi(a_i)$. Given a reduced  word $v=a_{i(1)}\dots a_{i(m)}$ 
we define the {\em naive expansion} of $\phi(v)$ to be  
the (unreduced) concatenation $u_{i(1)}\dots u_{i(m)}$. 
 
Note that if $v$ is the label on an interval of the bottom of a corridor, 
then the naive expansion of $\phi(v)$ is the label on the 
corresponding arc of the naive top of the corridor.

An edge $\e$ on the bottom of a  corridor $S$ is said to {\em die} in $S$ 
if the 2-cell containing that edge  does not contain any edge of  
$\top(S)$.  (Equivalently, if $w$ is the label on  $\bot(S)$ and $a_i$ is 
the label on $\e$, then the subword  $u_i=\phi(a_i)$ in 
the naive expansion of $\phi(w)$ is cancelled 
completely during the free reduction encoded in $\Delta$.) In Figure \ref{figure:Fold2}
the edge labelled $a$ on the bottom of the corridor dies.

\section{Singularities and Bounded Cancellation} \label{BCSection} 
 
We have noted that the structure  of a (folded, least-area disc) diagram 
over the natural presentation of a free-by-cyclic group 
is  the union of its 
corridors. 
In this section we pursue an 
understanding of how these corridors meet. 
\smallskip

\medskip

\begin{figure}[htbp] 
\begin{center} 
  
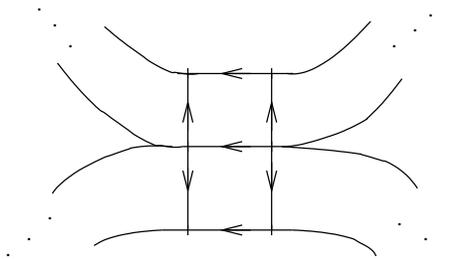 
  
\caption{Corridors cannot meet this way in a least-area diagram} 
\label{figure:wrongway} 
\end{center} 
\end{figure}

The first observation to make is that corridors cannot meet as in 
Figure \ref{figure:wrongway}. 
 
\begin{lemma} If $S\neq S'$, then $\bot(S)\cap\bot(S')$ 
consists of at most one point. 
\end{lemma} 
 
\begin{proof} For each letter $a$, there is only one type 
of 2-cell which has the label $a$ on its bottom side.  Thus, if two corridors 
were to meet in the manner of  Figure \ref{figure:wrongway}, then we would have a pair 
of 2-cells whose union was bounded by a loop labelled  
$u_it^{-1}tu_i^{-1}t^{-1}t$, which is 
freely equal to the identity. By excising this pair of 2-cells and 
filling the loop with a diagram of zero area, we would 
reduce the area of $\Delta$ without altering its boundary label --- 
but  $\Delta$  is assumed to be a least-area diagram. 
 
Thus $\bot(S)\cap\bot(S')$ contains no edges. To see that it cannot 
contain more than one vertex, follow the proof of Proposition 
\ref{SingularityProp}(1). 
\end{proof}

\begin{definition}  
A {\em singularity} in $\Delta$ is a non-empty connected  component of the intersection 
of the tops of two 
distinct folded corridors. A 2-cell  is said to {\em hit} the 
singularity if 
it contains an edge of the singularity.  

The singularity   is said to be degenerate if it consists of a single point, and
otherwise it is {\em non-degenerate}.
\end{definition}  
 
\medskip 
 
\begin{figure}[htbp] 
\begin{center} 
  
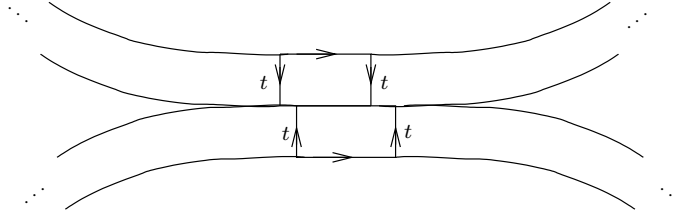 
  
\caption{A `singularity'} 
\label{figure:Singularity} 
\end{center} 
\end{figure}

Let $M$ be the maximum of the lengths of the words $u_i$ in 
our fixed presentation $\mathcal P$ of $F\rtimes_\phi \mathbb Z$. 
 
\begin{proposition}[Bounded singularities] \label{SingularityProp}$\ $

\begin{enumerate} 
\item[1.] If the tops of two corridors in a  least-area 
diagram meet, then their intersection is a singularity. 
\item[2.] 
There exists a constant $B$ depending only 
on $\phi$  such that less than $B$ 2-cells 
hit each singularity in a  least-area diagram over $\P$. 
\item[3.]  
If $\Delta$ is a least-area diagram over $\P$, 
then there are less than $2|\partial \Delta|$ non-degenerate singularities 
in $\Delta$, and each has length at most $MB$.
\end{enumerate} 
\end{proposition}

\begin{proof} Suppose that the intersection of the tops of two corridors $S$ and $S'$ 
contains two distinct vertices, $p$ and $q$ say. Consider the unique subarcs 
of $\top(S)$ and $\top(S')$ connecting $p$ to $q$. 
 Each of these arcs is labelled by a reduced word in 
the generators of $F$; since the arcs have the same endpoints in $\Delta$, 
these words must be identical.  If the arcs did not coincide, then 
we could excise the subdiagram that they bounded and replace it with 
a zero-area diagram, contradicting our least-area hypothesis. This proves 
(1).

\begin{figure}[htbp] 
\begin{center} 
  
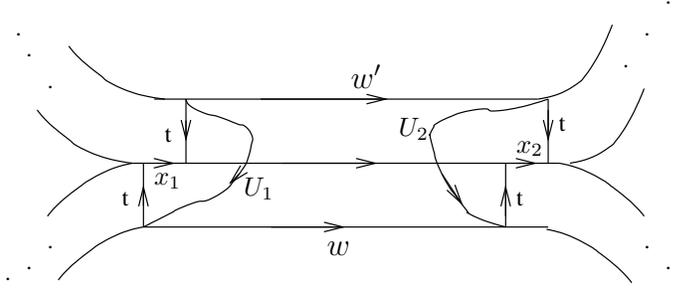 
  
\caption{The proof of Proposition \ref{SingularityProp}} 
\label{figure:BoundedSing} 
\end{center} 
\end{figure}

Figure \ref{figure:BoundedSing} portrays the argument we use to prove (2). In $S$ 
(respectively $S'$), we choose  
an outermost pair of oriented edges $\e_1, \e_2$ (resp. 
$\e_1',\e_2'$) labelled $t$ whose termini lie on the  
singularity. We then connect their endpoints by shortest 
arcs in the singularity as shown. Note that 
each of the arcs labelled $x_1$ and $x_2$ is contained in the top 
of a single 2-cell, and hence has length at most $M$. 
We write $\alpha_i$ to denote the concatenation of $\e_i$, the arc labelled $x_i$ 
and the inverse of $\e_i'$.

Let $U^{-1}_i\in F$ be the reduced word representing 
$\phi^{-1}(x_i)$. In $F\rtimes_\phi\mathbb Z$ we have $tx_it^{-1}U_i=1$; 
let  $\Delta_i$ be a least-area van Kampen diagram portraying  
this equality.   
 
Let $w$ (resp. $w'$) be the label on the edge-path 
in $\bot (S)$ (resp. $\bot(S')$) that connects 
the initial point of $\e_1$ (resp. $\e_1'$) to 
the initial point of $\e_2$ (resp. $\e_2'$).

If we excise from $\Delta$ the subdiagram bounded by the loop whose label 
is 
$t^{-1}wtx_2t^{-1}{w'}^{-1}tx_1^{-1}$, then we reduce the area of  
$\Delta$ by $|w| + |w'|$. (Recall that the edges on the bottom 
of a corridor are in 1-1 correspondence with the 2-cells of the 
corridor.) We may then attach a copy of $\Delta_i$ along $\alpha_i$ 
and fill the resulting loop labelled $U_1wU_2^{-1}{w'}^{-1}$ with 
a diagram of zero area, because this word is equal to $1$ 
in the free group $F$.  
Thus we obtain a new van Kampen diagram whose boundary label 
is the same as that of $\Delta$ and which has area 
$$ 
\area(\Delta) + \area(\Delta_1) + \area(\Delta_2) - |w| - |w'|. 
$$ 
Since $\Delta$ 
is assumed to be least-area, this implies that  
$ \area(\Delta_1) + \area(\Delta_2) \ge  |w| + |w'|.$

 Let $B_0$ be an upper bound on the area of 
all least-area van Kampen diagrams portraying equalities of the  
form  $txt^{-1}\phi^{-1}(x)^{-1}=1$ with $|x|\le M$. 
(It suffices to take $B_0=MM_{inv}$, where $M_{inv}$ is the maximum 
of the lengths of the reduced words $\phi^{-1}(a_i)$.) By definition,  
$ \area(\Delta_1) + \area(\Delta_2)\le 2B_0$, and hence 
$|w| + |w'|\le 2B_0$. Thus for (2) it suffices to let $B=2B_0 + 1$. 
 
The length of the singularity in the above argument 
is less than the sum of the lengths of  the naive 
expansions of $\phi(w)$ and $\phi(w')$. Since  $|w|+|w'|\le B$,  
the singularity has length less than $MB$.

It remains to bound the number of non-degenerate
singularities in $\Delta$. To this
end, we consider the subcomplex  $\Gamma\subset\Delta$ formed by the union of
the tops of all folded corridors. Arguing as in (1), we see that the
graph $\Gamma$ contains no non-trivial loops, i.e. it is a forest. Let 
$V$ denote the set of vertices in $\Gamma$ that have valence at least
3 or else lie on $\partial \Delta$. (Thus $V$ is the set of 
degenerate singularities, endpoints of non-degenerate singularities,
and endpoints of the tops of corridors.)
Let $E$ be the set of  connected components of
$\Gamma\smallsetminus V$.

$|V|-|E|$ is the number $\pi_0$  of connected components of the forest $\Gamma$.
The valence 1 vertices  $V^1\subset\Gamma$ are a subset of the endpoints of
the tops of corridors, so there are less than $|\partial\Delta|$
of them. One can calculate $|E|$ as half the sum of the valences of
the vertices $v\in V$, so $3(|V|-|V^1|) +|V^1| \le 2|E|$. 
Hence
$$
|E| = |V| - \pi_0 \le \frac 2 3 \big{(}|E| + |V^1|\big{)} -\pi_0
< \frac 2 3 \big{(}|E| + |\partial\Delta|\big{)}.
$$
Therefore $|E| < 2|\partial\Delta|$. 

Each non-degenerate singularity determines an element of $E$, so
the (crude) estimate in (3) is established. 
\end{proof}

\begin{lemma}[Bounded Cancellation Lemma] \label{BCL} There is a constant $B$, 
depending only on $\phi$, such that if 
$I$ is an interval consisting of $|I|$ edges 
on the bottom of a (folded) corridor $S$ in a least-area diagram over $\P$, 
and every edge of $I$ dies in $S$, then $|I| < B$. 
\end{lemma} 
\smallskip

\medskip

\begin{figure}[htbp] 
\begin{center} 
  
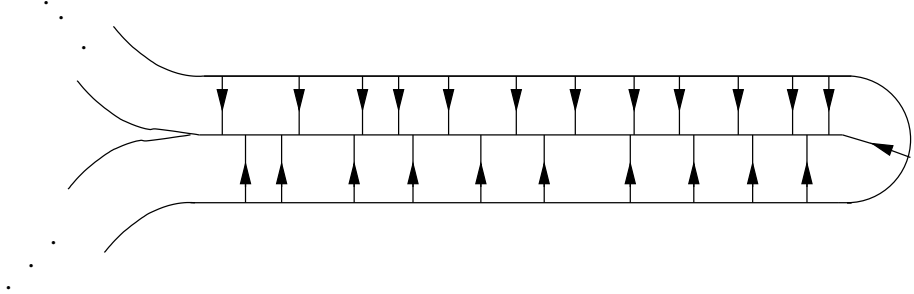 
  
\caption{Bounded cancellation lemma} 
\label{figure:BCL} 
\end{center} 
\end{figure} 
 
\begin{proof} The argument is entirely similar to that given for part (2) 
of the previous proposition. 
\end{proof}

The above lemma is a reformulation of the 
 Bounded Cancellation Lemma from \cite{Cooper}, 
which Cooper attributes to Thurston.

\begin{remark} {\em `Singularities are only 1 pixel large.'} 
The reader may find it useful to keep 
in mind the following picture: think of   
a least-area van Kampen diagram rendered on a computer 
screen and assume that the length of the boundary of  
the diagram is  large, so large that the constant $B$ 
in Proposition \ref{SingularityProp} has to be scaled to something less 
than 1 pixel in order to fit the picture on to the  computer's 
screen. 
In the resulting 
image one sees blocks of $t$-corridors as shown in Figure \ref{figure:singflow} 
below, and the singularities take on the appearance of classical $k$-prong
 singularities in the time-flow of $t$-corridors. 
\end{remark}

\medskip

\begin{figure}[htbp] 
\begin{center} 
  
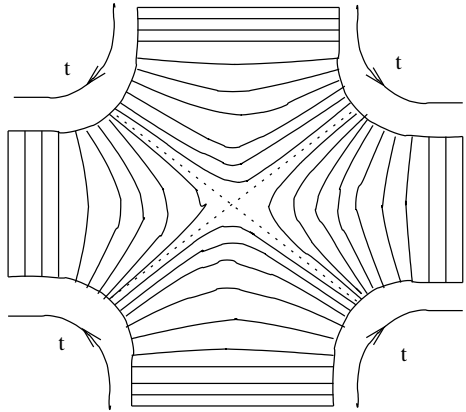 
  
\caption{Schematic depiction of a singularity} 
\label{figure:singflow} 
\end{center} 
\end{figure} 
 
\medskip

\section{Past, Future and Colour} \label{time} 
 
Our investigations thus far have led us to regard van Kampen diagrams 
over $\P$ as flows of corridors  
(at least schematically). We require some more vocabulary to pursue 
this approach. 
 
We continue to work with a fixed disc diagram $\Delta$ over $\P$.

\begin{definition}[Ancestors and Colour]\label{defMu} 
 Each edge $\varepsilon_1$ on the bottom of a corridor either 
lies in the boundary of $\Delta$, or else lies in the top of 
a unique 2-cell, the bottom of which we denote $\e_0$. We consider the 
partial ordering on the set $\mathcal E$ of edges from the bottom of all corridors 
generated by setting $\varepsilon_0 < \varepsilon_1$ whenever edges are related 
in this way.

If $\e'<\e$ then we 
call $\e'$ an {\em ancestor} of $\e$. The {\em past} of $\e$ 
is the set of its ancestors, and the {\em future} of $\e$ is  
the set of edges $\e''$ such that $\e<\e''$. 
 
Two edges are defined to be of the same {\em colour} if  
they have a common ancestor. Since every edge has 
a unique ancestor on the boundary, colours are in 
bijection with a subset\footnote{namely,  
those edges of $\partial\Delta$ that lie on the bottom of 
some 2-cell} of the edges in $\partial\Delta$ whose 
label is not $t$; in particular there are less than 
$|\partial\Delta|$ colours.

Each 2-cell in $\Delta$ has a unique edge 
in the bottom of a corridor. Thus 
we may also regard $\le$ as a partial 
ordering on the 2-cells 
of $\Delta$ and define the past, future and colour 
of a 2-cell. 
 
We define the past (resp. future) of a {\em corridor} 
 to be the union 
of the pasts (resp. futures) of its closed 2-cells. 
\end{definition} 
 
\begin{remark}\label{tree} 
Each $e\in\E$ and each 2-cell has at most one immediate 
ancestor (i.e. one that is maximal among its ancestors). 
Consider the graph $\mathcal F$ with vertex set $\E$ that has an edge 
connecting a pair of vertices if and only if 
one is the immediate ancestor of the other. Note 
that $\mathcal F$ is a forest 
(union of trees). 
 
The {\em colours} in the diagram correspond to the 
 connected components 
(trees) of this forest.  
 
 There is a natural embedding 
of $\mathcal F\hookrightarrow\Delta$: choose a point (`centre') 
 in the interior of each 2-cell 
and connect it to the centre of its immediate ancestor by an 
arc that passes through their common edge. 
 
\end{remark} 
 
If the future of a corridor $S'$ intersects a corridor $S$ then 
the intersection is connected: 
 
\begin{lemma}[Connected Pasts] \label{Connected}  
If a pair of 2-cells $\alpha$ and $\beta$ in a 
corridor $S$  have ancestors $\alpha'$ and $\beta'$ in a corridor $S'$, then every 
$2$-cell $\gamma$ that lies between $\alpha$ and $\beta$ in $S$  has 
an ancestor $\gamma'$ that lies between 
 $\alpha'$ and $\beta'$  in $S'$. 
 \end{lemma}

\begin{figure}[htbp] 
\begin{center} 
  
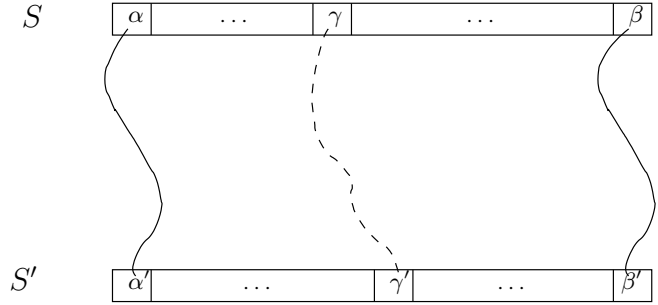 
  
\caption{The `loop' picture} 
\label{figure:loop} 
\end{center} 
\end{figure} 
 
\medskip 
 
\begin{proof} 
Connect the centres of  $\alpha$ and $\beta$ by an arc in the 
interior of $S$ that intersects only those 2-cells lying 
between $\alpha$ and $\beta$, and connect the centres of $\alpha'$ 
and $\beta'$ by 
a similar arc in the interior of $S'$. Along with these two 
arcs, we consider the embedded arcs 
connecting $\alpha$ to $\alpha'$ and $\beta$ to $\beta'$ in the forest 
$\mathcal F$ 
described in Remark \ref{tree}. 
These four arcs together form a  
loop, and the disc that this loop encloses does not intersect 
the boundary of $\Delta$. (Recall that $\Delta$ is a disc.) 
 
Consider the tree from $\mathcal F$ that contains $\gamma$. 
 We may assume that 
the arc in this tree 
that connects $\gamma$ to its ancestor on the boundary does 
not intersect the arc we chose in $S$. It must therefore intersect  
our loop either 
in $S'$, yielding the desired ancestor  $\gamma'$ in $S'$, or 
else in one of the arcs connecting $\alpha$ to $\alpha'$, or 
$\beta$ to $\beta'$. If the latter alternative pertains, $\alpha'$ or 
$\beta'$ is an ancestor of $\gamma$, and we are done. 
\end{proof} 
 
We highlight the degenerate case where the 2-cells $\alpha'$ and $\beta'$ 
are equal and have their bottom   on $\partial\Delta$: 
 
\begin{corollary}\label{muConn} 
Within a corridor, the 2-cells of each colour form a connected region.  
\end{corollary}

\section{Strategy, Strata and Conditioning}  
\label{StrategySection}

Everything that has been said up to  this point has 
 been true for  mapping tori of arbitrary automorphisms of 
 finitely generated free groups. {\em Henceforth, 
  we assume that the automorphism $\phi$ is positive.} 
 
A van Kampen diagram whose boundary cycle  has length $n$ contains at 
most $n/2$ corridors. Thus our Main Theorem is an immediate consequence of: 
 
\begin{theorem}\label{BoundS} 
 There is a constant $K$ depending only on $\phi$ 
such that each  corridor in a least-area diagram $\Delta$ over $\P$ 
has length at most $K\,|\partial\Delta|$. 
\end{theorem}

In order to establish the desired bound on the 
length of corridors, we must analyse how 
corridors grow as they flow into the future, and 
assess what cancellation can take place to inhibit this 
growth. In the remainder of this section we shall 
condition the automorphism to simplify 
the discussion of growth. 
 
\begin{remark}  
 The mapping torus $F\rtimes_{\phi^k}\mathbb Z$ is isomorphic to a  
 subgroup of finite index in  $F\rtimes_{\phi}\mathbb Z$, namely 
  $F\rtimes_{\phi}k\mathbb Z$. Thus, since the Dehn functions of 
  commensurable 
  groups are  
  $\simeq$ equivalent, we are free to replace $\phi$ by a convenient 
  positive power in our proof of the Main Theorem. 
  \end{remark}

\subsection{Strata} 
 
 In the following discussion we shall write $x$ to denote an 
arbitrary choice of letter from our basis $\{a_1,\dots,a_m\}$ 
for $F$.

Naturally associated to any positive automorphism one has 
{\em supports} and 
 {\em strata}.  The support 
  $\supp(x)$ associated to $x$ is  
   the set of all letters which appear in the freely 
   reduced word $\phi^j(x)$ for some $j \geq 0$.  
   The stratum $\S(x)\subset\supp(x)$ associated to $x$ consists 
   of those $y\in\supp(x)$ such that $\supp(x)=\supp(y)$. 
   
   Note that $y\in \supp(x)$ implies $\supp(y)\subseteq\supp(x)$, 
   and $y\in\S(x)$ implies $\S(y)=\S(x)$.

   There are two kinds of strata. 
     The first are {\em parabolic\footnote{Bestvina 
     {\em et al.} \cite{BFH} 
      use the terminology {\em non-exponentially-growing} strata} strata},  
      which are those of the 
     form $\S(x)$ with $x\notin \supp(y)$ 
      for all $y\in\supp(x)\ssm\{x\}$. 
     The second kind  
 are {\em exponential strata}, where one has $\S(x)=\S(y)$ for 
 some distinct 
 $x$ and $y$. The letter $x$ is defined to be {\em parabolic} 
 or {\em  exponential} according to the type of $\S(x)$.

        If $x$ is exponential then $|\phi^j(x)|$ grows exponentially with $j$.  If 
 all the edges of $\supp(x)$ are 
  parabolic then $|\phi^j(x)|$ grows polynomially 
  with $j$.  However, it may also happen that $x$ is a parabolic letter 
  but $|\phi^j(x)|$ 
  grows exponentially; this  
 will be the case if  $\supp(x)$ contains 
 letters $y$ such that $\S(y)$ is exponential.

\begin{example}  
 Define $\phi: F_3\to F_3$ by $a_1\mapsto a_1^2a_2,\ a_2\mapsto 
  a_1a_2,\ 
 a_3\mapsto a_1a_2a_3$. Then $\S(a_1)=\S(a_2) =  
 \{a_1, a_2\}$ is an 
 exponential stratum, while $\S(a_3)=\{a_3\}$ 
 is a parabolic stratum with $\supp(a_3)=\{a_1,a_2,a_3\}$. 
 \end{example} 
  
 \begin{remark}\label{induct} 
  The relation {\rm $[y< x$ if  
 $\S(y)\subset\supp(x)\ssm \S(x)]$} generates a partial 
 ordering on the letters $\{a_1,\dots,a_m\}$. For each 
 $x$, the subgroup of $F$ generated by $\pre(x)=\{y\mid y<x\}$ 
 is $\phi$-invariant. Let  $F\lfloor x \rfloor $ denote 
  the quotient of $\langle\supp(x)\rangle$ 
 by the normal closure of $\pre(x)\subset\supp(x)$, and 
 let $F\lceil x\rceil $  denote 
  the quotient of $F$ 
 by the normal closure of $\pre(x)\subset F$. Note 
that  $F\lfloor x \rfloor $  
is a free group with basis (the images of) the letters in $\S(x)$, and $F\lceil x\rceil $ is the free 
 group with basis $\{a_1,\dots,a_m\}\ssm\pre(x)$. 
 
The automorphisms of $\pre(x),\  F\lfloor x \rfloor$ 
and $F\lceil x\rceil $ induced by $\phi$ are positive 
with respect to the obvious bases, and their strata 
are images of the strata of $\phi$. 
 \end{remark} 
  
 \subsection{Conditioning the automorphism} 
  
 In the following proposition, the strata considered are those 
  of $\phi^k$. 
 (These may be smaller than the strata of $\phi$; consider 
 the periodic case for example.) 
  
  \begin{proposition}\label{power} There exists a positive 
   integer $k$ 
  such that $\phi_0:=\phi^k$ has the following properties: 
\begin{enumerate} 
\item[1.]       Each letter $x$ appears in its own image under $\phi_0$. 
\item[2.]       Each exponential letter $x$ appears 
 at least $3$ times in its own image under $\phi_0$. 
\item[3.]       For all $x$, each letter $y\in\supp(x)$ appears 
 in $\phi_0(x)$. 
\item[4.]       For all $x$ and all $j \geq 1$, the 
leftmost  and rightmost letters of $\phi_0^j(x)$ 
are the same as those of $\phi_0(x)$. 
\item[5.]  For all $x$, all $j\geq 1$ 
 and all strata $\S\subseteq\supp(x)$, 
  the leftmost  (respectively, 
  rightmost) letter 
 from $\S$ in the reduced word $\phi_0^j(x)$ is the same as 
  the leftmost (resp. 
 rightmost) 
 letter from $\S$ in $\phi_0(x)$. 
\end{enumerate} 
\end{proposition} 
 
\begin{proof} Items (1) to (3) can be seen as simple facts about 
positive integer matrices, read-off from the action of $\phi$ 
on the abelianization of $F$. 
(By definition $a_j\in \S(a_i)$ if and only if the $(i,j)$ 
 entry of some power 
of the matrix describing this action is non-zero.)

Assume that $\phi_1$ is a power of $\phi$ 
that satisfies (1) to (3). Note that (3) implies 
that the strata of $\phi_1$ coincide with those of any proper 
power of it.  
 
Replacing $\phi_1$ by a positive power if necessary, we may 
assume that if $\phi_1^j(x)$ begins with the letter $x$, for 
any $j\ge 1$, then $\phi_1(x)$ begins with $x$. This ensures  
that {\em{$[y\preceq_L x$ if some $\phi^j(x)$ begins with $y]$}} 
 is a partial 
ordering, for if $\phi_1^{j_k}(x_k)$ begins with $x_{k+1}$ for 
$k=1,\dots,r$ and if $x_{r+1}=x_1$,  
then $\phi_1^{\Sigma j_k}(x_1)=x_1$ 
and hence $x_1=x_2=\dots = x_r$. 

If $\phi_1(x)$ begins with $z$ then $z\preceq_L x$, so 
by raising $\phi_1$ to a suitable power 
we can ensure for all $x$ 
that $\phi_1(x)$ begins with a 
letter that is $\preceq_L$-minimal. The $\preceq_L$-minimal 
letters $y$ are precisely those such that $\phi_1(y)$ begins with $y$. 
An entirely similar argument applies to the relation 
 {\em{$[y\preceq_R x$ if some $\phi^j(x)$ ends with $y]$}}. 
 This proves (4). 
 
Now assume that $\phi_0$ satisfies (1) to (4). The assertion  
in (5) concerning leftmost letters from $\S$ is clear 
for those $x$ where  $\phi_0(x)$ begins with $x$. If $\phi_0(x)$ begins
with $y\neq x$, then either $\S\subset\supp(y)$ 
or else the occurrences 
of letters from $\S$ in $\phi_0^j(x)$ are 
in 1-1 correspondence with the occurrences in the  
image of $\phi_0^j(x)$ in $F\lceil y\rceil $.  (Notation of
Remark \ref{induct}.) In the latter case,
arguing by induction on the size of $\pre(y)$ we 
may assume that  the induced automorphism   $\lceil \phi_0\rceil_y 
:F\lceil y\rceil\to F\lceil y\rceil $ has 
the property asserted in (5); the desired conclusion  for $\phi_0^j(x)$
is then  tautologous. In the former case, if 
we replace $\phi_0$ by $\phi_0^2$ then the conclusion 
becomes as immediate as it was when $\phi_0(x)$ began with $x$.  
 
An entirely similar argument applies to rightmost letters. 
\end{proof} 
  
\begin{remark} 
Although we shall have no need of it here, it seems worth 
recording that item (5) of the above proposition 
remains true if one replaces strata $\Sigma \subset \supp(x)$ 
by supports $\supp(y)\subset\supp(x)$. 
\end{remark} 
 
\smallskip

 \begin{quote}{\em We now fix an automorphism $\phi=\phi_0$ and assume that 
 is satisfies conditions (1)-(5) above. 
 All of the constants discussed in the sequel 
 will be calculated with respect to this $\phi$.} 
\end{quote} 
 
\section{Preferred Futures, Fast Letters and Cancellation}  
\label{PrefFutSec} 
 
Having conditioned our automorphism appropriately, we 
are now in a position to analyse the fates of (blocks of) edges 
 as they evolve in time. 
  
 \begin{definition}[Preferred futures]\label{pref-fut} 
For each basis element $x\in\{a_1,\dots,a_n\}$,  
 we choose an occurrence of $x$ in the reduced word 
 $\phi(x)$ to be the (immediate) {\em preferred future of $x$}: 
   if $x$ is a 
  parabolic letter, there is only one possible choice; 
  if $x$ is an  
  exponential letter, 
   we choose an occurrence of $x$ that is neither  
  leftmost nor rightmost (recall that we have 
  arranged for $x$ to appear 
 at least three times in $\phi(x)$). More generally, we 
 make a recursive definition of the {\em preferred future 
 of $x$ in $\phi^n(x)$}: 
 this is the occurrence of $x$ in $\phi^n(x)$ that 
 is the preferred future of the 
 preferred future of $x$ in $\phi^{n-1}(x)$.  
  
 The above definition distinguishes an edge $\e_1$  on the top of 
 each 2-cell in our diagram $\Delta$, namely the edge 
 labelled by the preferred future of the label at 
 the bottom $\e_0$ of the 2-cell. We define $\e_1$ to 
 be the (immediate) 
  {\em preferred future} of $\e_0$. As with letters, 
 an obvious recursion then defines a preferred future of $\e_0$ 
 at each step in its future (for as long as it continues 
 to exist). 
  
 Note that $\e_0$ has at most  one preferred 
 future at each time. (It has exactly one until a preferred 
 future dies in a corridor, 
 lies on the boundary, or hits a singularity.) 
  
 If the bottom edge of a 2-cell is $\e_0$, then we define 
 the preferred future of that 2-cell at time $t$ to be the unique 2-cell 
 at time $t$ whose 
 bottom edge is the preferred future of $\e_0$. 
 \end{definition}

\subsection{Left-fast, constant letters, etc.} 
We divide the letters $x\in\{a_1^{\pm 1},\dots,a_m^{\pm 1}\}$ into 
classes according to the growth of the words 
$\phi^j(x), j=1,2,\dots$, and divide the edges of $\Delta$ into 
classes correspondingly. 
\begin{enumerate} 
\item[$\bullet$] 
If $\phi(x) = x$ then $x$ is called a {\em constant 
 letter}.  
  
\item[$\bullet$] If $x$ is a {\em non}-constant letter, then 
  the function $n\mapsto |\phi^n(y)|$ grows 
 like a polynomial of degree $d\in\{1,\dots,m-1\}$ or else as an exponential 
   function of $n$. 
 
\item[$\bullet$] Let $x$ be a non-constant letter. 
 If the distance between the preferred future of $x$ and the   
 beginning of the word $\phi^n(x)$ grows at least quadratically as 
 a function of $n$, we say that 
 $x$ is {\em left-fast}; if this is not the case, 
 we say that $x$ is 
 {\em left-slow}.  {\em Right-fast} and  
 {\em right-slow} are defined similarly. Note that $x$ is 
left-fast (resp. slow) if and only if $x^{-1}$ is right-fast (resp. slow). 
 
\item[$\bullet$] Let $x$ be a non-constant letter. If $\phi(x) = uxv$ (the shown occurrence of $x$ need not be the preferred future), where $u$ consists only of constant letters,   
 then we say that $x$ is {\em \lpl}. (We place no restriction on $v$; in particular 
 it may contain occurrences of $x$.) {\em Right para-linear} is defined 
 similarly. 
\end{enumerate}

\bd 
For \lpl letters, we define the {\em (left) para-preferred future} 
 (pp-future) to be the left-most occurrence of $x$ in $\phi(x)$. 
  The (right) pp-future of a \rpl letter is defined similarly, and 
  edges in $\Delta$ inherit these designations from their labels. 
  
  (It is possible that a letter  
   might be both \lpl  and right para-linear, and in such cases the 
   left and right 
    pp-futures need not agree. But when we discuss pp-futures, 
    it will always be clear  
     from the context whether we are favouring the left or the  right.) 
\ed

The following lemma indicates the origin of the 
terminology `left-fast' (cf.~\cite[Lemma 4.2.2]{BFH}).  
(A slight irritation arises from the fact that 
there may exist letters $x$  such that $x$ is not left-fast but  
 $\phi(x)$ contains left-fast letters; this difficulty accounts 
 for a certain clumsiness in the statement of the lemma.)

\begin{lemma}\label{C_0} There exists a constant $C_0$ with the 
following property: if $x\in\{a_1,\dots,a_n\}$ is such that 
 $\phi(x)$ contains a left-fast letter $x'$ 
and if $UVx\in F$ is a reduced word with $V$  positive\footnote{i.e. no 
 inverses $a_j^{-1}$ appear 
in $V$} and $|V|\ge C_0$, 
then for all $j\ge 1$, the preferred 
future of $x'$ is not cancelled 
when one freely reduces $\phi^j(UVx)$. Moreover, 
$|\phi^j(UVx)|\to\infty$  as $j\to\infty$. 
\end{lemma} 
 
\begin{proof} We factorize the 
reduced word $\phi^j(x)$ as $Y_{x,j}x'Z_{x,j}$ to emphasise the 
placement of the preferred future of a fixed left-fast letter 
$x'$ from $\phi(x)$. The fact that $x'$ is left-fast implies that $j\mapsto |Y_{x,j}|$ grows at least quadratically.  
 
Fix $C_0$ sufficiently large to 
ensure that for each of the finitely many possible 
$x\in\{a_1,\dots,a_n\}$, the integer 
$|Y_{x,j}|$ is greater than $Bj$ whenever $j\ge C_0/B$, 
where $B$ is the bounded cancellation constant. 
 
The Bounded Cancellation Lemma assures us that during the 
free reduction of the naive expansion of $\phi(UVx)$, 
at most $B$ letters of the positive word $\phi(Vx)$ will 
be cancelled. At most $B$ further letters will be cancelled 
when the naive expansion of $\phi^2(UVx)$ 
is freely reduced, and so on. Since $V$ and $\phi$ are positive and 
$ |V| \ge C_0$, it follows that $\phi^j(V)$ will 
not be completely cancelled during the free reduction  
of $\phi^j(UVx)$ if $j\le C_0/B$. When $j$ reaches $j_0:=\lceil C_0/B\rceil$ the 
distance 
from the preferred future of $x'$ to the left end of 
the uncancelled segment of $\phi^j(Vx)$ is 
at least  $|Y_{x,j_0}|$, which is greater than $Bj_0$ and hence $C_0$. 
Repeating the argument with $Y_{x,j_0}$ in place of $V$, we conclude that 
the length of the uncancelled segment of $\phi^j(Vx)$ in $\phi^j(UVx)$ 
remains positive and goes to infinity  with $j$. 
\end{proof}

Significant elaborations of the previous argument will be developed in 
Section \ref{ConstantSection}.  
 
\begin{definition}[New edges, cancellation and consumption]\label{new} 
Fix a 2-cell in $\Delta$. One edge in the top of 
the cell is the preferred future of the bottom 
edge; this will be called  {\em old} and the 
remaining edges will be called {\em new}. (These 
concepts are unambiguous relative to a fixed 2-cell or (folded) corridor, but `old edge' would be 
ambiguous if applied simply to a 1-cell of $\Delta$.) 
 
Two (undirected) edges $\e_1, \e_2$ 
in the naive top of a  
corridor are said to {\em cancel} each other if their images in the 
folded corridor coincide. If $\e_1$ lies to the 
left\footnote{Recall 
that corridors have a left-right orientation.} of 
$\e_2$, we say that $\e_2$ has been cancelled {\em 
from the left} and $\e_1$ has been cancelled {\em 
from the right}. 
If $\e_1$ is the preferred future of an edge $\e$ 
in the bottom of the corridor and $\e_2$ is a new 
edge in the 2-cell whose bottom is $\e'$, then we 
say that $\e'$ has {\em (immediately) consumed} $\e$ 
{\em from the right}. `Consumed 
from the left' is defined similarly. 
 
Let $e$ and $e'$ be edges in $\bot(S)$ for some 
corridor $S$, with $e$ to the left (resp. right) of $e'$. 
If an edge in the future of $e$ 
cancels a preferred future of $e'$, then we say 
that $e$ {\em eventually consumes} $e'$ {\em from 
the left (resp. right).}  
\end{definition} 
\begin{lemma} \label{NoOldCanc} 
A pair of old edges cannot cancel each other. 
\end{lemma} 
 
\begin{proof} 
Suppose that 
two old edges in the naive top of a corridor $S$ 
are labelled $x$ and cancel each other.  These 
edges are the  preferred futures of edges on $\bot(S)$ 
that bound an arc $\alpha$ labelled by a reduced word 
 $x^{-1}wx$. 
Consider the freely-reduced factorisation $\phi(x) = uxv$ where 
the visible $x$ is the preferred future. 
The arc in the naive top of $S$ corresponding 
to $\alpha$ is labelled  
 $v^{-1}x^{-1}u^{-1}Wuxv$, where $W$ is the naive 
 expansion of $\phi(w)$. The old edges that we are considering 
 are labelled by the visible occurrences of $x$ in this word and 
 our assumption that these edges cancel means that the subarc 
 labelled $x^{-1}u^{-1}Wux$ becomes a loop (enclosing a 
 zero-area sub-diagram) in the diagram $\Delta$. 
  
 But this is impossible, because $x^{-1}wx$ is freely reduced, 
which means that $W$ is not freely equal to the empty 
 word, and hence neither is $x^{-1}u^{-1}Wux$. 
 \end{proof}

 \begin{corollary}\label{parabolicC} 
  An edge labelled by a 
 parabolic letter $x$ 
 can only be consumed by an edge labelled $y$ with   
 $\supp(x)$ strictly contained in $\supp(y)$. 
 \end{corollary} 
  
\begin{remark} 
A non-constant letter can only be (eventually) consumed from the left (resp. right) by a right-fast  
(resp. left-fast) letter. 
\end{remark}

\begin{remark} The number of old letters in 
 the naive top of a corridor $S$ is $|S|$, so 
 the length of corridors in the future of $S$ 
 will grow relentlessly unless old letters are 
 cancelled by new letters or the corridor hits a 
 boundary or a singularity. 
 \end{remark}

An obvious separation argument provides 
us with another useful observation concerning cancellation:

\begin{lemma}\label{perfect} Let $\e_1,\ \e_2$ and $\e_3$ be three 
(not necessarily adjacent) edges that appear in 
order of increasing subscripts as one reads from 
left to right along the bottom of a corridor. If 
the future of $\e_2$ contains an edge of $\partial\Delta$ 
or of a singularity, then no edge in the future of 
$\e_1$ can cancel with any edge in the future of $\e_3$. 
\end{lemma}

\section{Counting Non-constant Letters} \label{NonConstantSection} 
 
In this section we fix a corridor $S_0$ in $\Delta$ and 
bound the contribution of non-constant letters to the 
length of $\bot(S_0)$.

\subsection{The first decomposition of $S_0$}\label{decomp} 
 
Choose an edge $\e$ on the bottom of $S_0$.  As we follow the 
preferred future of $\e$ forward one of the following (disjoint) events must 
occur: 
 
\begin{enumerate} 
\item[1.] The last 
  preferred future of $\e$ lies on the boundary of  
 $\Delta$. 
  
\item[2.]  The last 
  preferred future of $\e$ lies in a singularity. 
 
\item[3.]   The last 
  preferred future of $\e$ dies in a corridor $S$ (i.e. 
  cancels with another edge from the naive top of $S$). 
\end{enumerate} 
 
We shall bound the length of $S_0$ 
by finding a bound on the number of edges in each of these three 
cases.  
 
We divide Case (3)  into two sub-cases:

\begin{enumerate} 
\item[3a.] 
 The preferred future of $\e$ dies when it is cancelled by an edge 
that is not in the future of $S_0$.  
 
\item[3b.] 
The preferred future of $\e$ 
dies when it is cancelled by an edge 
that is in the future of $S_0$.  
\end{enumerate} 
 
\subsection{Bounding the easy bits} \label{EasyBounding} 
 
Label the sets of edges in $S_0$ which fall into the above classes 
$S_0(1), S_0(2), S_0(3a)$ and $S_0(3b)$ respectively.   
 We shall see that $S_0(3b)$ is by far the most troublesome 
 of these sets. 
    
The first of the bounds in the following lemma is obvious, and the 
second follows immediately from  Proposition \ref{SingularityProp}.

\begin{lemma}\label{bound1and2} 
 $|S_0(1)| \leq \n \text{   \rm{and}   } 
  |S_0(2)| \leq 2B\n$.    
  \end{lemma}

\begin{lemma}\label{bound3a} $|S_0(3a)|\le B\n$. 
\end{lemma} 
 
\begin{proof} The preferred future of each $\e\in S_0(3a)$ 
dies in some corridor in the future of $S_0$. Since 
there are less than $\n /2$ corridors, we will be done 
if we can argue that the preferred future of at most $2B$  
such edges can die in each corridor $S$. 
 
Lemma \ref{Connected} tells us that the future of $S_0$ 
intersects $S$ in a connected region, the 
bottom of which is an interval $I$. The Bounded Cancellation Lemma 
assures us that only the edges within a distance $B$ of the 
ends of $I$ can be consumed in $S$ by an edge from 
outside the interval. And by definition, if a preferred future 
of an edge from  $S_0(3a)$ is to die in $S$, then it must 
be consumed by an edge from outside $I$. 
\end{proof} 
 
We have now reduced Theorem \ref{BoundS} to 
the problem of bounding $S_0(3b)$, 
i.e. of understanding cancellation {\em within} the future of $S_0$. 
This will require a great deal of work. As a first step, 
we further decompose 
$S_0$, mingling the above decomposition based on the fates 
of preferred futures of  
edges with the natural decomposition of $S_0$ into 
colours, as defined in Definition \ref{defMu}.

\subsection{The chromatic decomposition of $S_0$} \label{chromatic} 
 
We fix a colour $\mu$ and 
write $\mu(S_0)$ to denote the interval of $\bot(S_0)$ 
consisting of edges coloured $\mu$.  
We shall abuse terminology to the 
extent of referring to $\mu(S_0)$ as {\em a colour}, evoking 
the mental picture of the 2-cells in $S_0$ being painted 
with their respective colours. (Recall 
that the 2-cells of $S_0$ are in 1-1 correspondence with 
the edges of $\bot(S_0)$.) 
 
We shall subdivide $\mu(S_0)$ into five subintervals 
according to the fates  
of the preferred futures of edges. To this end,  
we define  $l_{\mu}(S_0)$ to be the rightmost edge in $\mu(S_0)$   
whose immediate future contains  a left-fast edge that is 
ultimately consumed 
from the left by an edge of  $S_0$, and we define 
 $A_1(S_0,\mu)$ to be the set of edges in $\bot(S_0)$ 
 from the left end of $\mu(S_0)$ to $l_{\mu}(S_0)$, inclusive. 
We define $A_2(\mu,S_0)\subset\mu(S_0)$ 
 to consist of the remaining 
 edges in $\mu(S_0)$ whose 
 preferred futures  are  ultimately consumed 
from the left by an edge of  $S_0$. 
 
Similarly, we define $r_{\mu}(S_0)$ to be 
 the leftmost edge $\mu(S_0)$ that has a  right-fast edge in its immediate future 
 that is ultimately consumed 
from the right by an edge of  $S_0$, and we define 
 $A_5(S_0,\mu)$ to be the set of edges in $\bot(S_0)$ 
 from the right end of $\mu(S_0)$ to $r_{\mu}(S_0)$, inclusive. 
We define $A_4(\mu,S_0)\subset\mu(S_0)$ 
 to consist of the remaining 
 edges in $\mu(S_0)$ whose 
 preferred futures  are  ultimately consumed 
from the right by an edge of  $S_0$. 
 
Finally, we define $A_3(S_0,\mu)$ to be the 
 remainder of the edges in $\mu(S_0)$.

\begin{figure}[htbp] 
\begin{center} 
  
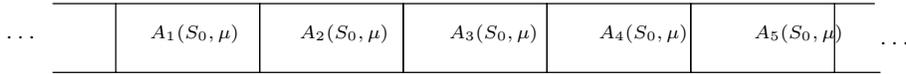 
 
\caption{The second decomposition of $S_0$} 
\label{figure:A1A5} 
\end{center} 
\end{figure} 
 
\medskip

Modulo the fact that any of the $A_i(S_0,\mu)$ 
might be empty, Figure 
10 is an accurate portrayal of  $\mu$: 
the $A_i(S_0,\mu)$ are  connected 
and they occur in 
ascending order of suffix from left to right. 
 
The chromatic decomposition of $S_0$ is connected to the 
decomposition of Subsection \ref{decomp} by the equality 
in the following lemma, which is a tautology. The  
inequality in this lemma is a restatement 
of Lemmas \ref{bound1and2} and \ref{bound3a}. 
 
\begin{lemma} \label{A3Lemma} 
$$ 
 \bigcup_{\mu} A_3(S_0,\mu) = S_0 \ssm S_0(3b)\ \ \  
 \text{  {\rm{and}}  }\ \ \  
 \sum_{\mu}|A_3(S_0,\mu)| < \left({3B} + 1\right)\n . 
 $$ 
 \end{lemma}    
 
Thus the following lemma is a step towards bounding the 
size of $S_0(3b)$. 
 
\begin{lemma} \label{A1A5Lemma} 
 
$$ 
|A_1(S_0,\mu)|  \leq  C_0 \ \ \  
 \text{  {\rm{and}}  }\ \ \  
|A_5(S_0,\mu)|  \leq  C_0. 
$$ 

\end{lemma} 
 
\begin{proof} 
We prove the result only for $A_1(S_0,\mu)$; 
the proof for $A_5(S_0,\mu)$ is entirely similar.

As in 
Lemma \ref{perfect}, we know that the entire  future of the edges of  
$A_1(S_0,\mu)$ to the left of $l_\mu(S_0)$ must 
eventually be consumed from the left  by edges of $S_0$. This means 
that we are essentially in the setting of Lemma \ref{C_0}, with  
$l_{\mu}(S_0)$ in the role 
of $x$ and $A_1(S_0,\mu)$ in the role of $Vx$.

Thus if the length of $A_1(S_0,\mu)$ were greater than $C_0$,  
then we would conclude that  no  left-fast edge in the  
immediate future of  $l_{\mu}(S_0)$ would be cancelled from the left by 
an edge  of $\bot(S_0)$, contradicting the definition of $l_{\mu}(S_0)$.   
\end{proof} 
 
\begin{corollary} 
$$ 
\sum_{\mu}|A_1(S_0,\mu)| \, \leq \, C_0\n 
 \ \ \  
 \text{  {\rm{and}}  }\ \ \  
\sum_{\mu}|A_5(S_0,\mu)| \, \leq \, C_0\n . 
$$ 
\end{corollary} 
 
\subsection{A further decomposition of $A_2(S_0,\mu)$ 
 and $A_4(S_0,\mu)$}

It remains to bound $A_2(S_0,\mu)$ and $A_4(S_0,\mu)$.  
We deal only with $A_4(S_0,\mu)$, the argument for $A_2(S_0,\mu)$ 
being entirely similar. 
 
\begin{figure}[htbp] 
\begin{center} 
  
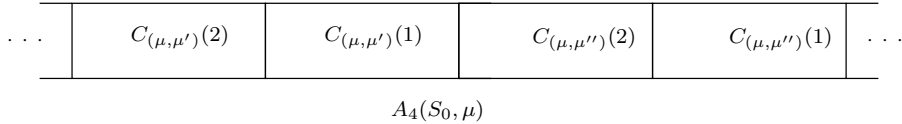 
  
\caption{$C_{(\mu,\mu')}(1)$ and $C_{(\mu,\mu')}(2)$.} 
\label{figure:Cmumu} 
\end{center} 
\end{figure}

First partition $A_4(S_0,\mu)$ into subintervals $C_{(\mu,\mu')}$ 
that consist of edges  that are eventually consumed by edges of a specified 
colour $\mu'$. Then partition $C_{(\mu,\mu')}$ into two subintervals:  
$C_{(\mu,\mu')}(1)$ begins at the   
right of  $C_{(\mu,\mu')}$ 
and ends with the last non-constant edge;  
 $C_{(\mu,\mu')}(2)$ consists of the remaining (constant) edges. 
See Figure \ref{figure:Cmumu}. 
\medskip

In the course of this section we will bound the size of the intervals $C_{(\mu,\mu')}(1)$ and during the following four sections we bound the sum over all pairs $(\mu,\mu')$
of the  sizes of the intervals $C_{(\mu,\mu')}(2)$ to get the desired bound on $|S_0(3b)|$.  In order to control this sum, we have 
to address the question of which colours can be adjacent. 
 
\subsection{Adjacent Colours}

In Corollary \ref{muConn} we saw that in any corridor 
$S$, the edges in $\bot(S)$ of a fixed colour form an interval. 
We say that two distinct colours $\mu$ and $\mu'$ are 
{\em adjacent} in $S$ if the closed intervals 
$\mu(S)$ and $\mu(S')$ 
have a common endpoint in $\bot(S)$. (Equivalently, 
there is a pair of 2-cells in $S$, one  coloured $\mu$ and 
the other $\mu'$, that share an edge labelled $t$.) 
We write 
 ${\vecZ}$ to denote the set of ordered pairs  
$(\mu,\mu')$ such that  
 $\mu$ and $\mu'$ are adjacent in some corridor $S$ 
with $\mu(S)$ to the left of  $\mu'(S)$  in $\bot(S)$, 
and we write $\vecZ$ to denote the set of unordered 
pairs.  
 
\begin{lemma} \label{NoOfAdjacencies} 
$$ 
|{\vecZ}| < 2\n -3 .     
$$ 
\end{lemma} 
 
\begin{proof}  
 We shall express this proof in the language 
of the forest $\mathcal F$ introduced in Remark \ref{tree}. Suppose 
that $\mu$ and 
$\mu'$ are adjacent in $S$. 
In $S$ we can connect the centre 
of some 2-cell coloured $\mu$ to the centre of some 2-cell 
coloured $\mu'$ by an arc contained in the union of 
the pair of 2-cells. The union of this arc and the trees in $\mathcal F$ 
corresponding to the colours $\mu$ and $\mu'$ disconnects 
the disc $\Delta$; each of the other trees in $\mathcal F$ 
is entirely contained in a 
component of the complement, and the 
colours with trees in different components can 
never be adjacent in any corridor. 
 
We can encode adjacencies of colours by a chord diagram: draw 
a round circle with marked points representing the colours of 
$\Delta$ in the cyclic order that they appear in $\partial\Delta$, 
then connect two points by a straight line if the corresponding 
colours are adjacent in some corridor.  
 The final phrase of 
the preceding paragraph tells us that the lines in this 
chord diagram do not intersect in the interior of the disc. 
A simple count shows that since there are less than 
$\n$ colours, there are less than $2\n -3$ lines in this diagram. 
\end{proof}

\subsection{Non-constant letters in 
 $C_{(\mu,\mu')}$ that are not left-fast} 
 \label{NonConstantSubsect} 
 
We stated in the introduction that a careful analysis of 
van Kampen diagrams would allow us to reduce the Main Theorem to the 
study of blocks of constant letters.  
In this section we achieve the last step of this reduction. 
 
\begin{lemma} \label{C1Lemma}  
There is a constant $C_1$ depending only on $\phi$ 
with the following property: 
 
Let $S$ be a corridor and let $\mu_1$ and $\mu_2$ be 
 colours that occur in $S$ with $\mu_1$ to the left of $\mu_2$ (but do not 
assume that $\mu_1(S)$ is adjacent to $\mu_2(S)$).  Let 
 $I\subset A_4(S,\mu_1)$ 
  be a sub-interval that satisfies the 
 following conditions 
\begin{enumerate} 
\item[1.] the left-most edge of $I$ is 
 non-constant \mbox{and } 
   
\item[2.]  the preferred future of each edge in $I$ 
 is eventually consumed by an edge of $\mu_2(S)$.\\ 
\end{enumerate} 
\noindent Then $|I| \leq C_1$. 
In particular, $|C_{(\mu,\mu')}(1)| \leq C_1$ 
for all $(\mu,\mu') \in {\vecZ}$. 
 
It suffices to take $C_1 = 2mB^2$, where $m$ is the rank of $F$, and 
$B$ is the constant from the Bounded Cancellation Lemma. 
\end{lemma}

\begin{proof} 
  The region $I$ being considered contains no edge with a right-fast  
letter in the $\phi$-image of its label. 
 Since all exponential letters 
 are both left-fast and right-fast, all non-constant edges in the future of $I$ are  parabolic.   
 
We begin the argument at the stage in time where $\mu_2$ 
starts cancelling $I$.  For 
notational convenience we assume that this time 
is in fact $\height(S)$. (If it is not, then the  
fact that the length of $I$ may 
have increased in passing from $\time(S)$ to 
this time adds greater strength to the 
bound we obtain.) 
  
 We focus on the leftmost 
 edge $\e_0$ of $I$ that is labelled 
 by a non-constant letter $x$ for which $\supp(x)$  is maximal  
 among the supports of all edge-labels 
 from $I$ (with respect to inclusion). 
 Let $y$ be the label on the edge 
 $\e_0'$ of $\mu_2(S)$ that 
 eventually consumes $\e_0$ (oriented as shown in Figure \ref{C0Pic}). 
 Note that $\supp(x)$ is strictly contained in $\supp(y)$, 
 by Corollary \ref{parabolicC}. If $\e_0'$ consumes $\e_0$ immediately, 
 then the bounded cancellation lemma tells us that 
$\e_0$ is a distance less than $B$ from the righthand end 
of  $I$. If not, then we 
 proceed one step into the future\footnote{proceeding one 
step into the future also allows us to assume that there 
are no letters coloured $\mu_1$ to the right of $I$} 
 and  appeal to the 
 conditioning done in Proposition \ref{power}(5) to assume 
 that for all $j\ge 1$, the rightmost letter in $\phi^j(y)$ 
 whose support includes $x$ is $y$. 
 We shall call the edge in the future of $\e_0'$ 
carrying the rightmost $y$  
 the {\em highlighted} future of $\e_0'$ (perhaps it is not 
 the preferred future).  

\begin{figure}[htbp] 
\begin{center} 
  
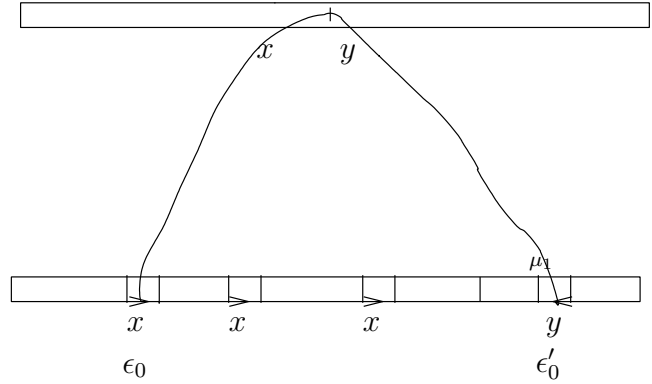 
  
\caption{The edge labelled $\epsilon_0'$ will eventually consume
$\epsilon_0$.}  
\label{C0Pic} 
\end{center} 
\end{figure} 
  
  The first important point to observe is that  
  the maximality of $\supp(x)$ ensures that 
  there will never be any new edges labelled $x$ 
  in the future of $I$ 
  (`new' in the sense of \ref{new}). 
   
 The second important point to note is that  
 the edges labelled $x$ in the future of 
 $\e_0'$ that are to cancel with 
 the futures of the 
 edges labelled $x$ in $I$ must all lie to the 
 left of the highlighted future of $\e_0'$. The point here is 
 that the highlighted future of $\e_0'$ 
 cannot be cancelled by an edge of $I$ (by the maximality 
 of $x$), and in order for it to be cancelled from the 
 other side, all the edges to 
 its right labelled $x$ 
 would have to be cancelled first, which would mean that they too 
 were cancelling with something not in the future of $I$. 
  
 We now come to the key observation of the proof: at 
 each stage $j$ steps into  the  
  future of $S$, the leftmost\footnote{we 
  have already noted that this is to 
  the left of the highlighted future of $\e_0'$} 
   edge $\e_j'$ in the future of $\e_0'$ 
  that is labelled  
  $x$ must be cancelled by an edge from 
  the future of $I$ {\em immediately}, i.e. in the corridor 
  where it appears at  $\height(S)+j$. 
   Indeed if this 
  were not the case, then  $\e_j'$ would develop a preferred 
  future which, being an old  
   edge (in the sense of Definition \ref{new}), could only  
  cancel with a  new edge (Lemma \ref{NoOldCanc}) 
  in the future of $I$. And since 
  we have arranged that there be no new edges labelled $x$, 
  the preferred future of $\e_j'$ would never cancel with 
  an edge in the future of $I$. But this cannot be, because 
  the continuing existence of a preferred future for $\e_j'$ would prevent 
  anything to its {\em right} consuming  an 
  edge in the future of $I$, and the penultimate sentence in the 
  third paragraph of this proof implies that no new 
edges labelled $x$ will ever appear to its {\em left} in the future of $\e_0'$. 
Thus if $\e_j'$ is not 
  cancelled immediately then we have a  contradiction 
  to the fact that $\e_0'$ must 
 eventually  consume $\e_0$.

  We have just proved that at $\height(S)+j$ the 
  edge $\e_j'$ must cancel with the preferred 
  future of an edge $\e_j$ in $I$ that is labelled $x$.  
  According to 
  the Bounded Cancellation Lemma, the preferred 
  future of $\e_j$ at $(\height(S)+j-1)$ must lie within 
  a distance $B$ of the right end of the future of $I$.  
  Since there is  no cancellation within the 
  future $I$, an iteration of this argument shows that 
   for as 
  long as there exist edges labelled $x$ in the future 
  of $I$, each successive pair of these edges is separated 
  by less than $B+|\phi(y)|\le 2B$ edges at each moment in time, 
  and the rightmost must be within a distance $B$ of the 
  right end of the future of $I$. 
   
  But since $\phi(x)$ contains at least 
  one letter other than  the preferred future of $x$, 
  it follows that there cannot be a pair 
  of  edges of $I$ labelled $x$ that remain unconsumed 
  at  $\height(S)+2B$, for otherwise 
  they would have grown a distance more than 
  $2B$  apart, contradicting the 
  conclusion of the previous paragraph. And proceeding 
  one more step into the future, the last edge labelled $x$ 
  must be consumed.  
   
  Since at most 
  $B$ letters of $I$ are cancelled at the right  
   at each stage in its future, all of the edges of $I$ labelled $x$ 
  are within a distance less than $2B^2$ of the right end of $I$, 
  and they are all consumed when $I$ has flowed $2B$ steps 
  into the future. 
  If no non-constant edges remain in the future of $I$ 
  at this stage, then we know 
  that $|I|\le 4B^2$. 
   
  If there do remain non-constant edges, we take the maximal interval of the 
   future of $I$ at   $\height(S)+2B$ whose leftmost 
   edge is non-constant, and we repeat the argument. (This 
   interval is obtained from the complete future of $I$ by 
   removing a possibly-empty collection of constant edges 
   at its left extremity.) 
    
  We proceed in this manner. The interval that 
  we begin with at each iteration has strictly fewer 
  strata than the previous one 
  and therefore the procedure 
  stops before $m=\text{\rm{rank}}(F)$ iterations. At 
  the time when it stops (at most $\height(S)+2mB$), the future 
  of $I$ has been cancelled entirely, except possibly for 
  a block of 
   constant edges at its left extremity.  
   With one final appeal to the bounded cancellation 
   lemma, we deduce that  $|I|\le 2mB^2$. 
  \end{proof}

  \begin{corollary}\label{C1Corollary} 
\[ \sum_{(\mu,\mu') \in {\vecZ}}|C_{(\mu,\mu')}(1)| < 2C_1\n . \] 
\end{corollary} 
 
\begin{proof} This follows immediately 
 from Lemmas \ref{NoOfAdjacencies} and \ref{C1Lemma}. 
\end{proof}

\section{The Bound on $\sum\limits_{\mu \in S_0}|A_4(S_0,\mu)|$ and $\sum\limits_{\mu \in S_0}|A_2(S_0,\mu)|$}\label{A4sec} 

The sum of our previous arguments has reduced us to the nub of the
difficulties that one faces in trying to prove the Main Theorem,
namely the possible existence of large blocks of constant letters in
the words labelling the bottoms of corridors.
Now we must obtain a bound on
 $$\sum\limits_{(\mu,\mu') \in {\vecZ}} |C_{(\mu,\mu')}(2)| 
$$ 
that will enable us to bound 
$\sum\limits_{\mu \in S_0}|A_4(S_0,\mu)|$ and\footnote{In 
practice 
 we  only need concern ourselves with $A_4$, the arguments for $A_2$ 
being entirely similar}   $\sum\limits_{\mu \in 
S_0}|A_2(S_0,\mu)|$ by a linear function of $\n$.  These are the final 
estimates required to complete the proof of the Main Theorem --- see 
Section \ref{summary} for a r\'esum\'e of the proof.  

 The regions $C_{(\mu,\mu')}(2)$ are static, in the sense 
that they do not change under iteration by $\phi$, so  
the considerations of future growth 
that helped us so much in previous 
sections cannot be brought to bear directly. Rather, we must 
analyse the complete history of blocks of constant letters, 
understand how large blocks come into existence, and 
 use global considerations to limit the sum of 
 the sizes of all such blocks. 
  
Because of the global 
nature of the arguments, 
 we shall not obtain bounds on the sizes of the individual sets $C_{(\mu,\mu')}(2)$.
 Instead, 
we shall identify an associated block of 
constant letters
elsewhere in 
the diagram  (a ``team") that is amenable to a delicate string of 
balancing arguments that facilitates a bound on a union of 
associated regions $C_{(\mu,\mu')}(2)$. 

Our strategy is motivated by the following considerations. 
Believing Theorem \ref{BoundS} to be true, we seek 
payment from the global geometry of $\Delta$ to compensate 
us for having to handle the troublesome blocks of constant 
edges $C_{(\mu,\mu')}(2)$;  the currencies of payment are 
{\em consumed colours} 
and dedicated subsets of 
edges on $\partial\Delta$ --- since $\Delta$ can have at 
most $\n$ of each, if we prove that adequate payment is available
then our troubles will be bounded and the Main Theorem  
will follow.  
The chosen currencies are apposite
because, as we shall  see in Section \ref{ConstantSection}, 
a large block of edges labelled by constant letters can only 
come into existence if  a colour (or colours) associated to a 
 component of this block in the past was consumed completely, 
  or else the boundary of $\Delta$ intruded into the past of 
   the block (or else something nearby) causing smaller regions of constant edges to elide.

In the remainder of this section we shall explain how various estimates on the behaviour of 
blocks of constant letters in $\Delta$ can be combined to obtain 
the bounds that we require on  
$\sum\limits_{\mu \in S_0}|A_4(S_0,\mu)|$ and  $\sum\limits_{\mu \in 
S_0}|A_2(S_0,\mu)|$. We hope that this explanation will provide the 
diligent reader with a useful road map and sufficient motivation to 
sustain them through the many technicalities needed to establish the 
estimates in subsequent sections.

In the following proposition, $M$ is the maximum length of the 
images $\phi(x)$ of the basis elements of $F$, while
$\ttt$ is the constant from the Pincer Lemma \ref{PincerLemma}, and
$C_1$ is the upper bound on the lengths of the intervals
$C_{(\mu,\mu')}(1)$ from Lemma \ref{C1Lemma}, $T_0$ comes from the Two
Colour Lemma \ref{TwoColourLemma} and $C_4$ comes from Lemma
\ref{G34pics}. The constant $\ll$ is
defined above Definition \ref{NestingDef}, and $B$ is the bounded cancellation
constant from Lemma \ref{BCL}.
\smallskip

\noindent{\bf{The Constant $K_1$ is defined to be}}
$$
 \AFourC.  
$$ 

\begin{proposition}\label{SummaryLemma} 
\[      \sum_{\mu \in S_0}|A_4(S_0,\mu)| \leq K_1  \n . 
 \]  
\end{proposition}

\subsection{Dramatis Personae} 
The ``proof" that we are about to present is essentially a scheme for
reducing 
the proposition to a series of technical lemmas that will be proved
in Sections \ref{teamSec} and \ref{BonusScheme}. These lemmas are
phrased in the language associated to {\em teams}, the precise
definition of which will
also be given  in  Section \ref{teamSec}. 
Many of the proofs involve global cancellation arguments
based on the {\em Pincer Lemma}, which will be proved in the
next section. 
Intuitively speaking, a {\em team}  
(typically denoted $\T$) is a  contiguous 
region of $\|\T\|$ constant letters all of which 
are to be consumed by a  fixed \lpl edge (the {\em reaper}). Notwithstanding 
this intuition, it is preferable for 
technical reasons  to define a team to be a set of pairs of colours $(\mu,\mu')\in\vecZ$, 
where $\mu'$ is fixed and the different {\em members} of the team correspond to  
different values of $\mu$. We write $(\mu,\mu')\in\T$ to denote 
membership. Teams also have {\em virtual members}, denoted $(\mu,\mu') \vin \T$ (see 
Definition \ref{Virtual}). There are less than $2\n$ teams (Lemma \ref{allIn}).
 
Each pair $(\mu,\mu')$ with $C_{(\mu,\mu')}(2)$ non-empty 
is either a member or a  virtual member  of a team (Lemma \ref{allIn}).  
There are {\em short} teams (Definition \ref{newTeams}) and long teams, 
 of which some are {\em distinguished} (Lemma \ref{Aget2Lemma}). 
There are four types of {\em genesis} of a team, (G1), (G2), (G3) and
(G4) (see Subsection \ref{genesis}).  Teams of genesis (G3) have
associated to them a pincer $\Pin_{\T}$ (Definition \ref{pl}) yielding
an auxiliary set of colours
$\subT$. There is also a set of colours $\chi_P(\T)$ associated to the
time before the pincer $\Pin_{\T}$ comes into play.  For long,
undistinguished teams, we also need to consider certain sets
$\CT$ and $\chi_{\delta}(\T)$ of colours consumed in the past of $\T$ (see
the proof of Lemma \ref{Aget2Lemma}). Such teams may
also have  three sets of edges in $\partial\Delta$ associated to  
 them: $\partial^\T$, $\down_1(\T)$ and  
$\down_2(\T)$. An important feature of the definitions of 
 $\partial^\T$ and $\down_1(\T)$ is that the sets associated 
  to different teams are disjoint.   This disjointness is crucial 
   in  the following proof, where we use the fact that the sum 
    of their cardinalities is at most $\n$. Similarly, the disjointness of
the sets $\chi_c(\T)$ is used to estimate the sum of their cardinalities by
$\n$ and likewise for $\chi_{\delta}(\T)$ and $\chi_P(\T)$.

It is not necessarily true that the sets $\down_2(\T)$ are disjoint
 for different teams, but we shall explain how to account for the
 amount of `double-counting' that can occur (see Lemma
 \ref{Aget2Lemma}).
 
Associated to every team one has  the time $t_1(\T)$ at 
which the reaper starts consuming the team (see Subsection
\ref{t1}). Teams genesis (G3) also have two
earlier times $t_2(\T)$ and $t_3(\T)$ associated to them as well as an
auxiliary set of edges $\QT$, the definitions of which
are somewhat  technical (see Definition \ref{PincerDef} {\em et seq.}). 

In Section \ref{BonusScheme} we describe a {\em bonus scheme} that
assigns a set  of extra edges, $\bonus(\T)$ to each team.  These
bonuses are assigned so as  to ensure that $|\bonus (\T)|+\|\T\|$
 dominates  the sum of 
 the cardinalities of the sets $\cmm$ 
associated to the  members and virtual members of $\T$. 
 
\smallskip 
\noindent{\bf Proof of Proposition \ref{SummaryLemma}.}

Recall that   $A_4(S_0,\mu)$ is partitioned  into disjoint regions $C_{(\mu,\mu')}$ 
which in turn are partitioned into $C_{(\mu,\mu')}(1)$ and 
$C_{(\mu,\mu')}(2)$.  
 
Given any $\mu_1$ and $\mu_2$, at most one  ordering of $\{\mu_1,\mu_2\}$ can 
arise in  $S_0$. Thus Lemma \ref{NoOfAdjacencies} 
implies that there are less than $2\n$ pairs $(\mu,\mu')\in\vecZ$ with 
$C_{(\mu,\mu')}\subset\bot(S_0)$  non-empty. 
It follows immediately from this observation and Lemma \ref{C1Lemma} that 
\[      \sum_{(\mu,\mu') \in {\mathcal Z}}|C_{(\mu,\mu')}(1)| \leq 2C_1\n . 
\]

Lemma  \ref{Aget2Lemma} accounts for the set of  distinguished 
 long 
teams $\dlong$: 
\[      \sum_{{\mathcal T} \in \dlong}\sum_{(\mu,\mu') \in  
{\mathcal T}}|C_{(\mu,\mu')}(2)| \leq 6B\n(T_1+T_0).   \]  
For all other teams $\T$ we rely on Lemma \ref{C1toTeamLength} which
states
\begin{equation}\label{goodEq} 
\sum_{(\mu,\mu') \in \T \mbox{ \tiny or } (\mu,\mu') \vin \T} |C_{(\mu,\mu')}(2)| \le
\|\T\|  + |\bonus(\T)| + B. 
\end{equation} 
We next consider the {\em genesis} of teams. All teams of genesis (G4)
are short (Lemma \ref{G4lemma}). And by Definition \ref{newTeams} for
the short teams $\T\in\Sigma$ we have
\[      \sum_{{\T} \in \Sigma}\sum_{(\mu,\mu') \in  
{\mathcal T}}|C_{(\mu,\mu')}(2)| \leq 2\ll\n + \sum_{\T\in\Sigma}
\big(|\bonus(\T)| + B\big).   \]

Lemma \ref{TeamAgeLemma} tells us that for teams of genesis (G1) and
(G2) we have 
\[      \|\T\| \leq 2MC_4|\down_1(\T)| + |\partial^{\T}|,  \]
whilst for teams of genesis (G3) we have
\[      \|\T\| \leq 2MC_4\big{(}|\down_1(\T)| +|\QT|\big{)}+
T_0\big{(}|\chi_P(\T)| + 1\big{)} +|\partial^{\T}| + \ll.       \]

Let $\Gthree$ denote the set of teams of genesis (G3) with $\QT$
non-empty. In Definition \ref{down2} we break $\QT$ into pieces so
that
$$
|\QT|\ =  t_3(\T) - t_2(\T)  
+ |\down_2(\T)|.
$$  
Making crucial use of the Pincer Lemma, in Corollary \ref{t1-t2Corr} we prove that
$$
\sum\limits_{\T \in \Gthree} 
 t_3(\T) - t_2(\T)  \ \le \ 3\ttt \, \n, 
$$ 
and in   Corollary \ref{downbound} we prove that
$$ 
\sum\limits_{\T \in \Gthree}|\down_2(\T)|
\leq (2 + 3\ttt + 5T_0)\n.
$$
This completes the estimate on $|\QT|$ and hence $\|\T\|$. 

Section 10 is dedicated to the proof of Proposition  \ref{BonusBound}, which states
\[      \sum_{\text{\small{teams}}}|\bonus(\T)| \leq \big( \Bb \big)\n. \]

Adding all of these estimates and recalling that there are less than $2\n$ teams, we deduce:
\[      \sum_{\mu \in S_0}|A_4(S_0,\mu)| \leq K_1 \n ,    \]
where $K_1$ is
$$ 
 \AFourC.
$$
Thus the proposition is proved.
\hfill$\square$
\smallskip 
\begin{remark} The stated value of the constant $K_1$ 
  is an artifact of our proof: we
have simplified the estimates at each stage for the sake
 of clarity rather than trying to optimise the
constants involved. Nevertheless, we have made some effort
 to make the arguments constructive
so as to prove that there exists an algorithm to calculate  
the Dehn function of $F\rtimes_\phi\mathbb Z$ directly from $\phi$.  
This is explained in some detail in  \cite{BGconstants}.
\end{remark} 
By a precisely analogous argument, we also have
\bp \label{A2Prop}
\[      \sum_{\mu \in S_0}|A_2(S_0,\mu)| \leq K_1 \n ,
 \] 
where $K_1$ is the constant defined prior to Proposition \ref{SummaryLemma}.
\end{proposition}

\section{The pleasingly rapid consumption of colours} \label{ConstantSection}

This section contains the  cancellation lemmas 
that we need to control the manner in which colours are consumed. 
The key result in this direction is the {\em Pincer Lemma} (Theorem \ref{PincerLemma}). 
   
\subsection{The Buffer Lemma}

\begin{lemma}\label{BufferLemma}  
Let $I\subset\bot(S)$ be an interval of edges labelled by constant letters, and 
suppose that the colours $\mu_1(S)$ and $\mu_2(S)$ lie either side 
of $I$, adjacent to it.  
 Provided that the whole of $I$ does not die in $S$,  no 
non-constant edge  coloured  $\mu_1$ will ever cancel with  
a non-constant edge coloured $\mu_2$. 
\end{lemma}  
 
\begin{proof} Suppose that the future of $I$ in $\top(S)$ 
is a non-empty interval labelled $w_0$. If $\mu_1(S)$ is to the left of $I$, 
then reading from the left beginning with the last non-constant 
edge coloured $\mu_1$, on the naive top of $S$ we have an interval labelled 
$x w_1 y$, 
where $y$ is a non-constant letter coloured $\mu_2$ and  
$w_1$ contains $w_0$ and perhaps some constant letters 
from $\mu_1$ and $\mu_2$.  
 
Our conditioning of $\phi$ (Proposition \ref{power}) ensures that, for all non-constant letters $z$, the rightmost non-constant letter in $\phi^j(z)$ is the same for all $j \geq 1$.  Therefore, in order for there to ever be cancellation between non-constant letters coloured $\mu_1$ and $\mu_2$, we must have $x = y^{-1}$.
Thus on $\top(S)$ there is an interval labelled  $xwx^{-1}$, where $w$ is the  
(non-empty) free-reduction 
of $w_1$. 
 
At times greater than $\height(S)$, the future of the interval that 
we are considering will continue to have a core subarc labelled $xw_jx^{-1}$,  
where $w_j$ is a conjugate of $w$ by a (possibly-empty) 
 word in constant letters (unless the interval hits a singularity or 
the boundary). In particular, no non-constant letters from $\mu_1$ 
and $\mu_2$ can ever cancel each other.  
\end{proof}

In the light of the Bounded Cancellation Lemma we deduce:

\begin{corollary} \label{BufferCorollary} 
Let $I\subset\bot(S)$ be an interval of edges labelled by constant letters, and 
suppose that the colours $\mu_1(S)$ and $\mu_2(S)$ lie either side 
of $I$, adjacent to it. 
If $|I|\ge B$ 
then there is never any cancellation between non-constant letters in $\mu_1$ and $\mu_2$. 
\end{corollary}

\subsection{The Two Colour Lemma} 
  
\bd \label{Neuters} 
Suppose that $U$ and $V$ are positive words\footnote{i.e. none of their letters are 
 inverses $a_j^{-1}$} 
 and that for some $k>0$ the only  negative exponents occurring   in $\phi^k(UV^{-1})$  
are on constant letters.  Then we say that $U$   
{\em $\phi$-neuters   $V^{-1}$ in at most $k$ steps}.  
\ed 
 
We shall also apply the term $\phi$-neuters to describe  the 
cancellation between  colours $\mu(S), \mu'(S) \subseteq \bot(S)$ that are adjacent in  corridors of van Kampen diagrams,  
and the following lemma remains valid in that context.

\begin{prop}[Two Colour Lemma] \label{TwoColourLemma} 
There exists a constant $\tz$ depending only on $\phi$ so that 
 for all positive words $U$ and $V$, if  
$U$ $\phi$-neuters  $V^{-1}$ then it does so in at most $\tz$ steps.  
\end{prop} 
 
\begin{proof} We express 
 $V^{-1}$ as a product of  three subwords:  reading from the left of
$V^{-1}$, the first subword ends with the last letter $y$ such that
$\phi(y)$ contains a left-fast letter; the second subword follows the
first and ends with the last non-constant letter in $V^{-1}$; the
remainder of $V^{-1}$ consists entirely of constant letters.
 
Lemma \ref{C_0} tells us that the length of the first subword is less
than $C_0$, and the proof of Lemma \ref{C1Lemma} provides a  bound of
$C_1$ on the length of the second subword.
 
Now consider the freely reduced form of $\phi^k(UV^{-1})$, and let $v_k$ denote its subword 
that begins with the first letter of negative exponent and ends with the final non-constant 
letter. The argument just applied to $V^{-1}$ shows that $v_k$ has length less than $C_0+C_1$ 
for all $k\ge 0$. 
 
Suppose that $U$ $\phi$-neuters $V^{-1}$ in exactly $N$ steps, let $\alpha_{N-1}$ 
be the letter of $\phi^{N-1}(UV^{-1})$ that consumes the last  letter of $v_{N-1}$, and 
let $\alpha_k$ be the ancestor of $\alpha_{N-1}$ in $\phi^k(UV^{-1})$. Write 
$\phi^k(UV^{-1}) = w_k\alpha_k u_k v_k w_k'$. 
 
Lemma \ref{C_0} shows that $|u_k| < C_0$ for all $k< N$, and we have just argued  
that $|v_k| < C_0+C_1$. Thus we obtain a bound  (independent of $U$ and $V$) on 
the number of words $\alpha_k u_k v_k$ that arise as $k$ varies --- call this 
number $\tz$. If $N$ were greater 
than $\tz$, then some configuration $\alpha_k u_k v_k$ with $v_k$ non-empty would recur. But 
this is nonsense, because once there is this repetition, the words $v_k$ will continue to repeat, and thus $V^{-1}$ will never be $\phi$-neutered, contrary to assumption.  
\end{proof} 

\begin{corollary} There exists a constant $\tz'$, depending only on $\phi$, with the following 
property: if $U$ and $V$ are positive words, $V$ begins with a non-constant letter 
and $\phi^k(UV^{-1})$ is positive for some $k>0$, then the least such $k$ 
is less than $\tz'$. 
\end{corollary}  
 
\begin{proof} The preceding lemma provides an upper bound on the least integer $N$ such that 
$\phi^N(UV^{-1})$ contains no non-constant  letters with negative exponent. Up to 
this point, the rightmost non-constant letter in $\phi^k(UV^{-1})$ may have been spawning 
constant letters to its right, and thus  $\phi^k(UV^{-1})$ may have a terminal segment consisting 
of constant letters. Since the rightmost non-constant letter of $\phi^k(V^{-1})$ does 
not vary with $k$ when $k<N$ (by Proposition \ref{power}),  the length of this segment  
grows at  a constant rate ($<M$) during each application of $\phi$. Similarly, its length 
changes at a constant rate after time $N$, decreasing until it is eventually cancelled. 
 
Since $N\le \tz$, this segment of constant letters has length less than $M\tz$ 
at time $N$, and hence is cancelled entirely before time $T_0(M+1)$. 
\end{proof}

\subsection{The disappearance of colours: Pincers and implosions} 
 
\begin{figure}[htbp] 
\begin{center} 

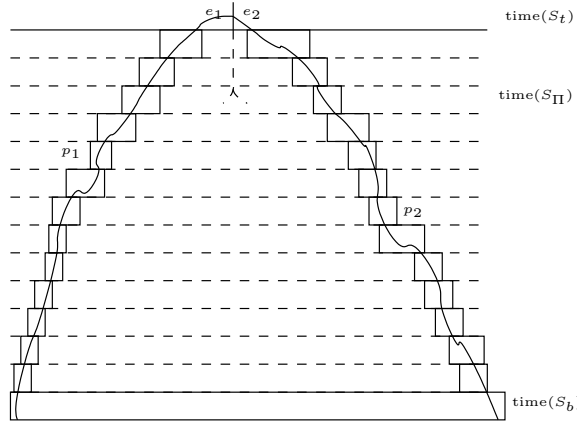 

\caption{A pincer.} 
\label{PincerPic} 
\end{center} 
\end{figure}

In this subsection we turn our attention to the detailed study of how non-adjacent colours 
along a corridor in $\Delta$ can come together solely as a result of the mutual 
annihilation of the intervening colours. Such an event determines a {\em pincer} (Figure \ref{PincerPic}), which is defined as follows.

\begin{definition}\label{pincerDef} 
Consider a pair of paths $p_1, p_2$ in $\mathcal F \subseteq \Delta$  
tracing the histories of $2$ non-constant edges $e_1, e_2$ that cancel in a corridor $S_t$.   
Let $\mu_i$ denote the colour of the 2-cells along  $p_i$. 
Suppose that at time $\tau_0$ these paths lie in a common corridor $S_b$.
Under these circumstances, we define 
the {\em pincer}  
$\Pin=\Pin( p_1, p_2, \tau_0 )$ to be the subdiagram of $\Delta$ enclosed by the chains of  
$2$-cells along $p_1$ and $p_2$, and the chain of $2$-cells connecting them in $S_b$. 
  
When it creates a desirable emphasis, we shall write $S_b(\Pin)$ and $S_t(\Pin)$ 
in place of $S_b$ and $S_t$. 
  
We define $S_\Pin$ to be the earliest corridor of the pincer in which $\mu_1(S_\Pin)$ and 
$\mu_2(S_\Pin)$ are adjacent. 
We define  $\tilde{\chi}(\Pin)$ to be the set of colours $\mu\notin\{\mu_1,\mu_2\}$ such 
that there is a 2-cell in $\Pin$ coloured $\mu$. And we define%\footnote{This
%definition has attracted significant social commentary \cite{pi}.}  
$$\life(\Pin) = \time(S_\Pi)-\time(S_b). 
$$ 
\end{definition} 
 
\begin{proposition}[Unnested Pincer Lemma]\label{prePincerLemma} 
There  exists a constant $\hat{T_1}$, depending only on $\phi$, such that for any pincer $\Pin$ 
\[      \life(\Pin) \leq \hat{T_1}(1+ |\tilde{\chi}(\Pin)|).        \] 
\end{proposition} 
 
Fix a pincer $\Pin$ and assume $\life(\Pin)\neq 0$.
 The idea of the proof of Proposition \ref{prePincerLemma} 
 is as follows:  we shall identify a constant $\hat{T_1}$ and argue 
 that if none of the colours $\mu\in\tilde{\chi}(\Pin)$ were consumed entirely  
  by  $\time(S_b) + \hat{T_1}$, 
the situation reached would be so stable  that no colours could be consumed in $\Pin$ at 
subsequent times,  
 contradicting the fact that all but $\mu_1$ and $\mu_2$ must be consumed by $\time(S_\Pi)$. 
 
With this approach in mind, we make the following definition: 
 
\begin{definition} Let $p$ be a positive integer.  
 A {\em $p$-implosive array} of colours  in a corridor $S$ 
is an ordered tuple $A(S)=[\nu_0(S),\dots,\nu_r(S)]$, with $r>1$, 
 such that: 
\begin{enumerate}
\item each  pair of colours $\{\nu_j, 
\nu_{j+1}\}$ is {\em essentially adjacent} in  $S$, meaning that there are no
non-constant edges of any other colour separating $\nu_j(S)$ from $\nu_{j+1}(S)$;
\item in each of the corridors $S=S^1,S^2,\dots, S^{p}$ in the future of $S$, 
 every $\nu_j(S^i)$ contains a non-constant edge;
\item in $S^p$,  {\em either} a non-constant edge coloured $\nu_0$ cancels a non-constant edge coloured $\nu_r$ 
(and hence the colours $\nu_j$ with $j=1,\dots,r-1$ are consumed 
entirely), {\em or else} all of the non-constant letters in $\nu_j(S^p)$, for $j=1,\dots, r-1$,
are cancelled in $S^p$ by edges from one of the colours of the array, while
$\nu_0(S^p)$ and $\nu_r(S^p)$ contain non-constant letters that survive in
the free-reduction of the naive future of the interval $\nu_0(S^p)\dots\nu_r(S^p)\subset\bot(S^p)$
(but  may nevertheless be cancelled in $S^p$ by edges from colours external to the array). 
\end{enumerate}
Arrays satisfying the first of the conditions in (3) are said to be of Type I, and those
satisfying the second condition are said to be of Type II. (These types are not
mutually exclusive.)

The {\em residual block} of an array of Type II is the interval of constant edges  between 
the rightmost non-constant letter of $\nu_0$ and 
the leftmost non-constant letter of $\nu_r$ in 
the free reduction of the naive future of  $\nu_0(S^p)\dots\nu_r(S^p)$.
The {\em enduring block} of the array  is the set of constant edges in 
$ \bot(S)$ that have a future  in the residual block.

 Note that there may exist {\em unnamed colours} between
$\nu_j(S)$ and $\nu_{j+1}(S)$ consisting entirely of constant edges.
\end{definition} 
 
\begin{remarks}\label{subarray} Let $[\nu_0(S),\dots,\nu_r(S)]$ be a $p$-implosive array.

\smallskip
(1) Any implosive subarray of $[\nu_0(S),\dots,\nu_r(S)]$ is $p$-implosive (same $p$).

(2) If an edge of $\nu_i$ cancels 
with an edge of $\nu_j$ and $j-i>1$, then this cancellation 
can only take place in $S^p$. If the edges cancelling are non-constant, 
 then the subarray $[\nu_i(S),\dots, \nu_j(S)]$ is   $p$-implosive of Type I.

(3) Given $x,y,w\in F$,
if the freely reduced words representing $x, y$ and $\phi(xwy)$ consist only 
of constant letters, then
so does the reduced form of $w$, since the subgroup generated by the constant
letters is invariant under $\phi^{\pm 1}$. It follows that the residual block
of any array of Type II contains edges from at most two of the colours $\nu_j$, and if
there are two colours they must be essentially adjacent, i.e. $\nu_j(S^p), \nu_{j+1}(S^p)$.

(4) For the same reason, the enduring block of an implosive
array of Type II is an interval involving at most two of the $\nu_j$, and if
there are two such colours then they must be essentially adjacent.
\iffalse
(5) The concept ``essentially adjacent" admits the possibility that there may be
several {\em unnamed colours} between $\mu_j(S)$ and $\mu_{j+1}(S)$, provided
that all of the edges in these unnamed colours are constant. Note that the
count defining the length of an array includes the edges in these unnamed colours.
\fi
\end{remarks}

\begin{lemma}\label{haveImp}  The ordered list of colours along each corridor before 
$\time(S_\Pi)$ in a pincer $\Pi$ must contain an implosive array.
\end{lemma}

\begin{proof} At the top of the pincer there is cancellation between non-constant
edges. Lemma \ref{BufferLemma} tells us that before $\time(S_\Pi)$ the colours
of these edges must have been separated by a non-constant letter of a different
colour, hence the list of non-constant colours along the bottom of $S_\Pi$ is a
1-implosive array. This same list of colours defines an implosive array at
each earlier time in the pincer until, going backwards in time, further non-constant
colours appear. Suppose $\mu$ has non-constant letters in $\Pi$ at
time $t$ but not time $t+1$. Let $\nu_0$ be  the first colour to the  left of $\mu$ that
contains non-constant letters at time $t+1$, and let $\nu_r$ be the first such colour
to the right. If $S_t$ is the corridor at time $t$, then the list of essentially-adjacent
non-constant colours $[\nu_0(S_t),\dots,\mu(S_t),\dots,\nu_r(S_t)]$ is a
1-implosive array. And $[\nu_0(S_{t'}),\dots,\mu(S_{t'}),\dots,\nu_r(S_{t'})]$
is a $(t'-t+1)$-implosive array for each earlier time $t'$ until (going backwards in
time) either further non-constant colours appear or else we reach the bottom of the
pincer.
\end{proof}

 If, further to the above lemma, we can argue that there is a constant $\hat{T_1}$ such
that each corridor before $\time(S_\Pi)$ contains a
$p$-implosive array with $p\le\hat{T_1}$, then we will know that at least one of the colours 
from $\tilde\chi(\P)$ is {\em essentially consumed} (i.e. comes to consist of constant
edges only)
 during each interval of $\hat{T_1}$ units in time during 
the lifetime of the pincer. Thus Proposition \ref{prePincerLemma} 
 is an immediate consequence of the following result, which will be proved
in (\ref{RIP}).

\begin{proposition}[Regular Implosions]\label{implosion} 
There is a constant $\hat{T_1}$ depending only on $\phi$ 
such that every implosive array in any minimal area diagram $\Delta$ is $p$-implosive 
for some $p\le\hat{T_1}$. 
\end{proposition} 
 
The first restriction to note concerning implosive arrays is this: 
 
\begin{lemma} \label{OnlyBColours} 
If $[\nu_0(S),\dots,\nu_r(S)]$ is implosive   %with $u$ unnamed colours.
of Type I, then $r\le B$. If
it is implosive of Type II, then $r< 2B$.
\end{lemma} 
 
\begin{proof} In Type I arrays,  the interval 
 $\nu_1(S^{p})\dots\nu_{r-1}(S^p)\subset\bot(S^p)$
  is to die in $S^p$, so $r-1<B$ by
  the Bounded Cancellation Lemma.
  For Type II arrays, one applies
the same argument to the intervals 
 joining $\nu_0(S^p)$ and $\nu_r(S^p)$ to the residual block of constant letters. 
\end{proof} 
 
\begin{remark} \label{shortisenough} 
In the light of Lemma \ref{OnlyBColours}, an 
obvious finiteness argument would provide the bound required for 
Lemma \ref{implosion} if we were willing 
to restrict ourselves to implosive arrays  with
 a uniform bound on their
length.
%the  length of the $\nu_j(S)$ and the length of the unnamed colours present in the array.
\iffalse
 Moreover, in the case of Type II arrays, we can exclude
the enduring block from the calculation of length, since it plays no role in any
cancellation.
\fi
Motivated by this observation, we seek to prove  that every implosive array contains an
implosive sub-array that is uniformly {\em short}.
\end{remark} 
 
\smallskip 
 
In order to identify a suitable notion of {\em short}, 
 we need to consider a further decomposition 
of the colours $\nu_j(S_b)$ in a $p$-implosive array $[\nu_0(S_b),\dots,\nu_r(S_b)]$. 
 
Previously (Subsection \ref{chromatic}) we partitioned each colour 
$\nu_j(S_b)$ into five intervals $A_1(S_b,\nu_j),\dots, 
A_5(S_b,\nu_j)$ and then further decomposed $A_4$ into subintervals 
$C_{(\nu_j,\nu')}(1)$ and $C_{(\nu_j,\nu')}(2)$ according to the 
colours of the edges that were going to consume these subintervals in 
the future. There is a corresponding decomposition of $A_2$ into 
intervals which we denote $C^2_{(\nu_j,\nu')}(1)$ and 
$C^2_{(\nu_j,\nu')}(2)$ (where $\nu'$ is now to the left of $\nu_j$ in 
$S_b$).  
 
\begin{figure}[htbp] 
\begin{center} 

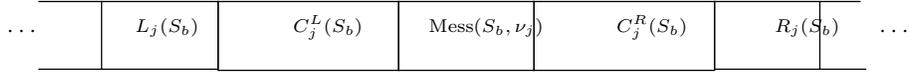 

\caption{The decomposition of the colour $\nu_j$} 
\label{PincerDecomp} 
\end{center} 
\end{figure}

Adapting to our new focus, we now define $R_j(S_b)=A_5(\nu_1,S_b)\cup 
C_{(\nu_j,\nu_{j+1})}(1)$, and $L_j(S_b)=A_1(\nu_1,S_b)\cup 
C^2_{(\nu_j,\nu_{j+1})}(1)$. We also define  $C_j^R(S_b)$ to be 
$C_{(\nu_j,\nu_{j-1})}(2)$ minus any edges from the excluded block, and $C_j^L(S_b)$ to 
be $C^2_{(\nu_j,\nu_{j-1})}(2)$ minus any edges from the excluded block.
Thus we obtain a decomposition of 
$\nu_j(S_b)$ into five intervals (see Figure \ref{PincerDecomp})  
$$ 
 L_j(S_b),\ C_j^L(S_b), \ \mess(S_b,\nu_j),\    C_j^R(S_b),\ R_j(S_b) 
$$ 
where $\mess(S_b,\nu_j)$ contains the edges 
whose preferred future dies at the time of  implosion together
with edges from the excluded block\footnote{At this point the reader may 
find it helpful to recall that only arrays of Type II have excluded
blocks, and such a block is either contained in a single colour,
or in adjacent colours $\nu_j(S_b)\cup\nu_{j+1}(S_b)$ with
the intervening intervals $R_j(S_b) \dots L_{j+1}(S_b)$ empty.}.

The terminal colours in our array, $\nu_0$ and $\nu_r$, play a special 
role. This is reflected in the fact that we shall only need to consider 
the segment of $\nu_0$ from its right end 
 up to and including the edge one to the left of $\mess(S_b,\nu_0)$. And 
 in $\nu_r$ we shall only need to consider the segment from its left end 
up to and including the edge one to the right of $\mess(S_b,\nu_r)$. 
 We write $\mathcal L(\nu_j,S_b)$ and $\mathcal R(\nu_j,S_b)$, 
  respectively, to denote these sub-intervals of $\nu_j(S_b)$. 
  
 \begin{definition} The length of $A(S)=[\nu_0(S),\dots,\nu_r(S)]$, written $\|A(S)\|$,
 is the number of edges in the interval
$\mathcal L(\nu_0,S)\dots\mathcal R(\nu_r,S)\subset \bot(S)$. (Note that $\|A(S)\|$ takes account of the
unnamed colours.)
\end{definition}

In keeping with the notation in the definition of $p$-implosive, we shall 
write $S^t$ for the corridor $t$ steps into the future of $S_b$; in particular  $S^0=S_b$ and
each $\nu_j$ with $j=1,\dots,r-1$ essentially vanishes in $S^p$. 
 
By definition, no preferred future of any 
 edge in  $\mess(\nu_j,S_b)$ is cancelled 
before $S^p$. Hence these intervals do not shrink in length before 
that time, and as in the proof of Lemma \ref{OnlyBColours} we can use
the Bounded Cancellation Lemma to bound the sum of their
lengths:

\begin{lemma} \label{2Bcols}
After excluding the edges of the enduring block, the sum of the lengths of the 
intervals  $\mess(\nu_j,S_b)$ is   at most $2B$. 
\end{lemma}

Combining this estimate with the bounds from Lemmas \ref{C_0} 
and \ref{C1Lemma}, we deduce that for $j=1,\dots,r-1$  
 
\[   |\nu_j(S_b)| \le    |C^L_j(S_b)| + |C^R_j(S_b)|  + 2C_0 + 2C_1 + 2B +\mathcal E_j,  \] 
where $\mathcal E_j$ is the number of edges from the excluded block coloured $\nu_j$.

Similarly, 
$$ 
|\L(\nu_0,S_b)|\le   
|C^R_0(S_b)| + C_0 + C_1 + B +\mathcal E_0  
$$ 
and  
$$ 
|\R(\nu_{r},S_b)|\le  
|C^L_{r}(S_b)| +C_0 + C_1 + B + \mathcal E_r. 
$$ 
This motivates us to define an array of colours $[\nu_0(S),\dots,\nu_r(S)]$ 
to be {\em very short} if for $j=1,\dots,r-1$ we have 
\[      |\nu_j(S)|  \leq 2C_0 + 2C_1 + 5B  + 1,    \]  
and 
\[      |\L(\nu_0,S)|  \leq C_0 + C_1 + 5B  + 1,    \] 
 and  \[      |\R(\nu_r,S)|  \leq C_0 + C_1 + 5B + 1,        \] 
and for $j=0,\dots,r-1$ the interval formed by the unnamed colours between $\nu_j(S)$
and $\nu_{j+1}(S)$ has total length at most $B$.

An implosive array is said to be {\em short} if it satisfies the weaker 
inequalities obtained by increasing  each of these bounds by $2B{\tz}$.   
 
\begin{lemma} \label{vshort} Let  $A=[\nu_0(S^0),\dots,\nu_r(S^0)]$ be  a 
$p$-implosive array 
with $p\ge \tz$. 
\begin{enumerate}
\item If $[\nu_0(S^\tz),\dots,\nu_r(S^\tz)]$ is very short, 
then $A$ is short. 
\item If $A$ is short, then $\|A\|\le 2B(2C_0+2C_1 + 5B +1+2BT_0) +2B^2(1+2T_0).$
\end{enumerate}
\end{lemma} 
 
\begin{proof} Item (1) is an immediate consequence of the Bounded Cancellation 
Lemma \ref{BCL}. The (crude) bound in (2) is an immediate consequence of Lemma \ref{2Bcols}
and the inequalities in the definition of {\em short}; the first summand is an estimate
on the sum of the lengths of the named colours, and the second summand accounts for
the unnamed colours.
\end{proof} 
   
The following 
lemma is the key step in  the proof of Proposition \ref{prePincerLemma}. 
 
\begin{lemma} \label{shorty} 
If $A(S^0)=[\nu_0(S^0),\dots,\nu_r(S^0)]$
 is a $p$-implosive array,
then at least one of the following statements is true:
\begin{enumerate}
\item $p\le 2T_0$; 
\item $A(S^0)$ is short;
\item  $p > 2T_0$ and $A(S^{T_0})$ contains an 
implosive sub-array $[\nu_k(S^\tz),\dots,\nu_l(S^\tz)]$ that is very short.
\end{enumerate}
\end{lemma} 
 
\begin{proof} Assume $p > 2T_0$ and that
 $[\nu_0(S^0),\dots,\nu_r(S^0)]$ is not short. We claim that there is a 
 block of at least $B+1$ constant letters in the interval determined by the
array
 $\mathcal L(\nu_0,S^{T_0})\dots\mathcal L(\nu_r,S^{T_0})$.
Indeed, by definition, if an array
 is not short then either one of the $\mathcal E_j$ has
length at least $B+1$, or one of the
 blocks of unnamed colours has length at
least $B(2T_0+1)+1$, or
else at least one of the intervals of 
constant letters $C^L_j(S^{0})$ 
or $C^R_j(S^{0})$ has length at least $B(T_0+1)+1$. 
In the first case, since  $\mathcal E_j$  is in
the excluded block, none of its edges are cancelled
 before the moment of implosion, and
hence it contributes a block of at least $ B+1$ constant
 letters to $A(S^{T_0})$; in the second
case, the Bounded Cancellation 
Lemma assures us that the length of the appropriate block of unnamed colours  can decrease by 
at most $2B$ at each step before the implosion of the array, 
and hence it still contributes
a block of  at least $ B+1$ constant edges to $A(S^{T_0})$;
and similarly, in the third case,
 $C^{\ast}_j(S^{0})$  can decrease by 
at most $B$ at each step before the implosion of the array.

Let $\beta$ be a block of at least $ B+1$ constant edges in $A(S^{T_0})$ with non-constant
edges $e_l$ and $e_{\rho}$ immediately to its left and right, respectively. 
\iffalse
The Two Colour Lemma 
\ref{TwocolourLemma} assures us that the letters labelling the non-constant edges adjacent
to the future of $\beta$ do not change between  $\time(S^{T_0})$ and $\time(S^p)$,
the moment of implosion. 
\fi
The Buffer Lemma \ref{BufferLemma} assures us that the
non-constant edges in the future of $e_l$ will never interact with the non-constant edges in the
future of $e_{\rho}$. Thus at least one of $e_l$ or $e_{\rho}$ must be {\em stabbed in the back}, i.e. 
its entire non-constant future must be consumed by edges on its own side of $\beta$. Suppose,
for ease of notation, that it is $e_l$ and let $\nu_i$ be the colour of $e_l$. We claim that if $\nu_k$
is the colour of the letter that ultimately consumes $e_l$, then $k\le i-2$.

We shall derive a contradiction
from the assumption that the edge which ultimately
 consumes  $e_l$ is  coloured $\nu_{i-1}$.
There are two cases to
consider according to whether $e_{\rho}$ is also coloured $\nu_i$. If it is, then we consider the
word $V$ labelling the arc of $\bot(S^0)$
 from the left end of $\nu_i(S^0)$ to the past of
$e_l$; the consumption
of the non-constant future of 
$e_l$ completes the $\phi$-neutering  of $V$
 by the word labelling $\nu_{i-1}(S^0)$, 
in particular this neutering will have taken more than $T_0$ steps in time, contradicting
the Two Colour Lemma \ref{TwoColourLemma}. If $e_{\rho}$ is not coloured $\nu_i$, then the
consumption of the non-constant future of 
$e_l$ results in a new essential adjacency of colours and hence can only be complete
 at the moment of
implosion, i.e. $\time(S^p)$. But this consumption constitutes the neutering of $\nu_i(S^{T_0})$
by $\nu_{i-1}(S^{T_0})$, and according to the Two Colour Lemma this neutering 
must be accomplished in at most $T_0$ units of time. Thus $p\le 2T_0$,
contrary to our hypothesis. 

\def\kill{\!\!\!\searrow\!}

Thus we have proved that the edge which ultimately consumes $e_l$ is coloured
$\nu_k$ where $k\le i-2$. Under these circumstances (or the 
symmetric situation with $e_{\rho}$ in place of $e_l$) we say that {\em $\nu_k$
neuters $\nu_i$ from behind} and write  $\nu_k \kill \nu_i$.
 
%***  We draw a `cancellation diagram' FIGURE to SHOW WHAT HAS HAPPENED... *** 
 
\begin{figure}[htbp] 
\begin{center} 
  
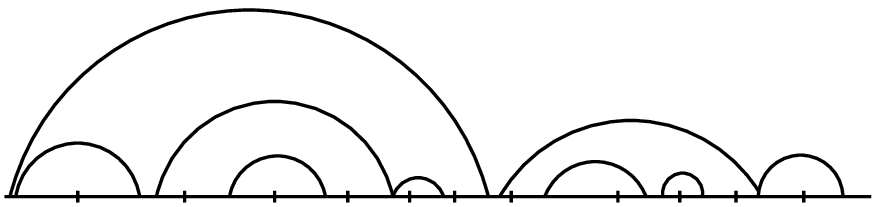 
  
\caption{The nesting associated to $\kill$} 
\label{figure:Arc} 
\end{center} 
\end{figure} 
 
\medskip 
 
There is a natural {\em nesting} among the $\kill$-related pairs of colours from the
array: 
% if $\{\nu_{k_1},\nu_{i_1}\}$ and $\{\nu_{k_2},\nu_{i_2}\}$ are two 
% such pairs, 
 $(\nu_{k_1},\nu_{j_1}) < (\nu_{k_2},\nu_{j_2})$ if $\nu_{k_1}$ and 
$\nu_{j_1}$   
both lie between $\nu_{k_2}$ and $\nu_{j_2}$ in $S^0$. See Figure 
\ref{figure:Arc}. 
 
We focus our attention on an innermost (i.e. minimal) 
pair with  $\nu_k\kill\nu_i$.  By definition $|k-i|\ge 2$. If there were
a block of at least $ B+1$ constant letters between the closest non-constant
letters of $\nu_k(S^{T_0})$ and $\nu_i(S^{T_0})$, 
then the preceding argument 
would yield a neutering from behind that contradicted the innermost nature
of $\nu_k\kill\nu_i$. Thus  
$[\nu_j(S^\tz),\dots,\nu_k(S^\tz)]$ is a very short array, and we are done.  
\end{proof}

\begin{RIproof}\label{RIP}{\em Proof of Regular Implosions (Prop.\ref{implosion}):} Given the bound in Lemma \ref{vshort}(2), an obvious finiteness
argument provides a constant $\tau$ such that every short implosive array is $p$-implosive
with $p\le \tau$. And the same bound applies to implosive arrays that contain a short
sub-array (Remark \ref{subarray}(1)). So in the light of Lemmas  \ref{shorty} and   \ref{vshort}(1),
it suffices to let $\hat{T_1} = \max\{2T_0, \tau\}$.\hfill $\square$
\end{RIproof}
 
\subsection{Super-Buffers}

In this subsection we prove an important cancellation lemma based on
Proposition \ref{prePincerLemma}, this lemma involves the following
constant.

\begin{definition} \label{T1'Lemma}
We fix an integer $T_1'$ such  that one gets repetitions  in all $T_1'$-long subsequences of $5$-tuples of reduced words
\[	U_k:=\Big( u_{k,1}, u_{k,2}, u_{k,3}, u_{k,4}, u_{k,5}	 \Big) \ \ \ \ k=1,2,\dots \]
with $|u_{k,1}| $ and $|u_{k,1}| $ at most $ C_0+ C_1 + 2B +1$, while
 $|u_2^k|$ and $ |u_4^k|$ are at most $ C_0 + C_1$, and $|u_3^k| \leq 4B+1$. That is,
for some $t_1\le t_2\leq T_1'$ and 
\[	\Big( u_{t_1,1}, u_{t_1,2}, u_{t_1,3}, u_{t_1,4}, u_{t_1,5}	 \Big) = \Big( u_{t_2,1},  u_{t_2,2}, u_{t_2,3}, u_{t_2,4}, u_{t_2,5}	 \Big) .	\]
\end{definition}

\newtheorem{stipulation}[theorem]{Stipulation}

\begin{stipulation} Assume $T_1' \ge \hat{T_1}$.
\end{stipulation}
 
The cancellation lemma we need is most easily phrased in terms of 
colours of subwords, which we define as follows, keeping firmly in mind
the example of a stack of partial corridors excised from the interior of a van Kampen
diagram, retaining their memory of the colours to which the edges belong.

We have a word $W$ with a  decomposition into preferred subwords
$V = V_1 V_2
\cdots V_k$, where each $V_i$ is either positive or negative;
we think of these subwords as having colours $\mu_1, \ldots \mu_k$.
Take the freely reduced words $\phi(V_i)$, concatenate them, then
cancel to form a freely reduced word. There is some freedom in the
choice of cancellation scheme, as in the folding of corridors, but we fix
a choice, thus assigning to each letter of the freely reduced form of
$\phi(V)$ the colour  $\mu_i$ of its ancestor. We repeat this process,
thus assigning colours to the letters in the reduced form of $\phi^k(V)$ for
each integer $k>0$.

The process that we have just described is an algebraic description of
a choice of  minimal area van Kampen diagram for $t^{-k}Vt^k\phi^k(V)^{-1}$.
Thus the following lemma is a comment on the form of  such
diagrams.

\begin{proposition} \label{NoDoubleNeuter}  
Let $V=V_1V_2V_3$ be a concatentation of words (coloured $\nu_1, \nu_2, \nu_3$)
each of which is either positive or negative.
If $W$ is a subword of the reduced
form of  $\phi^{T_1'}(V)$ and   $W$ has a non-constant
letter coloured $\nu_i$ for each $i\in\{1,2,3\}$, 
then for all $k \geq 0$ there are
non-constant letters  in $\phi^k(W)$ coloured $\nu_2$.
\end{proposition}

\begin{proof} Let $\nu_i(W)$ denote the subword of $W$
coloured $\nu_i$, and let $\nu_i^j$ denote  the maximal subword coloured $\nu_i$ in
(the reduced word representing) $\phi^i(V_1V_2V_3)$ .
Note that $\nu_2(W)=\nu_2^{T_1'}$, and more generally $\nu_2^{T_1'+j}$
is the maximal word  in $\phi^j(W)$ coloured $\nu_2$.

Fix $k>T_1'$ and consider the diagram formed by the stack of corridors
described prior to the proposition. The bottom
of the first corridor is labelled $V$, and we regard it as being divided into
three coloured intervals according to the decomposition $V_1V_2V_3$.
Since $\nu_2(W)$ contains non-constant letters and $T_1'>\hat{T_1}$, 
the array formed by these colours is not implosive (Proposition \ref{prePincerLemma}),
and hence  
$\nu_1(W)$ and $\nu_3(W)$ will never essentially consume $\nu_2(W)$.
However, the proposition is not yet proved because there remains
the possibility
that  $\nu_2$ may essentially vanish because it  neuters $\nu_1(W)$, say, and is then
neutered by $\nu_3(W)$. We proceed under this assumption, seeking
a contradiction. (The case where the roles of 
$\nu_1$ and $\nu_3$ are reversed is entirely similar.)

For each $1 \leq i \leq T_1'$, we have
$\phi^i(V_1V_2V_3)=\nu_1^i, \nu_2^i$ and $\nu_3^i$.
Write $\nu_2^i\equiv V^i(1)  V^i(2)  
V^i(3)$, where $V^i(1)$ ends with last  non-constant letter in
$\nu_2^i$ whose entire non-constant future is eventually consumed by
letters coloured $\nu_1$, and $V^i(3)$ begins with the leftmost
non-constant letter whose entire non-constant future is  
eventually consumed by letters coloured $\nu_3$.
Lemmas \ref{C_0} and \ref{C1Lemma} tell us that $V^i(1)$ and $V^i(3)$ have
length at most
$C_0 + C_1$.  

\noindent{\em Claim:} $V^i(2)$  contains exactly one non-constant edge
and has length no more than $4B+1$.

We are assuming that $\nu_2(W)$ neuters $\nu_1(W)$. Consider the
(non-constant) edge $\e_i$ in
$\nu_2^i$ that will eventually consume the final non-constant edge in
$\nu_1(W)$. Note that $\e_i$ is  the
leftmost non-constant edge in $V^i(2)$. Moreover, we are assuming 
that  $\nu_3(W)$ ultimately neuters
$\nu_2(W)$, so in particular  it consumes the entire future of
any edge to the right of $\e_i$, which
forces $\e_i$ to be the rightmost non-constant edge in $V^i(2)$. The Buffer
Lemma tells us that $\e_i$ must lie within $2B$ of both ends of $V^i(2)$,
and hence the claim is proved.

Looking to the left
of $V^i(1)$, we now consider the  subword    $L^i$ of $\nu_1^i$  
that begins with the leftmost  non-constant edge in the future of
which there is a non-constant letter that cancels with a letter
coloured $\nu_2$. And looking to the right of  $V^i(3)$, we consider
the subword that ends with the rightmost non-constant letter  in the future of
which there is a non-constant letter that cancels with a letter
coloured $\nu_2$. 
any of whose non-constant future cancels
with an edge painted $\nu_2$.  As in previous arguments, The Buffer Lemma and
Lemmas \ref{C_0}, \ref{C1Lemma} tell is that  $|R^i|, |L^i| \le C_0 +
C_1 + 2B +1$, for all $i$.

We have already bounded the lengths of $V^i(1), V^i(2)$ and $V^i(3)$
by $C_0+C_1, 4B+1$ and $C_0+C_1$, respectively. Thus we are
now in a position to invoke the repetitive behaviour described in Definition
\ref{T1'Lemma}:
for some positive integers $i $ and $t$ with $i+t\le T_1'$, we get a repetition 
\[	\Big( R^i, V^i(1), V^i(2), V^i(3), L^i \Big) =  
\Big( R^{i+t}, V^{i+t}(1), V^{i+t}(2), V^{i+t}(3), L^{i+t} \Big).
\]
For as long as we are assured of the continuing presence
of $\nu_1^{i+s}$ and $\nu_3^{i+s}$,
 the fate of $\nu_2^i=V^i(1)V^i(2)V^i(3)$ under $s$ iterations of $\phi$ depends
only on $(R^i, V^i(1), V^i(2), V^i(3), L^i)$. Thus
$$
\Big( V^j(1), V^j(2), V^j(3) \Big) = \Big( V^{j+t}(1),
V^{j+t}(2), V^{j+t}(3) \Big)
$$
for all $j\ge i$ within the time scale of this assurance. However this leads us
to an absurd conclusion, because once $ \nu_1$ has become constant, 
at all subsequent time,
the surviving word coloured  $\nu_2$ contains as a proper subword, the
$\nu_2$ word that existed at the corresponding times in the  cycles (of
period $t$) before $T_1'$, and in particular they can never essentially
vanish, contrary to our assumption that $\nu_3$ eventually neuters
$\nu_2$. 
\end{proof}

\subsection{Nesting and the Pincer Lemma}

In subsequent sections we would like to
 bound the life of pincers by arguing
that during the lifetime of a pincer,
 colours must be consumed at a predictable rate (appealing
to Proposition \ref{prePincerLemma}),
noting that there  are only a limited number of colours. However, the bounds
we need will require us to ascribe each consumed colour to a {\em unique}
pincer. Thus we encounter problems whenever one pincer is contained in another.
For reasons that will become apparent in subsequent sections,
%(and  were foreshadowed in the proof of Lemma \ref{shorty})
in situations where
we must confront this problem, the inner of the two pincers will have a long block
of constant edges along the corridor immediately above its peak. More precisely,
we will find ourselves in the situation described in the following definition. The
appearance of the constant $\ll := 2B(T_0+1)+1$ in the following definition is
explained by the role that this constant played in the course of Lemma \ref{shorty}.

\begin{figure}[htbp] 
\begin{center} 
  
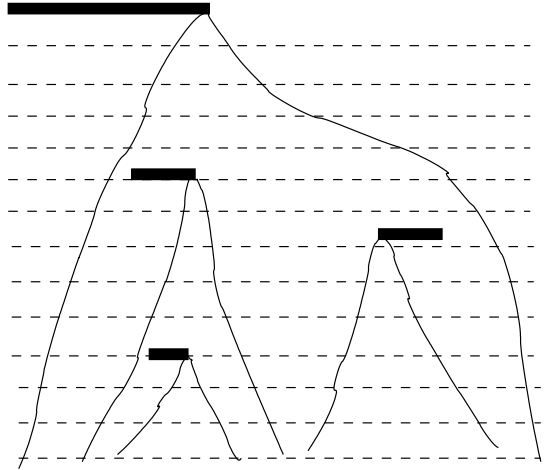 
  
\caption{A depiction of nesting} 
\label{figure:Nest} 
\end{center} 
\end{figure}

\begin{definition} \label{NestingDef}
Consider one pincer $\Pin_1$ contained in another $\Pin_0$.  Suppose
that in the corridor $S \subseteq \Pin_0$ at the top of $\Pin_1$
(where its boundary paths $p_1(\Pin_1)$ and $p_2(\Pin_1)$ come
together) the future in $\top(S)$ of at least one of the edges
containing $p_1(\Pin_1) \cap \bot(S)$ or $p_2(\Pin_1) \cap \bot(S)$
contains no non-constant edges, and this future\footnote{We allow this
future to be empty, in which case ``contained in" means that the immediate past
of the long block of constant edges is not separated from $\Pin_1$ by any
edge that has a future in $\top(S)$.}
lies in an interval of at least $\ll$ constant edges contained in
$\Pin_0$.  Then we say that $\Pin_1$ is {\em nested in} $\Pin_0$. (in
Figure \ref{figure:Nest}, the $\ll$-long block of constant edges are shown in
black.) We say that $\Pin_1$ is {\em left-loaded} or {\em right-loaded}
according to the direction in which  the $\ll$-long block of constant edges
extends from the peak of $\Pin_1$.
\end{definition} 

\begin{remark}  A nested pincer cannot be both left-loaded and right-loaded (cf.
Remark \ref{subarray}(3)).

If $\Pin_1$ is left-loaded, then  the future of  
$p_1(\Pin_1) \cap \bot(S)$ contains no non-constant edges.
It may happen that the future of $p_2(\Pin_1)$ also contains no 
non-constant edges; in this case the colour $\mu$ of $p_2(\Pin_1)$
essentially vanishes in $S$ due to   cancellation between non-constant edges 
of $\mu$ and some colour to its right. Symmetric considerations apply
to right-loaded pincers.
\end{remark}

\begin{definition} \label{chiP}
For a pincer $\Pin_0$, let $\{ \Pin_i \}_{i \in I}$ be the set of all
pincers nested in $\Pin_0$.  Then define
\[	\chi(\Pin_0) = \tilde{\chi}(\Pin_0) \smallsetminus \bigcup_{i \in
I} \tilde{\chi}(\Pin_i).	\]
\ed

\bl \label{nestLife}
If the pincer $\Pin_1$ is nested in $\Pin_0$ then 
$\time(S_t(\Pin_1)) < \time(S_{\Pin_0}).	$
\end{lemma}
\begin{proof}
The presence of the hypothesised block of constant letters in
$\top(S_t(\Pin_1))$ makes this an immediate consequence of the Buffer
Lemma \ref{BufferLemma}.
\end{proof}

Define $T_1 :=  T_1' + 2T_0$.
The following theorem is the main result of this section.

\begin{theorem}[Pincer Lemma] \label{PincerLemma} 
For any pincer $\Pin$
\[	\life(\Pin) \leq T_1(1 + |\chi(\Pin)|).	\]
\end{theorem}

\begin{proof}
The heart of our proof of Proposition \ref{prePincerLemma} was that
in each block of $\hat T_1$  steps in time between  $\time(S_b)$ and
$\time(S_\Pin)$ at least one colour essentially disappears. Our proof
of the present theorem is an elaboration of that argument: we must
argue for the essential disappearance of a  colour that is not
contained in any of pincers nested in $\Pin$.  Thus we concentrate
on that region of the pincer $\Pin$ that is exterior to the set of 
{\em co-level\footnote{i.e. those that are maximal
with respect to inclusion among the pincers nested in $\Pin$} 1} pincers nested in it;
let $\{ \Pin_j\}, \ j=1,\dots,{J}$ be the set of such, indexed in order of appearance from
left to right.

For $j=1,\dots,J-1$, let $\Sigma_j$ denote the set of  colours  along the bottom of $\Pin$
that have a non-constant edge strictly
between $\Pin_j$ and $\Pin_{j+1}$; if $\Pin_j$ is left-loaded, then we include
the colour of $p_2(\Pin_j)$ in $\Sigma_j$, and if $\Pin_j$ is right-loaded, then we include
the colour of $p_1(\Pin_j)$ in $\Sigma_{j-1}$. Likewise, we define $\Sigma_0$ to be
the set of non-constant colours that lie to the left of $\Pin_1$ together with the
colour of $p_1(\Pin)$, and we define
$\Sigma_{J}$ to be
the set of non-constant colours that lie to the right of $\Pin_J$ together with the
colour of $p_2(\Pin)$. 

In order to prove the theorem, we derive a contradiction from the assumption that
in the first $T_1$ units of time in the life of $\Pin$ no colours in the union of the
$\Sigma_j$ essentially vanish. (There is no loss of generality in starting at the
bottom of the pincer, since given any other starting time, one can discard the
pincer below that level.) We label the corridors, beginning at the bottom of $\Pin$
and proceeding in time as $S^0,S^1,\dots$

We focus on a single $\Sigma_j$, and write its colours in order as    $\nu_1, \ldots , \nu_r$.
We analyse how the colours in $\Sigma_j$ come to vanish.  
The first important observation is that  $2 \le i \le r-1$,
 it is not possible for the colour $\nu_i$ to essentially vanish (at any time)
due to cancellation merely between the colours in $\Sigma_j$.  
For if this happened,  there would be an implosive array in $S^0$
containing $\nu_i(S^0)$  and so, by Proposition \ref{prePincerLemma}, $\nu_i$ would vanish before $S^{T_1}$, contrary to our assumption.

There remains the possibility that $\nu_2$ may neuter  $\nu_1$ (after
$S^{T_1}$).  This can happen in two ways.  The first is that $\Pin_{j-1}$ is left-loaded: in
this case
the neutering happens within time $T_0$ of the top of $\Pin_{j-1}$ (by Two Colour Lemma),
and we are then in a stable situation in the sense that $\nu_3$ cannot subsequently neuter $\nu_2$,
by Proposition \ref{NoDoubleNeuter}.  Now suppose that $\Pin_{j-1}$ is right-loaded.
Consider the earliest time $t_0$ at which there is a block of  at least $B+1$ constant edges in the
past of the $\lambda_0$-long block associated to $\Pin_{j-1}$. If $\nu_2$ is to neuter
$\nu_1$, then it must do so within   $T_0$ steps of this time. Indeed, within $T_0$ steps,
if the non-constant edges of $\nu_1$ to the right of the block have not been consumed
by $\nu_2$, 
then they will never be consumed by a colour from $\Sigma_j$. 

There is a further event that we must account for, which is closely related to
neutering: it may
happen that $\nu_1$ is the colour of $p_2(\Pin_{j-1})$ and that $\nu_2$ consumes
all of the non-constant edges to the right of the block of constant edges discussed above;
this is not a neutering but nevertheless the Two Colour Lemma applies. We would like
to apply Proposition \ref{NoDoubleNeuter} in this situation to conclude that
$\nu_3$ cannot subsequently neuter $\nu_2$.
This is legitimate provided
$t_0\ge\time (S^{T_1'})$.  If $t_0< \time (S^{T_1'})$, then we still know that $\nu_3$
cannot  neuter $\nu_2$ before $S^{T_1}$, because by hypothesis no colour from
$\Sigma_j$ essentially vanishes before this time. On the other hand, the Two Colour Lemma
tells us that if $\nu_3$ is to neuter $\nu_2$, then it must do so within $T_0$ steps
from $t_0$, and $t_0+T_0\le  \time (S^{T_1})$. Thus, once again, we conclude that
$\nu_3$ can never neuter $\nu_2$.

Entirely similar arguments show that it cannot happen that $\nu_r$ is neutered
by $\nu_{r-1}$ and that subsequently $\nu_{r-2}$ neuters $\nu_{r-1}$.

We have established the existence of a stable situation: proceeding past the point where
the restricted amount of possible neutering within $\Sigma_j$
has occurred, we may assume that the next
essential disappearance of a colour from $\Sigma_j$ can only occur as a result of
cancellation with a colour from some $\Sigma_i$ with $i\neq j$. Such further cancellation
must occur, of course, because all but two\footnote{Degenerate cases with few
colours are covered by the Two Colour Lemma and the Buffer Lemma.} of the
 colours in $\bigcup_j\Sigma_j$ must be consumed within $\Pin$.

Passing to innermost pair of interacting $\Sigma_k$
we may assume $i=j-1$ (cf. proof of Lemma \ref{shorty}). Thus our proof will be
complete if we can argue that cancellation between non-constant edges
from $\Sigma_{j-1}$ and $\Sigma_j$ is impossible. We have
argued that the colours which are to cancel will be essentially adjacent within
time $T_0$ of 
the top of $\Pin_{j-1}$. On the other hand, there is a block of $\ll$ constant
edges separating $\Sigma_{j-1}$-nonconstant edges and $\Sigma_{j}$-nonconstant edges
at the top of $\Pin_{j-1}$. Since $\ll > 2B(T_0+1)$ at least $B+1$ of these constant edges
remain $T_0$ steps later. The Buffer Lemma now obstructs the supposed
cancellation between non-constant edges in $\Sigma_{j-1}$ and $\Sigma_j$. 
\end{proof}

\section{Teams and their Associates}\label{teamSec} 
 
We begin the process of grouping pairs of colours $(\mu,\mu')$ into 
teams.  
 
\subsection{Pre-teams} \label{t1} 

The whole of $C_{(\mu,\mu')}(2)$ will ultimately be  consumed by  
a single edge $\e_0\in\mu'(S_0)$. 
We consider the time $t_0$ at which the future of $\e_0$ 
starts consuming the future of  $C_{(\mu,\mu')}(2)$.  
If $|C_{(\mu,\mu')}(2)|> 2B$, then this consumption will not be completed in 
three steps of time  (Lemma \ref{BCL}). We claim that in  this circumstance, the 
leftmost $\mu'$-coloured edge after the first two steps of the
cancellation must be left para-linear. Indeed it is not left-constant since it must consume edges in 
the future of $C_{(\mu,\mu')}(2)$, and since no non-constant $\mu'$-edges 
are cancelled by $\mu$ in passing from the first to the second stage of  
cancellation, the leftmost non-constant $\mu'$-label must remain the same (Proposition 
\ref{power}). We denote this \lpl edge at time $t_0+2$ by $\e^\mu$. 
 
Let $\e_\mu$ be the rightmost edge in the future of $C_{(\mu,\mu')}(2)$  at time $t_0$. 
We trace the ancestry of $\e_\mu$ and $\e^\mu$ in the trees of $\F\subset\Delta$ corresponding 
to the colours $\mu$ and $\mu'$ (as defined in \ref{tree}). 
 We go back to the last point in time 
$\ptmm$ at which  
both ancestors  lay in a common corridor 
 {\em and} the interval on the bottom of this corridor between the pasts of 
  $\e_{\mu}$ and $\e^\mu$  
is comprised entirely of constant edges whose future is eventually  
 consumed by the ancestor of 
$\e^\mu$ at this time. We denote this corridor  $S_{\uparrow}$. 
 
\begin{definition}\label{preteam}  The ancestor of $\e^\mu$ at time $\ptmm$ is called the 
{\em reaper} and is denoted $\prmm$.  The set of edges in $\bot(S_{\uparrow})$ 
 which are eventually consumed by $\prmm$ is denoted $\pEmm$.  
  This is a contiguous set of edges. 
The {\em pre-team} $\pTmm$ is defined to 
be the set of pairs $(\mu_1,\mu')$ such that $\pEmm$ contains 
edges coloured $\mu_1$. The number of edges in $\pEmm$ is denoted 
$\|\hat \T\|$. 
\end{definition} 
 
In a little while we shall define {\em teams} to be   
pre-teams satisfying a certain maximality condition (see Definition \ref{newTeams}).

\begin{remark} If $\ptmm<\time(S_0)$ then 
near the right-hand end of $\pEmm$ one may have  an interval of 
colours  $\nu$ such that $\nu(S_0)$ is empty.  
\end{remark} 
 
\smallskip 
 
In the proof of Proposition \ref{SummaryLemma} we saw that it would be
desirable if (whatever our final 
definition of {\em team} and $\bonus$ may be) the following inequality
(\ref{goodEq}) should hold for all teams: 
\begin{equation} \label{preTeamInequality} 
\sum_{(\mu,\mu') \in \T \mbox{ \tiny or } (\mu,\mu') \vin \T} |C_{(\mu,\mu')}(2)| \le 
\|\T\|  + |\bonus(\T)| + B. 
\end{equation}

The following lemma shows that, even without introducing a  bonus scheme 
or virtual members, the 
desired inequality is straightforward for pre-teams with $\ptmm \geq \time(S_0)$. 
 
\begin{lemma} \label{t1high} 
 If $\ptmm \geq \height(S_0)$ then $\pTmm$ satisfies 
\[ \sum_{(\mu,\mu') \in \pTmm}|\cmm| \leq \|\pTmm\| + B.      \] 
\end{lemma} 
 
\begin{proof} By definition   
$\mu'(S_0)$ does not start consuming any  of the
 $C_{(\mu_1,\mu')}(2)$ with $(\mu_1,\mu')\in\pT$ before $\ptmm$  
(apart from a possible nibbling of  length $< B$ from the rightmost team 
member  at time $\ptmm -1$).  
Since each $C_{(\mu_1,\mu')}(2)$ consists only of edges consumed 
by $\mu'(S_0)$, the future  of each $C_{(\mu_1,\mu')}(2)$ at time $\ptmm$ 
will have the same length as $C_{(\mu,\mu')}(2)$ 
(except that the rightmost may have lost these $< B$ edges).  
And these futures are contained in $\pEmm$. 
\end{proof}

The case where $\ptmm <\height(S_0)$ is more troublesome. As $\pEmm$ flows forwards in 
time, the number of constant letters in the future of $\pEmm$ that are  
consumed by $\prmm$ between $\ptmm$ and $\time (S_0)$ may be outweighed by the number of 
constant letters generated to the left of the future of $\pEmm$ that will 
ultimately be consumed by $\prmm$. 
 
It is to circumvent the failure of inequality 
(\ref{preTeamInequality}) in this setting that we are 
obliged to instigate the bonus scheme described in Section \ref{BonusScheme}.

\subsection{The Genesis of pre-teams} \label{genesis} 
 
We fix $\pTmm$ with $\ptmm < \time(S_0)$ and consider 
the various events that occur  at  $\ptmm$ to prevent 
us pushing the pre-team back one step in time.  We write $S_\omega$ 
to denote the corridor at time $\ptmm$ containing $\pEmm$. 
 
\begin{figure}[htbp] 
\begin{center} 
  
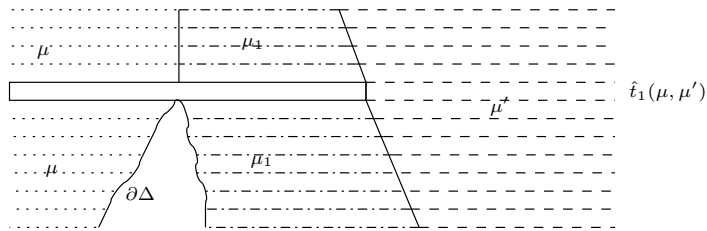 
  
\caption{A team of genesis (G1)}
\label{G1Pic} 
\end{center} 
\end{figure}

\begin{figure}[htbp] 
\begin{center} 
  
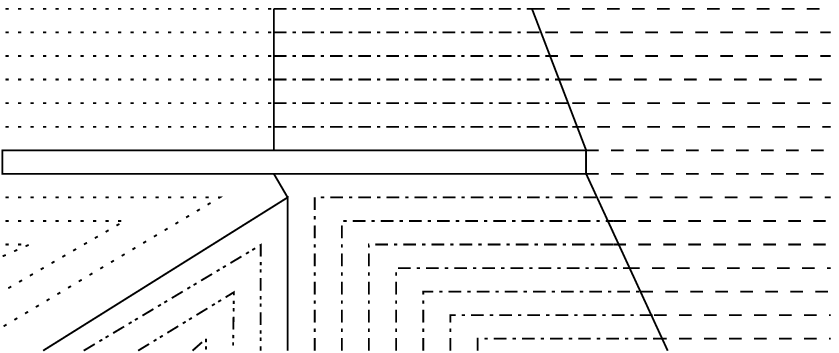 
  
\caption{A team of genesis (G2)}
\label{G2Pic} 
\end{center} 
\end{figure} 

\begin{figure}[htbp] 
\begin{center} 
  
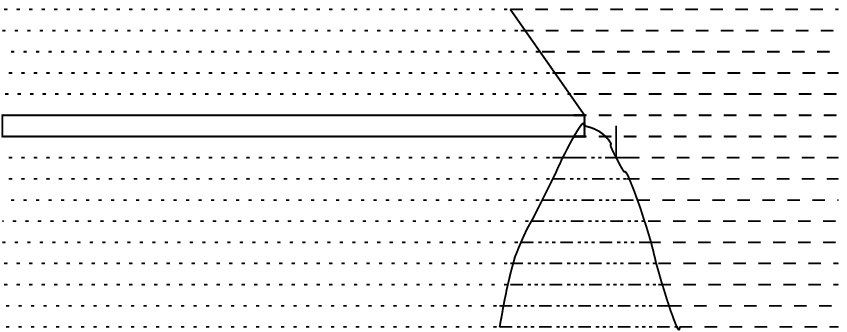 
  
\caption{A team of genesis (G3)}
\label{G3Pic} 
\end{center} 
\end{figure} 

\begin{figure}[htbp] 
\begin{center} 
  
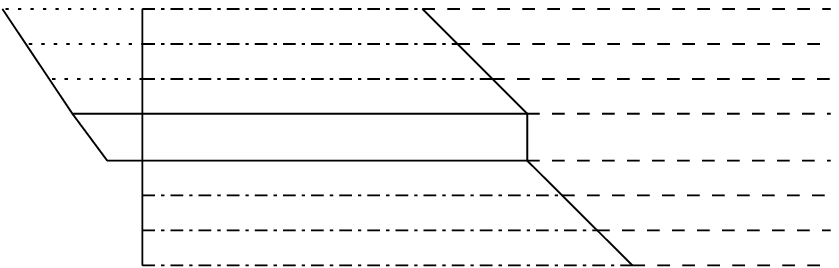 
  
\caption{A team of genesis (G4)}
\label{G4Pic} 
\end{center} 
\end{figure}

\def\CC{C}
 
There are four types of events: 
 
\begin{enumerate} 
\item[(G1)] The immediate past  of $\CC_{(\mu,\mu')}(S_\omega) 
$ is separated from the 
past of $\prmm$ by an intrusion of $\partial\Delta$ (Figure \ref{G1Pic}). 
\item[(G2)] We are not in case (G1), but the immediate past  of
$\CC_{(\mu,\mu')}(S_\omega)$ is separated from the 
past of $\prmm$ because of a singularity (Figure \ref{G2Pic}). 
\item[(G3)] The immediate past of $C_{(\mu,\mu')}(S_\omega)$ is  still
in the same corridor as the past of $\prmm$, but it is  separated from
it by a non-constant letter (Figure \ref{G3Pic}). 
\item[(G4)] We are not in any of the above cases,  
but the immediate past of the rightmost letter in  
 $C_{(\mu,\mu')}(S_\omega)$ is not constant (Figure \ref{G4Pic}). 
\end{enumerate} 
 
\smallskip  
 
\def\ST{S_\T} 
\def\STmm{\S_{\T(\mu,\mu')}}

The following lemma explains why Figures \ref{G3Pic} and \ref{G4Pic} are an
accurate portrayal  of cases (G3) and (G4). 
 
Let $M_{inv}$ be the maximum length of $\phi^{-1}(x)$ over generators
$x$ of $F$, and $C_4 = M_{inv}.M$.  
 
\begin{lemma}\label{G34pics} 
If $I$ is an interval on $\top (S)$ labelled by a word $w$ in constant
letters
then the reduced word labelling the past of $I$ in $\bot(S)$ is of the 
form $u\alpha v$, where $\alpha$ is a word in constant letters and 
$|u|$ and $|v|$ are less than $C_4$. Moreover, if the past of the leftmost  
(resp. rightmost) letter 
in $w$ is constant, then $u$ (resp. $v$)  is empty. 

In particular, $|I| \leq |\alpha| + 2MC_4$.
\end{lemma}

\begin{figure}[htbp] 
\begin{center} 
  
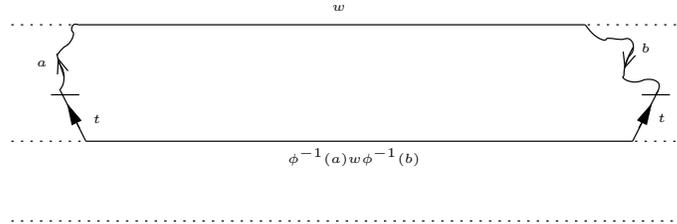 
  
\caption{The proof of Lemma \ref{G34pics}}
\label{bottom-constant} 
\end{center} 
\end{figure} 

\begin{proof} See Figure \ref{bottom-constant}.  Follow the path from
the left end of $I$ to $\bot(S)$.  This passes through a (possibly
empty) path $a^{-1}$, followed by an edge labelled $t^{-1}$, where the
length of $a$ is less than $M$ (since it can be chosen to be on the
top of a $2$-cell which has an edge in $I$).
Similarly, at the right end of $I$ we have a path labelled
$bt^{-1}$, where the length of $b$ is less than $M$. The path along
$\bot(S)$ joining the two endpoints of these paths is labelled by the
reduced word freely equal in $F$ to $\phi^{-1}(awb) =
\phi^{-1}(a)w\phi^{-1}(b)$.  The only non-constant edges in this word
come from $\phi^{-1}(a)$ and $\phi^{-1}(b)$, which have lengths at
most $M.M_{inv}$.  This proves the assertion in the first sentence.
 
The assertion in the second sentence follows from the observation that
if $x,\, y$ and $\phi(x\beta y)$ consist only of constant letters,
then so does the reduced form of $\beta$, and the assertion in the
final sentence follows immediately from the first.
\end{proof} 
 
\begin{remark}\label{decreell} 
 It is convenient to assume that $MC_4 < \ll$. (In the unlikely 
event that this is not the case, we simply increase $\ll$.) 
\end{remark}

We are finally in a position to make an appropriate definition of a team.

\begin{definition}\label{newTeams} \label{shortDef}
All pre-teams $\pTmm$ with $\hat 
t_1(\mu_1,\mu')\ge\time(S_0)$ 
are defined to be teams, but the qualification criteria for pre-teams with 
 $\hat t_1(\mu_1,\mu')<\time(S_0)$ are 
more selective. 
 
If the genesis of  $\pTmm$ is of type  (G1) or (G2), then 
 the rightmost component of the pre-team may  form a pre-team 
at times before $\ptmm$.  In particular,  it may happen  
that $(\mu_1,\mu')\in\pTmm$ but $\ptmm > \hat t_1(\mu_1,\mu')$ and hence 
$(\mu,\mu')\not\in\pT(\mu_1,\mu')$. To avoid double 
counting in our estimates on $\|\T\|$ we disqualify the  
(intuitively smaller) pre-team $\pT(\mu_1,\mu')$ in these settings.  
 
If the genesis of $\pTmm$ is of type (G4), then again it may happen 
that what remains to the right of $\pTmm$ at some time  before $\ptmm$ is a pre-team. 
In this case, we disqualify the (intuitively larger) pre-team  $\pTmm$.   
  
The pre-teams that remain after these disqualifications 
are now defined to be {\em teams}.  
 
A typical team will be denoted $\T$ 
and all hats will be dropped from the notation for their associated objects 
(e.g. we write $\Emm$ instead of $\pEmm$).  
 
A team is said to be {\em short} if $\|\T\|\le \ll$
or $\sum\limits_{(\mu,\mu')\in\T} |C_{(\mu,\mu')}(2)| \le \ll$. Let
$\S$ denote the set of short teams.
\end{definition} 
 
\begin{lemma} \label{G4lemma} Teams of genesis (G4) are short. 
\end{lemma} 
 
\begin{proof} Lemma \ref{G34pics} implies that $\ET$ is in the immediate 
future of an interval of length at most $C_4$. And we have decreed (Remark 
\ref{decreell}) that $MC_4< \ll$. 
\end{proof} 
 
We wish our ultimate definition of a team to be such that every  pair  $(\mu,\mu')$ 
with $C_{(\mu,\mu')}(2)$ non-empty is assigned to a team. The above definition 
fails to achieve this because of two phenomena: first, a pre-team 
 $\pTmm$ with genesis of type (G4) may
have been disqualified, leaving $(\mu,\mu')$ teamless; second, in our initial discussion of  
pre-teams (the first paragraph of Section \ref{t1}) we excluded pairs $(\mu,\mu')$ 
with $|C_{(\mu,\mu')}(2)|\le 2B$. The following definitions remove these difficulties. 
 
\bd[Virtual team members] \label{Virtual} 
If a pre-team $\pTmm$ of type (G4) is disqualified under the terms of Definition \ref{newTeams} 
and the smaller team necessitating disqualification is $\pT(\mu_1,\mu')$,  
then we define $(\mu,\mu')\vin\pT(\mu_1,\mu')$ and $\pTmm\subset_v\pT(\mu_1,\mu')$. 
We extend the relation $\subset_v$ to be transitive and extend $\vin$ correspondingly. 
If $(\mu,\mu')\vin\T$ then $(\mu_2,\mu')$ is said to be a {\em virtual member} of  
the team $\T$. 
\ed 
 
\bd If  $(\mu,\mu')$ is such that  $1\le |C_{(\mu,\mu')}(2)|\leq 2B$ and 
$(\mu,\mu')$ is neither a member nor a virtual member of any previously 
defined team, then we define $\T_{(\mu,\mu')}:=\{(\mu,\mu')\}$ to be a
(short) team with $\|\T_{(\mu,\mu')}\|=|C_{(\mu,\mu')}(2)|$.  
\ed 
 
\begin{lemma}\label{allIn}  
Every   $(\mu,\mu')\in\vecZ$ with $C_{(\mu,\mu')}(2)$ non-empty is a member 
or a virtual member of exactly one team, and there are less than $2\n$ teams. 
\end{lemma} 
 
\begin{proof} The first assertion is an immediate consequence of the preceding 
three definitions, and the second  follows 
 from the fact that $|\vecZ| < 2\n$. 
\end{proof}

\subsection{Pincers associated to teams of Genesis (G3)}
 
In this subsection we describe the pincer $\Pin_\T$ canonically
associated to each team of genesis $(G3)$. 
The definition of $\Pin_\T$
involves the following concept which will prove important also for teams
of other genesis.

\begin{definition}\label{narrowPast}
We define the {\em narrow past} of a team 
$\T$ to be the set of constant edges that have a future in $\ET$. The narrow 
past may have several components at each time, the set of which  
are ordered left to right according to the ordering in $\ET$ of their futures.  We call these components {\em sections}. 
\end{definition}

{\center{\em{For the remainder of this subsection we 
consider only long teams of genesis (G3). }}}

\begin{definition}[The Pincer $\tilde\Pin_\T$] \label{pl}\label{t2}
The paths labelled $\hat p_l$ and $\hat p_r$ in Figure \ref{G3Pic}
 determine a pincer and are defined as follows. Let $\xT$ be the
 leftmost non-constant edge to the right of $\mu$ in the immediate
 past  of  $\T$, and let  $x_1(\T)$  be the edge that consumes it.   
Define $\tplT$ to be the path in $\F$ that traces the history of $\xT$
to the boundary,  and let $\tprT$ be the path that traces the history
of  $x_1(\T)$.
(Note that  $x_1(\T)$ is left-fast.)

Define $\tilde t_2(\T)$ to be the earliest time at which the  
paths $\tplT$ and $\tprT$ lie in the same corridor.
The segments of the paths $\tplT$ and $\tprT$ after this time, together
with the path joining them along the bottom of the
corridor at  time $\tilde t_2(\T)$ form a pincer. We denote this pincer
$\tilde\Pin_{\T}$. 
\end{definition}

The Pincer Lemma argues for the regular disappearance of colours
within a pincer during those times when more than two colours continue
to survive along the corridors of $\tilde\Pin_\T$. 
However, when there are only two colours the situation  is
more complicated.  

We claim that the following situation cannot arise:
$\time(S_{\hat\Pin_\T}) \leq \tone - T_0$, the path
$\tplT$ and the entire narrow past of $\T$ are in the same corridor at
time $\tone - T_0$, and at this time they 
are separated only by constant edges. For if this were the case,
then the colour of $\tprT$ would $\phi$-neuter the colour of $\tplT$
but would take more than $T_0$ steps to do so, contradicting
the Two Colour Lemma.  Thus at least one of the three hypotheses in
the first sentence
of this paragraph is false; we consider the three possibilities. The
troublesome case (3) leads to a cascade of pincers as depicted in
Figure \ref{cascade}.

\begin{definition}[The Pincer $\Pin_{\T}$ and times $t_2(\T)$ and
$t_3(\T)$] \label{PincerDef} 

\ 

\begin{enumerate}
\item {\em Some section of the narrow past of $\T$ is not in the same corridor as
 $\tplT$ at time $\tone - T_0$:} In this case\footnote{this includes the
possibility that $\tplT$ does not exist at time  $\tone - T_0$}
we define $t_2(\T)=t_3(\T)$ to be the earliest time at which the entire
narrow past of $\T$ lies in the same corridor as $\tplT$ and has length at least
$\ll$.
\item {\em Not case (1), there are no non-constant edges between $\tplT$
and the narrow past of $\T$ at time $\tone - T_0$:} In this case
$\time(S_{\tilde\Pin_{\T}}) > \tone - T_0$.  We define
$\Pin_{\T} = \tilde\Pin_{\T}$ and  $t_3(\T) =
\time(S_{\Pin_{\T}})$. If the narrow past of $\T$ at time $\tone - T_0$ 
has length less than $\ll$, we define
$t_2(\T) =  t_3(\T)$, and otherwise $t_2(\T) =\tilde t_2(\T)$. 
\item 
{\em Not in case (1) or case (2):}
In this case there is at least one non-constant edge between the
narrow past of $\T$ and $\tplT$ at 
$\tone - T_0$. We  pass to the latest time at which there is such an
intervening
non-constant edge and consider  the path $\tilde p_l^\prime(\T)$
that traces the history of the
leftmost intervening non-constant edge $x'(\T)$ and the path $\tilde
p_r^\prime(\T)$
that traces the history of the edge $x_1^\prime(\T)$ that cancels
with $x'(\T)$.
We define  $\tilde t_2'(\T)$ to be the earliest time at which the  
paths $\tilde p_l^\prime(\T)$ and $\tilde p_r^\prime(\T)$  lie in the
same corridor
and consider  the pincer formed by the  segments of the paths $\tilde
p_l^\prime(\T)$ and $\tilde p_r^\prime(\T)$
after  time $\tilde t_2'(\T)$ together
with the path joining them along the bottom of the
corridor at  time $\tilde t_2^\prime(\T)$.

We now repeat our previous analysis with the primed objects $\tilde p_l^\prime(\T), \tilde t_2^\prime(\T)$ {\em etc.} in
place of $ \tilde p_l(\T), \tilde t_2(\T)$ {\em etc.}, checking whether we now fall into case (1) or (2);
if we do not then we pass to $ \tilde p_l''(\T), \tilde t_2''(\T)$ {\em etc.}, and iterate the analysis until
we do indeed fall into case (1) or (2), at which point we acquire the desired definitions of 
$\Pin_\T,\, t_2(\T),\,  t_3(\T)$.
\end{enumerate}

Define $p_l(\T)$ (resp. $\prT$) to be the left (resp. right)
boundary path of the pincer $\Pin_\T$ extended backwards
in time through $\F$ to $\partial\Delta$. Define $p_l^+(\T)$ to be the
sequence of non-constant edges (one at each time) lying immediately to
the right of the narrow past of $\T$ from the top of $\Pin_{\T}$ to
time $\tone$. (These are edges of the leftmost of the primed $\tplT$
considered in case (3).)
\ed

\begin{definition}
Let $\T$ be a long team of genesis (G3).  Let $\chi_P(\T)$ be the set
of colours containing the paths $\tilde p_l(\T), \tilde p_l'(\T),\tilde p_l''(\T),\dots$ that
arise in (iterated applications of) case (3) of Definition \ref{PincerDef} but
do not become $p_l(\T)$.
\ed 

\begin{figure}[htbp] 
\begin{center} 
  
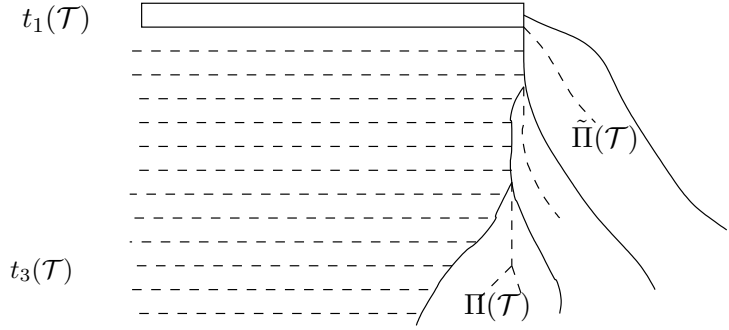 
  
\caption{The cascade of pincers.}
\label{cascade} 
\end{center} 
\end{figure}

The preceding definitions are framed so as to make the following important
facts self-evident.

\begin{lemma} 
 \label{t1t3forTeam}\label{disj1}

\ 

\begin{enumerate}
\item
If $\T$ is a long team of genesis (G3),  
\[      t_1(\T) - t_3(\T) \leq T_0(|\chi_P(\T)| + 1).   \] 
\item
If $\T_1$ and $\T_2$ are disjoint then  $\chi_P(\T_1)\cap \chi_P(\T_2)=\emptyset$.
\end{enumerate}
 \end{lemma}

\subsection{The length of teams} \label{TeamLemmas}

\begin{definition}\label{down1}  
Define $\down_1({\mathcal T})\subset\partial\Delta$ to consist of those 
edges $e$ that are labelled $t$ and satisfy one of the following conditions: 
\begin{enumerate} 
\item[1.] $e$ is at the left end of a corridor containing a section of the narrow 
past of $\T$ that is not leftmost at that time; 
\item[2.]  $e$ is at the right end of a corridor containing a section of the narrow 
past of $\T$ that is not rightmost at that time;  
\item[3.]  $e$ is at the right end of a corridor which 
contains the rightmost section of the narrow past of $\T$ at that time but which does 
not intersect $\plT$.\\ 
\end{enumerate} 
\end{definition} 
 
All of the edges shown on the boundary in
Figure \ref{figure:TeamAge} are contained in $\down_1(\T)$. 

\begin{definition} Define $\partial^\T\subset\partial\Delta$ to be the
set of (necessarily constant) edges that have a preferred future in
$\ET$.
\end{definition}

We record an obvious disjointness property of the sets defined above.

\begin{lemma}\label{disj2}

\ 

\begin{enumerate}
\item For distinct teams $\T_1$ and $\T_2$, $\partial^{\T_1}$ and
$\partial_{\T_2}$ are disjoint.
\item For distinct teams $\T_1$ and $\T_2$, $\down_1(\T_1)$ and
$\down_1(\T_2)$ are disjoint. 
\end{enumerate}
\end{lemma}

\bd\label{QT}
Suppose that $\T$ is a team of genesis (G3).  We define  $\QT$ be the
set of edges $\e$ with the following properties:
$\plT$ passes through $\e$  before time $t_3(\T)$, and the corridor $S$
with $\e\in\bot(S)$ contains the entire narrow past of $\T$ and
this narrow past has length at least $\ll$.
\ed

The following lemma gives us a bound on $|\ET|$, which will reduce 
our task to that of bounding $|\QT|$ for teams of genesis (G3). 
 
\begin{lemma}  \label{TeamAgeLemma}  

\ 
 
\begin{enumerate} 
\item[1.] If the genesis 
of $\T$ is of type (G1) or (G2), then 
$$
\|\T\| \leq 2MC_4\,|\down_1({\mathcal T})| +  |\partial^\T|  .
$$ 
\item[2.] If the genesis of $\T$ is of type (G3), then 
$$
\|\T\| \leq 2MC_4\,|\down_1({\mathcal T})| +|\partial^\T| + 2MC_4\,|\QT| +
2MC_4T_0\big(|\chi_P(\T)| +1\big) + \ll .
$$  
\end{enumerate} 
\end{lemma} 
 
\medskip 
 
\begin{figure}[htbp] 
\begin{center} 
  
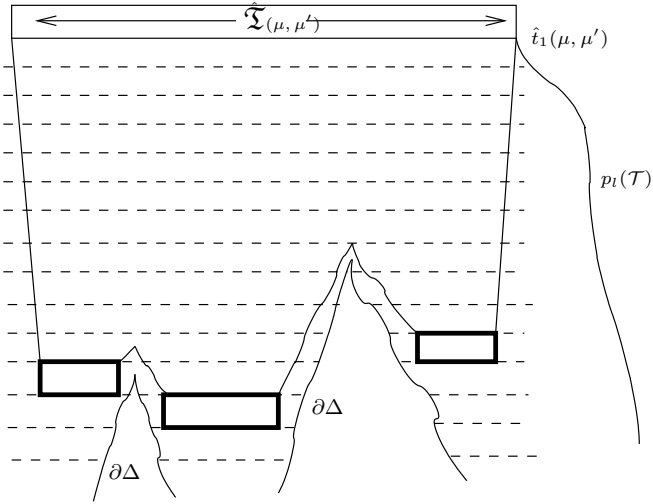 
 
\caption{Bounding the size of a team in terms of $|\down_1|$ and $|p_l|$} 
\label{figure:TeamAge} 
\end{center} 
\end{figure} 
 
\begin{proof}     The first thing to observe is that at any 
stage in the past of $\ET$ the set of letters lying in a 
single corridor form a connected region. As in Lemma \ref{G34pics}, 
this is simply a matter of noting that if $\phi(aub)=w$ where $w, a$ and $b$ 
consist only of constant letters, then $u$ must equal a word in constant letters. 
 
Consider the   past of $\ET$ at a time $t$. Write $k_t$ for the number of 
corridors that contain a non-trivial component of this past. 
The total 
increase in length of these components  when one goes forward to time $t+1$ 
is bounded by $2MC_4k_t$, since the connectedness of the past 
implies that  the only growth that can happen for existing components occurs 
at their extremities, where a block of at most $MC_4$ constant letters
may be added. This follows from Lemma \ref{G34pics}.  Also at time $t+1$,
constant letters from $\partial\Delta$ may join the past of $\ET$, and
there may be new components of constant letters (each of length less
than $2MC_4$) whose ancestors at time $t$ were non-constant
letters. Thus we have three possible causes of increase. The first and
third account for growth of at most $2MC_4k_{t+1}$ and the second
(boundary) contribution is the number of elements of $\partial^\T$
that occur at time $t+1$.  
If the genesis of $\T$ is of type (G1) or (G2), then at least $k_{t+1}$ 
edges of $\down_1(\T)$ occur at time $t$, compensating us for the growth summand 
$2MC_4k_{t+1}$. If the genesis of $\T$ is of type (G3) then we still have the 
above compensation {\em except} at those times where no edges of $\down_1(\T)$ occur. 
At these latter times the whole of the narrow past of $\T$ 
lies in a single corridor through which $\plT$ passes. Since the
narrow past lies
in a single corridor, it is connected and  grows at most $2MC_4$ when moving 
forward one unit 
of time (unless added to by $\partial^\T$).  

The summands $2MC_4\,|\QT|$ and  $2MC_4T_0\big(|\chi_P(\T)| +1\big) $ in item (2) of the lemma
account for the growth of the narrow past in the intervals of time
below  $t_3(\T)$, and from $t_3(\T)$ to $t_1(\T)$, respectively. The additional summand $\ll$
allows us to desist from our estimating if the narrow past of $\T$ ever shrinks to have
length less than $\ll$. 
\end{proof}

\subsection{Bounding the size of $\QT$}\label{Proofs} 

For the remainder of this section we concentrate exclusively on long teams of genesis (G3)
with $\QT$ non-empty. We denote the set of such teams by $\Gthree$.  
Our goal is to bound $|\QT|$. (In the light of our
previous results, this will complete the required analysis of the
length of teams.)  

Recall from  Definition \ref{PincerDef} that for teams of genesis 
(G3), the paths $\plT$ and $\prT$ and the chain 
of 2-cells joining them in the corridor at time $\ttwo$ form a pincer denoted $\Pin_{\T}$. 
The set $\subT$ was defined in Definition \ref{chiP}.

An important
feature of  teams in $\Gthree$ is:

\begin{lemma}\label{GotBlock}
If $\T \in \Gthree$ then there exists a block of  at least $\ll$
constant edges   immediately adjacent to $\Pin_{\T}$ at each time from
$t_3(\T)$ to the top of $\Pin_{\T}$, and adjacent to $p_l^+(\T)$ from
then until $\tone$. (At time $\tone$ this block contains $\ET$.)
\end{lemma}

\begin{proof} The hypothesis
that $\QT$ is non-empty means that the narrow past
of $\T$ at  some time before $t_3(\T)$ has length at
least $\ll$ and is contained in the same corridor as 
$\plT$ (see Definition \ref{QT}).  The definition of $t_3(\T)$ implies that
the  narrow past of $\T$ is contained in a block of constant letters 
immediately adjacent to $\plT$ or $p_l^+(\T)$ from time $t_3(\T)$
until $\tone$. Since the  length of the narrow past of $\T$ does not
decrease before $\tone$, these blocks of constant letters must have
length at least $\ll$.
\end{proof}

The following is an immediate consequence of the Pincer Lemma. 

\begin{lemma} \label{t1-t2Lemma}  
For all $\T \in \Gthree$, 
$$ 
  t_3(\T) - \ttwo = \life(\Pin_{\T}) \leq \ttt (|\subT|+1) .   
$$
\end{lemma}

\begin{lemma} \label{nesters}
If $\T_1, \T_2 \in \Gthree$ are distinct teams then $\chi(\Pin_{T_1})
\cap \chi(\Pin_{T_2})  = \emptyset$.
\end{lemma}  
 
\begin{proof} The pincers $\Pin_{\T_i}$ are either disjoint or else
one is contained in the
other. In the latter case, say  $\Pin_{\T_1}\subset\Pin_{\T_2}$,  
the existence of the  block of $\ll$ constant edges established in
Lemma \ref{GotBlock} means that $\Pin_{\T_1}$ is  
actually nested in $\T_2$ in the sense of Definition \ref{chiP}. Thus 
$\chi(\Pin_{\T_1}) \cap \chi(\Pin_{\T_2})  = \emptyset$ (by 
Definition \ref{chiP}). 
\end{proof}
 
\begin{corollary} \label{t1-t2Corr}
$\sum\limits_{\T \in \Gthree} t_3(\T) - t_2(\T) \leq 3\ttt \n$.
\end{corollary}

It remains to bound the number of edges in $\QT$ which occur before $\ttwo$;
this is  cardinality of the following set.
 
\begin{definition}\label{down2} 
For $\T\in\Gthree$ we define $\down_2(\T)$ to be the set of edges in $\partial\Delta$ that 
lie at the righthand end of a corridor containing an edge in $\QT$
before time $\ttwo$.
\end{definition}

The remainder of this section is dedicated to obtaining a bound on
$$
\sum\limits_{\T\in\Gthree}|\down_2(\T)|,
$$ 
(see Corollary \ref{downbound}).
 
At this stage our task of bounding $\|\T\|$ would be complete if
the the sets $\down_2(\T)$ associated to distinct teams 
were disjoint --- unfortunately they need not be, because of the possible 
nesting of teams as shown in Figures \ref{figure:Nest}
and \ref{figure:doublecount}. Thus we shall 
be obliged to seek further pay-off for our troubles. To this end we 
shall identify two sets of consumed colours $\chi_c(\T)$ and 
$\chi_{\delta}(\T)$ that arise from 
the nesting of teams. 
 
In order to analyse  the 
effect of nesting we need the following vocabulary.

There is an obvious left-to-right ordering of those paths in the forest 
$\F$ which begin on the arc of $\partial\Delta\ssm\partial S_0$ that commences 
at the initial vertex of the left end of $S_0$. (First one orders the trees, then 
the relative order between paths in a tree is determined by the manner in 
which they diverge; the only paths which are not ordered relative to each 
other are those where one is an initial segment of the other, and this 
ambiguity will not concern us.) 

\smallskip

\noindent{\bf Notation:} We write $\Gthree'$ for the set of teams $\T
\in \Gthree$ such that $\down_2(\T) \neq \emptyset$.

\smallskip

We shall need the following obvious separation property.

\begin{lemma}\label{separate} 
Consider $\T \in \Gthree'$.  If a
path $p$ in $\F$ is to the left of $\plT$ and a path $q$ is the right
of $\prT$,
then there is no corridor connecting $p$ to $q$ at any time $t<\ttwo$. 
\end{lemma} 

\begin{proof}  The hypothesis $\down_2 (\T) \neq \emptyset$ implies
that before $\ttwo$ the paths $\plT$ and $\prT$ are not in the same
corridor.
\end{proof}

\begin{definition} \label{depthDef} 
$\T_1\in\Gthree'$ is said to be {\em below} $\T_2\in\Gthree'$ if
$p_l(\T_2)$ and $p_r(\T_2)$ both lie  between $p_l(\T_1)$ and
$p_r(\T_1)$ in the left-right ordering described above.

$\T_1$ is said to be {\em to the left} of $\T_2$ if both  $p_l(\T_2)$ and $p_r(\T_2)$ 
lie to the right of $p_r(\T_1)$.

We say that $\T$ is at {\em depth} $0$ if there are no teams above it. 
Then, inductively, we say that a team is at depth $d+1$ if $d$ is the maximum 
depth of those teams above $\T$.   
 
A {\em final depth} team is one with no teams below it. 
 
Note that there is a complete left-to-right ordering of teams $\T\in\Gthree$ at any given depth. 
\end{definition}

\begin{lemma}\label{trapS_0} If there is a team from $\Gthree'$ below
$\T\in\Gthree'$, then $\tone \ge \time (S_0)\ge \ttwo$. 
\end{lemma}

\begin{proof} The first thing to note is that  if
 $\time (S_0)$ were less than $\ttwo$, then  the narrow
past of $\T$  at time $t_2(\T)$ must contain at least
$\ll$ edges. This is because the length of the narrow past of $\T$ cannot
decrease before $\tone$, and
at $\time(S_0)$ the narrow past
is the union of the intervals $C_{(\mu,\mu')}(2)$ with $(\mu,\mu')\in\T$,
which has length at least $\ll$ since $\T$ is assumed not to be short.

Thus if $\time(S_0)<\ttwo$ then  we are in the non-degenerate
situation of Definition \ref{PincerDef} and the defining property of $\ttwo$
means that  before time $\ttwo$  no edge to the right of $\prT$ lies
in the same
corridor as all the colours of $\T$ (cf. Lemma \ref{separate}).  In
particular this is true of
the past of the reaper of $\T$ (assuming that it has a past at time
$\ttwo$). On the
other hand,   the reaper of $\T$ has a past in $S_0$ (by the very
definition of a team), as do all of the colours of $\T$. And since they
lie in a common corridor at $\time(S_0)$, they must also do so
at all times up to $\tone$. This contradiction implies
that in fact $\time(S_0)\ge\ttwo$. 
 
Consider Figure \ref{figure:Nest}.  
Suppose that $\T'\in\Gthree'$ is below $\T$. The proof of Lemma \ref{GotBlock}
tells us that there is a block of constant edges extending from the
top of $\Pin_{\T'}$ containing the narrow past of $\T'$, and there is
a similarly long block
extending from the path $p_l^+(\T)$ at each subsequent time until
$t_1(\T')$. Thereafter the future of the block is contained in the
block of constant edges that evolves into the union of the
$C_{(\mu,\mu')}(2) \subseteq \bot(S_0)$ with $(\mu,\mu') \in \T'$,
which is long by hypothesis.
 
 At no time can this evolving block extend across $\plT$ 
because by  definition the edges along $\plT$ are labelled by non-constant
letters.  Thus the evolving
block is trapped to the right of $\plT$ and to the left of 
$\prT$. In particular, it must vanish entirely before the time
at the top of the pincer $\Pin_\T$, which is no later than
$\tone$ and therefore $\tone \ge \time (S_0)$. 
\end{proof} 
 
The following is the main result of this section.

\begin{lemma}   \label{Aget2Lemma} \label{chiT} There exist sets of
colours $\chi_c(\T)$ and $\chi_{\delta}(\T)$ associated to each team
$\T\in\Gthree'$ such that the sets associated to
distinct teams are disjoint and the following inequalities hold. 

For each fixed team $\T_0 \in \Gthree'$ (of depth $d$ say),  
the teams  of depth $d+1$ that lie below   $\T_0$ 
may be described as follows: 
\begin{enumerate} 
\item[$\bullet$] There is at most one {\em distinguished team} $\T_1$,
and
$$ 
\|\T_1\|\le 2B\Big(\ttt (1+ |\chi(\Pin_{\T_0})|) + T_0(|\chi_P(\T_0)| +
1)\Big).
$$
\item[$\bullet$] There are some number of final-depth teams. 
\item[$\bullet$] For each of the remaining teams $\T$ we have 
  $$|\down_2(\T_0) \cap\down_2(\T)| \le \ttt  \Big( 1 
+ |\chi_c(\T)| \Big) + T_0 \Big( |\chi_{\delta}(\T)| + 2 \Big).
$$ 
\end{enumerate}
\end{lemma}  

\begin{proof}  The first thing to note is that if two teams $\T, \T'
\in \Gthree'$ are at the same depth, then $\down_2(\T)$ and
$\down_2(\T')$ are disjoint. Indeed if $\T$ is to the left of $\T'$,
then at times before $\ttwo$ the paths $\plT$ and $p_l(\T')$ never lie in
the same corridor.  Let $\T \in \Gthree'$ be a team of level $d+1$
that is below $\T_0$ and consider the edge $e$  at the right end of a
corridor earlier than $\ttwo$ that contains an edge in $\QT$. We are
concerned with the fact that this edge may be in $\down_2(\T_0)$.  In
this situation we say that $\T_0$ and $\T$ {\em double count} $e$.

\begin{figure}[htbp] 
\begin{center} 
  
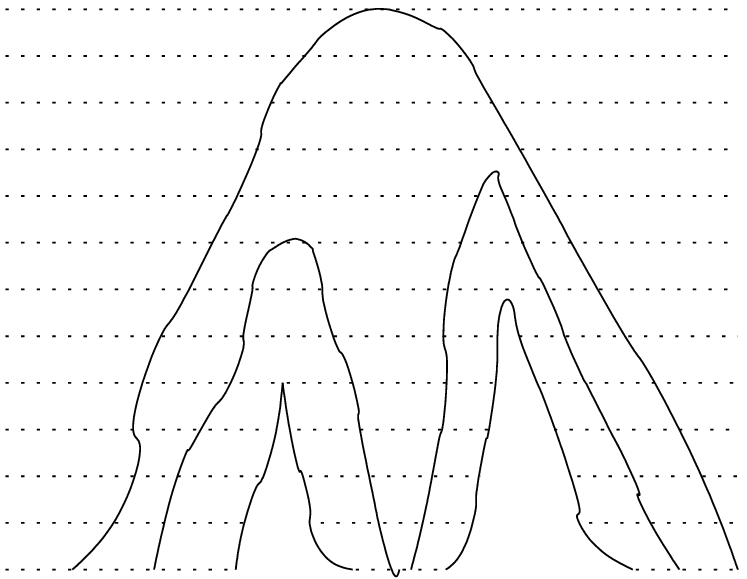 
  
\caption{A depiction of double-counting} 
\label{figure:doublecount} 
\end{center} 
\end{figure} 
 
Let $\T_1, \dots ,\T_r$ be the teams in $\Gthree'$ of depth $d+1$
which double-count with $\T_0$, ordered from left to right, with the
final-depth teams deleted.  We define $\chi_c(\T)$ to be empty for
teams not on this list. $\T_1$ will be the distinguished team.

Since there is no double-counting between teams of the same level, 
 the sets of times  at which 
$\T_1, \dots,\T_r$ double-count with $\T_0$ must be disjoint. Indeed if 
 $i < j$ then the set of times at which $\T_i$ double-counts 
  with $\T_0$ is earlier than the 
set of times at which $\T_j$ 
double-counts with ${\T_0}$ (Lemma \ref{separate}). Moreover, 
 the times for each $\T_i$ form an interval, which we denote $\I_i$.  
 
We assume $r\ge 2$ and describe the construction of 
the sets $\chi_c(\T_i)$ and $\chi_{\delta}(\T_i)$ that account for
double-counting.

The first thing to note is that each $\I_i$ must be later than $t_2(\T_1)$, 
 by Lemma \ref{separate}. 
The second thing to note is that   the entire interval of time $\I_i$ 
must also be earlier than $t_1(\T_1)$. Indeed if  some double-counting by 
$\T_i$  and $\T_0$  were to occur  after $t_1(\T_1)$, then we would
have $t_2(\T_k) > t_1(\T_1)$. But then
$\time (S_0) > t_1(\T_1)$, so Lemma \ref{trapS_0} would imply  that
there was no team below $\T_1$, contrary to hypothesis.  

We separately consider the intervals $\I_i \cap [ t_2(\T_1),t_3(\T_1)
]$ and  $\I_i \cap [ t_3(\T_1),t_1(\T_1) ]$, whose union is all of
$\I_i$.

For that part of $\I_i$ before $t_3(\T_1)$, the proofs of
the Pincer Lemma (Theorem \ref{PincerLemma}) and Proposition 
\ref{prePincerLemma} tell us that colours in $\chi(\Pin_{\T_1})$
will be consumed at the rate of at
least one per $\ttt$ units of time. Define
$\chi_c(\T_i)$ to be this set of consumed colours. We have
$$
\Big|\, \I_i \cap [ t_2(\T_1),t_3(\T_1) ]\, \Big| \leq \ttt (1 +
|\chi_c(\T_i)|) .
$$

Now consider  $\I_i \cap [ t_3(\T_1),t_1(\T_1)]$. Define
$\chi_{\delta}(\T_i)$ as follows. The discussion in Definition
\ref{PincerDef} shows that in any period of time of length $T_0$ in
the interval $[t_3(\T_1),t_1(\T_1) ]$ at least one colour in
$\chi_P(\T_1)$ disappears.  Let $\chi_{\delta}(\T_i)$ be
the set of colours in $\chi_P(\T_1)$ which disappear during $\I_i
\cap [t_3(\T_1), t_1(\T_1)]$ (these disappearances correspond to the
discontinuities in the `path' $p_l^+(\T_1)$).  By construction, we
then have\footnote{There is a 2 rather than the familiar 1 on the
right to account for the colour containing
$p_l(\T_1)$, which is not included in
$\chi_P(\T_1)$; there might be up to $T_0$ corridors between
$t_3(\T_1)$ and the top of $\Pin_{\T_1}$.}

\[	\Big| \, \I_i \cap  [ t_3(\T_1),t_1(\T_1)]\, \Big| \le
T_0(|\chi_{\delta}(\T_i)| + 2),	\]
and combining these estimates we have
$$
|\I_i| \le \ttt  \Big( 1 + |\chi_c(\T_i)| \Big) + T_0
\Big(  |\chi_{\delta}(\T_i)| +  2\Big) , 
$$
as required. 
Since the intervals $\I_i$ are disjoint, 
the sets $\chi_c(\T_i),\, i=2,\dots,r$ are mutually disjoint. 
And by construction, these sets are also disjoint from the sets 
associated to teams other than 
the $\T_i$ under consideration (i.e. those under other depth $d$  
teams, or those of different 
depths). The same considerations hold for the sets
$\chi_{\delta}(\T_i),\, i=2,\ldots, r$.

In Figure \ref{figure:double}, the shaded region is where we recorded
the regular disappearance of the colours forming $\chi_c(\T_i)$,
whilst in Figure \ref{figure:doubletwo}, the shaded region is where we
recorded the regular disappearance of the colours forming
$\chi_{\delta}(\T_i)$.

\begin{figure}[htbp] 
\begin{center} 
  
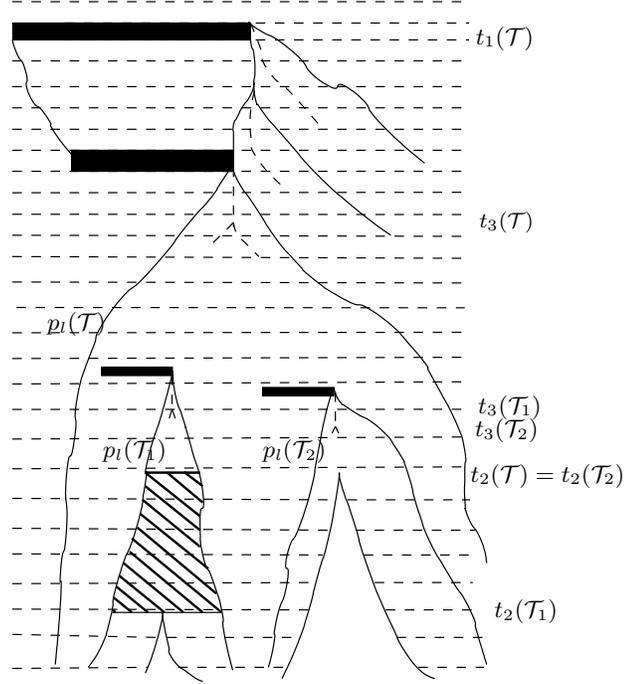 
  
\caption{Finding the colours $\chi_c(\T_i)$} 
\label{figure:double} 
\end{center} 
\end{figure}

\begin{figure}[htbp]
\begin{center}

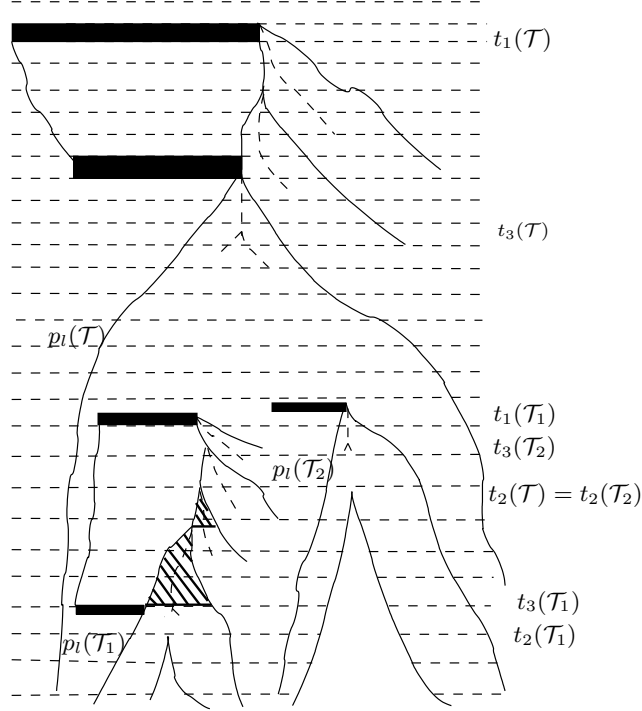

\caption{Finding the colours $\chi_{\delta}(\T_i)$}
\label{figure:doubletwo}
\end{center}
\end{figure}
 
It remains to establish the inequality 
$$ 
\|\T_1\|\le 2B\Big( \ttt (|\chi(\Pin_{\T_0})| + 1) +
(|\chi_P(\T_0)|+1) \Big). 
$$ 
We first note (as in the proof of Lemma \ref{trapS_0})  that 
$\ET_1$ is trapped between $\plT$ and $\prT$, so it must be consumed
entirely between  the times $t_1(\T_1)$ and  $t_1(\T_0)$. But by the
Bounded Cancellation Lemma, the length of  the future of $\ET_1$ can
decrease by at most $2B$ at each step in time. Therefore  $\|\T_1\|\le
2B (t_1(\T_0) - t_1(\T_1))$.
 
$\T_1$ is assumed not be final-depth, so from Lemma \ref{trapS_0} we have 
$t_2(\T_0) \le \time (S_0) \le t_1(\T_1)$. By combining these  
 inequalities with Lemmas \ref{t1-t2Lemma} and \ref{t1t3forTeam}
we obtain: 
\begin{eqnarray*} 
\|\T_1\| &\le& 2B\, \Big(t_1(\T_0) - t_1(\T_1)\Big)\\  
& \leq & 2B\, \Big(t_1(\T_0) - \time (S_0)\Big)\\  
& \leq & 2B\, \Big(t_1(\T_0) - t_2(\T_0)\Big)\\  
& \leq & 2B \Big[ \ttt \Big(1+|\chi(\Pin_{\T_0})|\Big) + T_0
\Big( |\chi_P(\T_0)| + 1 \Big) \Big]. 
\end{eqnarray*} 
\end{proof} 
 
\begin{corollary} \label{Firstdown2sum} Summing over the set of teams $\T \in \Gthree'$ that
are not distinguished, we get 
$$ 
\sum_\T \Big|\down_2 (\T)\Big| \le2\,\Big| \bigcup_\T \down_2(\T)\Big|
+\sum_{\T}  \ttt  \Big( 1 
+ |\chi_c(\T)| \Big) + \sum_{\T} T_0 \Big( |\chi_{\delta}(\T)| + 2
\Big) .
$$ 
\end{corollary} 
 
\begin{proof} Suppose $\T \in \Gthree'$ of depth $d+1$ is not
final-depth and not distinguished, and that $\T$ double-counts with
some $\T_0$ of depth $d$ above it.  Then, by Lemma \ref{Aget2Lemma},
we have
\begin{eqnarray*}
|\down_2(\T)| & = & |\down_2(\T) \smallsetminus \down_2(\T_0)| +
|\down_2(\T) \cap \down_2(\T_0)| \\
& \leq & |\down_2(\T) \smallsetminus \down_2(\T_0)| + \ttt (1 +
|\chi_c(\T)|) + T_0 (2 + |\chi_{\delta}(\T)|).
\end{eqnarray*}
Suppose that $\T' \in \Gthree'$ is a team of depth $k < d$ and that
$\T'$ is above $\T$. If $\T$ double-counts with $\T'$ at time $t$,
then $\T$ double-counts with $\T_0$ at time $t$, by Lemma
\ref{separate}.  Therefore, the set of edges that $\T$ double-counts
with any team of lesser depth is exactly $\down_2(\T) \cap
\down_2(\T_0)$.

Thus we have accounted for all double-counting other than than
involving final depth teams.  The factor $2$ in the statement of the
corollary accounts for this.
\end{proof} 

And summing over the same set of teams again, we obtain: 
\begin{corollary} \label{downbound} 
$$ 
\sum_\T |\down_2(\T)| \ \le \ \n (2 + 3\ttt + 5T_0). 
$$
\end{corollary} 
\begin{proof} The sets of colours $\chi_c(\T)$ and $\chi_{\delta}(\T)$
are disjoint. And the
union of the sets $\down_2(\T)$  is a subset of $\partial\Delta$. The
set of all colours and the  set of edges in $\partial\Delta$ each have
cardinality at most $\n$.   And the number of teams is less than $2\n$
(Lemma \ref{allIn}).
\end{proof}

\section{The Bonus Scheme} \label{BonusScheme}

We have defined teams and obtained a global bound on 
$\sum\|\T\|$. 
If $\cmm$ is non-empty then $(\mu,\mu')$ is a member or  
 virtual member of a  unique team. 
If this team is such that $\tone \ge \time (S_0)$, then no member of 
the team is virtual  and we have the inequality
$$\|\T\|>\sum\limits_{(\mu,\mu') \in \T}|\cmm| - B$$
 established in Lemma \ref{t1high}.
 We indicated following this lemma how this inequality
 might fail in the case where
  $\tone < \time (S_0)$. In this section we take up this
  matter in detail
   and introduce a {\em bonus scheme} that
 assigns additional edges to teams in order to compensate for the possible failure
 of the above inequality when   $\tone < \time (S_0)$.
 
By definition, at time $\tone$ the reaper $\r=\rT$ lies
 immediately to the right of $\ET$.  The edges of 
$\ET$ not consumed from the right by $\r$ by  $\time (S_0)$
 have a preferred future in $S_0$ 
that lies in $\cmm$ for some member $(\mu,\mu')\in\T$. 
However, not all of the edges of 
$\cmm$ need arise in this way:
 some may not have  a constant ancestor at time $\tone$.
And  if $(\mu,\mu')$ is only a virtual member of $\T$,
 then no edge of $\cmm$ lies in the 
future of $\ET$. The  {\em bonus} edges in $\cmm$
are a certain subset of those  that do not have a constant 
ancestor at time $\tone$. They are defined
as follows.

\begin{definition} Let $\T$ be a team with $\tone < \time(S_0)$
and consider a time $t$ with $\tone < t < \time(S_0)$.

The {\em swollen future} of $\T$ at time $t$ is the interval 
of constant edges beginning immediately to the left of the pp-future of $\rT$. 
 
Let $e$ be a non-constant edge that lies immediately to the left of the 
swollen future of $\T$ but whose ancestor is not a 
right para-linear edge in this position. If $e$ is a right para-linear and
 the (constant) rate 
at which $e$ adds letters to the swollen future of $\ET$ is greater 
than the (constant) rate at  which the future of the reaper cancels letters 
in the future of $\ET$, then we define $e$ to be 
a {\em rascal}; if $e$ is right-fast then we define it to be a {\em terror}.
In both cases,
 we define the
{\em bonus provided by $e$} to be the set of edges in the swollen future
of $\T$ in $S_0$ that have $e$ as their most recent non-constant
ancestor, and are eventually consumed by $\rT$. 

The set $\bonusT$ is the union of the bonuses provided to $\T$ by all
rascals and terrors.
\end{definition} 
 
\begin{lemma}  \label{C1toTeamLength} 
For any team $\T$, 
\[      \sum_{(\mu,\mu') \in \T \mbox{ \tiny{or} } (\mu,\mu')
 \vin T}|\cmm| \leq \|\T\| + |\bonus(\T)| + B.    \] 
\end{lemma}

\begin{proof}
If $\tone \geq \time(S_0)$, this follows immediately from Lemma
\ref{t1high}.  If $\tone < \time(S_0)$ then at each step in time
between $\tone$ and $\time(S_0)$ the only possible cause of growth in
the length of the swollen future of the team is the possible action of
a rascal or terror if such is present at that time.  (There is no
interaction of the swollen future with the boundary or singularities,
because of the exclusions in the second paragraph of Definition
\ref{newTeams}.)

The swollen future has length $\|\T\|$ at time $\tone$ and length at
least  $\sum |C_{(\mu,\mu')}(2)|$ at $\time(S_0)$. By definition,
$|\bonus(\T)|$ is a bound on the growth in length between these
times. (The summand $B$ is thus unnecessary in the case  $\tone <
\time(S_0)$.)
\end{proof}
 
 The following lemma shows that our main task in this
 section will be to analyse the behaviour of rascals.

\begin{lemma} The sum of the lengths of the bonuses 
provided to all teams by terrors is less than $2M\n$.
\end{lemma}  

\begin{proof} Since it is right-fast, a terror will be separated
from the team to which it is associated after one unit of
time, and hence the bonus that it provides is less than $M$.
There is at most one terror for each possible adjacency of colours
and hence the total contributions of all terrors is less than
$2M\n$.
\end{proof}
 
The typical pattern of influence of rascals on a team 
is shown in Figure \ref{belowS0};
there may be several times at which rascals
appear at the left of $\T$ and provide a  
bonus for the team before being consumed from the left (or otherwise detached 
from the team).  
 
\begin{figure}[htbp] 
\begin{center} 
  
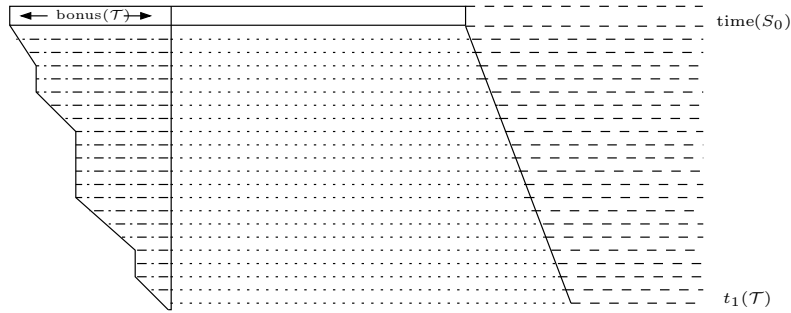 
  
\caption{The generic situation below $\time(S_0)$.} 
\label{belowS0} 
\end{center} 
\end{figure}

\begin{definition}[Rascals' Pincers] \label{RascalPin}
We fix a team $\T$ with $\tone < \time(S_0)$
 and consider the interval of time $[\tau_0(e),\tau_1(e)]$,
where $\tau_0(e)$ 
is the time at which a rascal $e$ appears at the left end of the swollen 
future of $\T$, and $\tau_1(e)$ is the time at which its future is no 
longer to the immediate left of the future of the swollen future of $\T$.  
 
In the case where the pp-future $\hat e$ of $e$ at time  $\tau_1(e)$
is cancelled from the left by an edge $e'$, we define
$\tau_2(e)$ to be the earliest time when the pasts of $\hat e$ and $e'$ are
in the same corridor. The path in $\F$ that traces the pp-future of
$e$ up to  $\tau_1(e)$ is denoted $p_e$ and the path following through the
ancestors of  $e'$ from
 $\tau_2(e)$ to $\tau_1(e)$ is denoted  $p_e'$.
 The pincer\footnote{to lighten the terminology, here we allow the
 degenerate case where the ``pincer" has
  no colours  other than those of $e$
 and $e'$}
  formed by $p_e$ and $p_{e}'$  with base at time  $\tau_2(e)$
is denoted  $\Pin_e$. 
\end{definition}
 
\begin{lemma} \label{lowRascals} 
The total of all bonuses provided to all teams by  rascals $e$ with 
 $\tau_1(e)\le\time(S_0)$ is less than $(3\ttt + 2T_0 +1)M\n$.  
\end{lemma} 
 
\begin{proof} Consider a rascal $e$. We defer the case where $e$ hits
a singularity or the  boundary. If this does not happen, the
pp-future $\hat e$ of $e$ at time $\tau_1(e)$ is cancelled from the left
by an edge $e'$ (which  is right-fast since $e$ is not
constant). We consider the pincer $\Pin_e$ defined above.
The presence of the swollen future of $\T$ at the top of the
pincer allows us to apply the  Two Colour Lemma to conclude that
$\tau_1(e) - T_0 \geq \time(S_{\Pi_e})$ (in the degenerate case
discussed in the footnote, $ \time(S_{\Pi_e})$ is replaced by
$\tau_2(e)$). And the
Pincer Lemma tells us that 
\[      \tau_1(e) - \tau_2(e) \leq \ttt \Big( 1 + |\chi(\Pin_e)| \Big)
+ T_0.       \]
In fact, we could use $\tilde \chi(\Pin_e)$ instead of
$\chi(\Pin_e)$ in this estimate because there cannot be any nesting
amongst the pincers $\Pin_e$ with $\tau_1(e)\le \time(S_0)$,
because nesting would imply that the swollen future of $\T$, which is
immediately to the right of the lower rascal, would be trapped beneath
 the upper pincer,
  contradicting the fact that the team has a non-empty future in $S_0$. 
 
In the case where $e$ hits the boundary or is separated from the team
 by a singularity (at time $\tau_1(e)$) we 
define  $\tau_2(e)=\tau_1(e)$. No matter what the
fate of $e$, we define
 $\partial^e$ to be the set of edges in
 $\partial\Delta$ at the left ends of corridors  
containing the future of $e$ between 
$\tau_0(e)$ and $\tau_2(e)$.
The sets $\partial^e$ assigned to different rascals are disjoint,
so summing over all rascals with $\tau_1(e)\le \time(S_0)$
we have
\begin{eqnarray*}
\sum_e \Big(\tau_1(e) - \tau_0(e)\Big) & = & \sum_e (\tau_1(e) - \tau_2(e)) +
 (\tau_2(e)-\tau_0(e))\\
& \le & \sum_e \ttt\Big(1 + |\chi (\Pin_e)|\Big) + T_0 + |\partial^e|.
\end{eqnarray*}
Since the sets $\chi(\Pin_e)$ and $\partial^e$ are disjoint,
the terms $\ttt |\chi (\Pin_e)|$ and $ |\partial^e|$ 
contribute less than $(\ttt + 1)\n$ to this sum. And since the number of
 rascals is bounded by the number of possible adjacencies of colours, the remaining
terms contribute at most $(\ttt + T_0)2 \n$. Thus
$$
\sum_e  \Big(\tau_1(e) - \tau_0(e)\Big)\ \le\ (3\ttt + 2T_0 + 1) \n .
$$
The bonus produced by each rascal in each unit of time is less than $M$, so  
the lemma  is proved.
\end{proof} 

It remains to consider the size of the bonuses provided by rascals $e$ with 
 $\tau_1(e)>\time (S_0)$.

The bonuses that are not accounted for in Lemma \ref{lowRascals}
reside in blocks  of constant edges along $\bot(S_0)$ each of which is the swollen 
future of some team, with
 a \rpl letter at its left-hand end (the pp-future of
  a rascal) and a \lpl letter at its left-hand end (the pp-future of
  the team's reaper).  
 
 \begin{definition}
A {\em left-biased} rascal $e$ is one with
 $\tau_1(e) > \time(S_0)$ that satisfies the following properties: 
\begin{enumerate} 
\item[1.] the pp-future of the rascal  is (ultimately) 
consumed from the left by an edge of $S_0$,  
\item[2.] the swollen future of  $\T$ at 
time $\tau_1(e)$ has length at least $\ll$ and  
the pp-future of the reaper $\rT$ is still immediately to its right.
\end{enumerate} 
 \end{definition}
 
\def\life{\text{\rm{life}}}  
\def\B{\text{\euf{B}}}

\begin{definition} Let $\B\subset\bot(S_0)$ be an interval of constant edges with a  
right para-linear letter at its left-hand end and a left-linear letter $\r$ at its right-hand end. We 
say that $\B$ is {\em right biased} if $\r$ is ultimately consumed by an edge (to its right) 
in $S_0$. We define $\life(\B)$ to be the difference  between $\time(S_0)$ and  the time at which 
the \lpl letter $\rho$ is consumed. 
And we define the {\em effective volume} of $\B$ to be the number of edges in $\B$ 
that are ultimately consumed by $\rho$. 
\end{definition} 
 
We have the following tautologous tetrad of possibilities covering the swollen teams whose bonuses are 
not entirely accounted for by Lemma \ref{lowRascals}. 
 
\begin{lemma} \label{tetrad} Let $\B\subset\bot(S_0)$ be an interval of constant edges that is 
the swollen future of a team with a rascal 
at its left-hand end and a \lpl letter $\r$ at its right-hand end.  
Then at least one of the following holds: 
\begin{enumerate}  
\item[\rm{(i)}] the length of $\B$ is at most $\ll$;
\item[\rm{(ii)}] $\B$ is the swollen future of a team with a left-biased rascal; 
\item[\rm{(iii)}] $\B$ is right-biased;
\item[\rm{(iv)}] neither of the non-constant letters at the ends of $\B$
 is ultimately consumed by an edge of $S_0$. 
\end{enumerate} 
\end{lemma} 

We note here that when the length of $\B$ is at most $\ll$ then we
have a short team, and we have already accounted for short teams.
The following three lemmas correspond to eventualities (ii) to (iv).
  
\begin{lemma} \label{highRascals} 
The sum of the bonuses provided to all teams by left-biased rascals is
less than $(2M + 6M\ttt + 4MT_0 + 2\ll + 6B\ttt + 4BT_0) \n$.
\end{lemma}

\begin{proof} The proof of this result is similar to the work done in
the previous section.  We have a pincer $\Pin_e$ associated to the
rascal $e$.  Since we are only concerned with the times when the
rascal is immediately adjacent to a block of constant letters, it must
be that at time $\tau_1(e) - T_0$ either we are below $\tau_0(e)$ or
$\time(S_{\Pin_e})$ (cf. Definition \ref{PincerDef}). Therefore the
following is an immediate consequence of the Pincer Lemma.
\[	\tau_1(e) - \tau_2(e) \leq \ttt (1 + |\chi(\Pin_e)|) + T_0 .
\]
It now suffices to bound the amount of time for which $e$ is adjacent
to the narrow past of $\B$ before $\tau_2(e)$.  We define $\tau_0'(e)$
to be the latest time when the rascal $e$ has contributed less than
$\ll$ edges to $\bonus(\T)$.  Then the bonus provided by $e$ is at
most $M(\tau_1(e) - \tau_0'(e)) + \ll$. As in the previous section, we
define $\down_2(e)$ to be those edges on the left end of corridors
containing $e$ at times before $\tau_2(e)$ but after $\tau_0'(e)$.
Just as in Lemma \ref{Aget2Lemma} and the corollaries immediately
following it, we then have a notion of {\em depth} of
rascals describing the nesting of the pincers $\Pin_e$\footnote{One
extends the paths $p_e$ and $p_e'$ of Definition 
\ref{RascalPin} back in time to $\partial\Delta$ so as to define the
order definining depth}. We also have
{\em distinguished} rascals (corresponding to the distinguished teams
in Lemma \ref{Aget2Lemma}), and proceeding as in the proof of Lemma
\ref{Aget2Lemma} we get the following estimates:

if $e_1$ is a distinguished rascal of depth $d+1$ and $e_0$ is the
rascal of depth $d$ above it, then the bonus provided by $e_1$ is at
most $2B\Big( T_1 (1 + |\chi(\Pin_{e_0}))|) + T_0 \Big)$, since all of
the bonus provided by $e_1$ must disappear before $\tau_1(e_0)$;

for other rascals $e$ of depth $d+1$ which are below $e_0$ we have a
set of colours $\chi_c(e)$, disjoint for distinct teams such that
\[	|\down_2(e) \cap \down_2(e_0)| \leq T_1(1 + |\chi_c(e)|) +
T_0.	\]
Therefore, summing over the set of rascals which are not distinguished
we get (cf Corollary \ref{Firstdown2sum})
\[ \sum_e |\down_2(e)| \leq 2\Big| \bigcup_e \down_2(e)\Big| + \sum_e \Big
( \ttt (1 + |\chi_c(e)|) + T_0 \Big).	\]
And summing over the same set of rascals, we get 
\[ \sum_e |\down_2(e)| \leq (2 + 3\ttt + 2T_0)\n .	\]
Therefore, for undistinguished rascals, we have
\begin{eqnarray*}
\sum_e \tau_1(e) - \tau_0'(e) & = & \sum(\tau_1(e) - \tau_2(e)) +
\sum(\tau_2(e) - \tau_0'(e))\\
& \leq & (3\ttt + 2T_0)\n + (2 + 3\ttt + 2T_0)\n,
\end{eqnarray*}
and so the contribution of all left-biased rascals is at most
\[	\Big( (2 + 6\ttt + 4T_0)M + 2\ll + 6B\ttt + 4BT_0 \Big) \n,	\]
as required.
\end{proof}

\begin{lemma} \label{RightBiased}  The sum $\sum \life(\B)$ over those
$\B$ that are right-biased but  do not satisfy conditions (i) or (ii)
of Lemma \ref{tetrad} is at most $(3\ttt B + 2T_0B)\n$.
\end{lemma} 

\begin{figure}[htbp] 
\begin{center} 
  
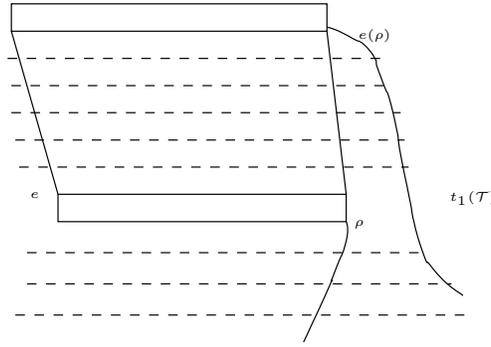 

\caption{A depiction of a right-biased team.} 
\label{r-teams} 
\end{center} 
\end{figure}

\begin{proof}
Once again, as in Lemmas \ref{lowRascals} and \ref{highRascals}, we
obtain compensation for the continuing existence of a non-constant
letter by using the Pincer Lemma to see that colours must be consumed
at a constant rate in order to facilitate the life of $\r$. Thus we
consider the left-fast edge that consumes the pp (i.e. left-most
non-constant) future of $\r$; this edge is denoted $e(\r)$ in Figure
\ref{r-teams}. The Pincer Lemma and the 2 Colour Lemma tell us that if
$\Pin_{e(\r)}$ is the
pincer associated to these paths (with $S_0$ at the bottom) then 
\[      \life(\B) \leq \ttt (1 + |\chi(\Pin_{e(\r)})|) + T_0.   \]
Suppose that $\B$ and $\B'$ are two right-biased blocks with
associated edges $e(\r)$ and $e(\r')$ consuming their reapers. We claim
that the sets $\chi(\Pin_{e(\r)})$ and $\chi(\Pin_{e(\r')})$ are
disjoint.  The key point to observe is that since we are not in case
(ii) of Lemma \ref{tetrad} the length of the swollen future of $\B$
increases from $\time(S_0)$ to the top of $\Pin_{e(\r)}$; since $\B$
had length at least $\ll$, we therefore have a block of more than
$\ll$ of more than $\ll$ constant edges at the top of
$\Pin_{e(\r)}$. Thus the pincers associated to $\B$ and $\B'$ are
either disjoint or nested. Hence $\chi(\Pin_{e(\r)})$ and
$\chi(\Pin_{e(\r')})$ are disjoint. Thus summing over all right-biased
blocks $\B$ we obtain
\[      \sum_{\B \mbox{ right-biased}} \life(\B) \leq (3\ttt B + 2T_0B)
\n, \]
as required.
\end{proof} 
Since any letter consumes less than $M$ constant letters in any unit of time, we conclude: 
\begin{corollary} \label{rightCor} The sum  of the effective volumes of  
all blocks that are right-biased but  do not satisfy conditions (i) and (ii) 
of Lemma \ref{tetrad} is at most $(3M\ttt B + 2MT_0B)\n$.  
\end{corollary} 
 
\begin{lemma} The sum of all blocks that satisfy condition (iv) of  
Lemma \ref{tetrad} is at most $(2B +1)\n$. 
\end{lemma} 
 
\begin{proof} Possibility (iv) involves several subcases: the key
event which halts the growth of the swollen future of $\B$ may be a
collision with $\partial\Delta$ or a singularity;  it may also be that
the key event is that the future of the rascal or reaper adjacent to
$\B$ is cancelled by an edge that is not in the future of $S_0$.  
 
But no matter what these key events may be, since we are in not in
cases (ii) or (iii), associated to the blocks in case (iv) we have the
following set of paths partitioning that part of the diagram $\Delta$
bounded by $S_0$ and the arc of $\partial\Delta$ connecting the
termini of the edges at the ends of $S_0$:

The path $\pi_l$ begins at $\time (S_0)$ and follows the pp-future of
the rascal at the right-end of the future of $\B$ until it hits the
boundary, a singularity, or else is cancelled by an edge $\e_l$ not in
the future of $S_0$; if it hits the boundary, it ends; if it hits a
singularity, $\pi_l$ crosses to the bottom of the corridor $S$ on the
other side of the singularity, and turns left to follow $\bot(S)$ to
the boundary (see Figure \ref{Pi_lOne}); if $\e_l$ cancels with the
pp-future of the rascal, then $\pi_l$ follows the past of $\e_l$
backwards in time to the boundary (see Figure \ref{Pi_lTwo}).

\begin{figure}[htbp] 
\begin{center} 
  
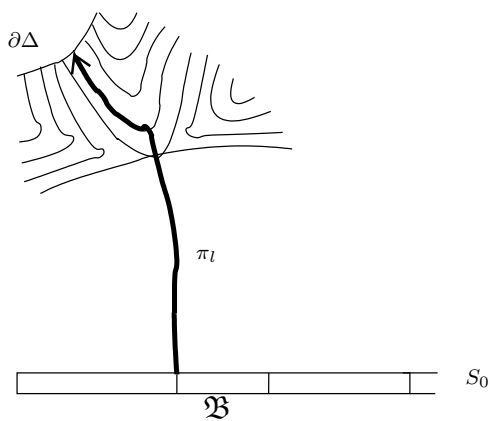 
  
\caption{The path $\pi_l$ hits a singularity.} 
\label{Pi_lOne} 
\end{center} 
\end{figure} 

\begin{figure}[htbp]
\begin{center}

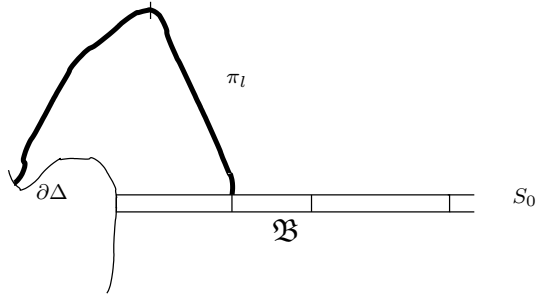

\caption{The path $\pi_l$ in cancelled from outside of the future of
$S_0$}
\label{Pi_lTwo}
\end{center}
\end{figure}

The path $\pi_r$ describing the fate of $\r$ is defined similarly
(except that it turns right if it hits a singularity). 
 
It is clear from the construction that no two of these paths can cross, thus we have the  partition  
represented schematically in Figure \ref{partition}. 
 
\begin{figure}[htbp]
\begin{center}

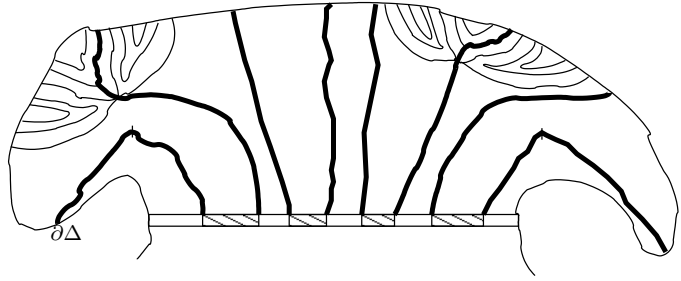

\caption{The schematic partition of $\Delta$ by the paths $\pi_l$ and
$\pi_r$.}
\label{partition}
\end{center}
\end{figure}

\def\bdy{\text{\rm{bdy}}} 
\def\up{\text{\rm{up}}} 
Given a swollen team $\B$ of type (iv), we follow the swollen future of $\B$ until its flow is 
interrupted (at time $\iota(\B)$, say) 
by meeting a singularity, the boundary of $\Delta$, or else its rascal or reaper is cancelled. 
Consider the set of corridors that contain some component of the swollen future of $\B$ 
after $\iota(\B)$. Consider also the set of edges $\bdy(\B) \subseteq
\partial\Delta$ 
that lie in the swollen future of $\B$. We keep account of the set of corridors by recording 
the set of their ends on $\partial\Delta$, except that we ignore an end if we have to cross 
a path $\pi_l$ or $\pi_r$ to reach it. Note that at least one end of each corridor is recorded. 
Let $\up(\B)\subset\partial\Delta$ denote the  set  of ends recorded.  
  
Since the sets $\bdy(\B)$ and $\up(\B)$ are contained in the portion of $\partial\Delta$ accorded 
to $\B$ by the partition formed by the paths $\pi_l$ and $\pi_r$, the sets associated to different 
$\B$ are disjoint. In each unit of time beyond $\iota(\B)$ each
component of the swollen future of $\B$ can shrink by at most $2B$ (by
Lemma \ref{BCL}).  The set $\up(\B)$ measures the sum of the number of
components over all such times, and $|\bdy(\B)|$ is the number of
uncancelled edges. Thus we see that the length of
the swollen future of $\B$ at time $\iota(\B)$ is at most $2B|\up(\B)|
+ |\bdy(\B)|$. Finally, the continued presence of the rascal ensures
that the swollen future of $\B$ grows in each interval of time from
$\time(S_0)$ to $\iota(\B)$.  Thus it follows that the
length of $\B$ is also bounded by this number. So summing over all
$\B$ of type (iv) we have: 
$$ 
\sum |\B| \,\le\, \sum \Big(2B|\up(\B)| + |\bdy(\B)|\Big) \le (2B +1)\n , 
$$ 
as required. 
\end{proof} 
 Summarising the results of this section we have 
 
\bl \label{BonusBound} Summing over all teams that are not short, we
have 
\[      \sum_{\T}|\bonus(\T)| \leq \Big( \Bb \Big) \n .  \]
\end{lemma}

\section{The proof of the Main Theorem} \label{summary} 
 
Pulling all of the previous results together, define  
$$
K_1 = \AFourC,
$$
and
$$ 
K = \K. 
$$ 
 
\begin{theorem} \label{S0<=Kn}
$|S_0| \leq K\n$. 
 
\Prf: 
The corridor $S_0$ can be subdivided into distinct colours which form connected regions. 
  Each colour $\mu$ can be partitioned into connected (possibly empty) regions $A_1(S_0,\mu),  
  A_2(S_0,\mu), A_3(S_0,\mu), A_4(S_0,\mu)$ and $A_5(S_0,\mu)$.  By Lemma \ref{A1A5Lemma}, 
  Proposition \ref{SummaryLemma}, Lemma \ref{A3Lemma}, Proposition \ref{A2Prop} and Lemma \ref{A1A5Lemma}, respectively, 
\begin{eqnarray*} 
\sum_{\mu \in S_0}|A_1(S_0,\mu)|  & \leq & C_0\n ,\\ 
\sum_{\mu \in S_0}|A_2(S_0,\mu)|  & \leq &  K_1  \n ,\\ 
\sum_{\mu \in S_0}|A_3(S_0,\mu)|  & \leq & (2B+1)\n ,\\ 
\sum_{\mu \in S_0}|A_4(S_0,\mu)|  & \leq &  K_1 \n,\mbox{ and} \\ 
\sum_{\mu \in S_0}|A_5(S_0,\mu)|  & \leq & C_0\n . 
\end{eqnarray*} 
Summing completes the proof of Theorem \ref{S0<=Kn}.
\et
Since there are at most $\frac{\n}{2}$ corridors in $\Delta$,  
\[      \mbox{Area}(\Delta) \leq \frac{K}{2}\n^2 ,       \] 
which proves the Main Theorem.

\section{Glossary of Constants}

$B$ -- the Bounded Cancellation constant (Lemmas \ref{BCL} 
and \ref{SingularityProp}). 
 
$C_0$ -- maximum distance a left-fast (right-fast)  
letter can be from the left (right) edge of its colour if it is to be
cancelled from the
left (right) within the future of the corridor.  See Lemma \ref{A1A5Lemma}. 
 
$C_1$ -- an upper bound on the   
lengths of the subintervals 
 $C_{(\mu,\mu')}(1)$ of $A_4(S_0,\mu)$. By definition, $C_{(\mu,\mu')}(1)$ 
 is consumed by $\mu'(S_0)$; it begins at the right end of
$A_4(S_0,\mu)$ and ends   at the last non-constant letter. See Lemma
\ref{C1Lemma}. Note that one can take $C_1=2mB^2$.

$M$ -- the maximum of the lengths of the images  $\phi(a_i)$ of the basis 
elements $a_i$, i.e. the maximum length of $u_1, \ldots, u_m$ in
the presentation $\P$ (see equation \ref{presentation}).

$M_{inv}$ -- the maximum of the lengths of $\phi^{-1}(a_i)$. 
 
$T_0$ -- the constant from the 2-Colour Lemma (Lemma
\ref{TwoColourLemma}). For all positive words $U$ and $V$, if
$U$ neuters  $V^{-1}$ then it does so in at most $T_0$ steps. 

$\hat{T_1}$ -- the constant from the Unnested Pincer Lemma, Theorem \ref{prePincerLemma}.

$T_1'$ -- the constant from Definition \ref{T1'Lemma}.  Recall that we stipulate that $T_1' \ge \hat{T_1}$.

$T_1 := T_1' + 2T_0$ -- $T_1$ is the constant from the Pincer Lemma, Theorem \ref{PincerLemma}. 

$C_4 := MM_{inv}$

$\ll := {\rm max} \{ 2B(T_0 + 1)+1, MC_4 \}$

Finally, $K_1$ is defined to be
\[	\AFourC	,	\]
and $K = 2C_0 + 2K_1 + 2B + 1$.

\end{document}